\documentclass[preprint,12pt]{elsarticle}

\RequirePackage{amsmath, amsthm, amsfonts, amssymb}
\RequirePackage{graphicx, color}
\RequirePackage{bm}
\RequirePackage{pifont}
\RequirePackage[english]{babel}
\RequirePackage[per-mode=symbol]{siunitx}
\usepackage{tikz}
\usepackage{subcaption}

\graphicspath{{./Parts/Images/}}

\newcommand{\labelstyle}[1]{\textsc{#1}}

\DeclareMathOperator{\tr}{{\rm tr}}
\newcommand{\jumpS}[1]
    {\mathop{ \left[\!\!\left[#1\right]\!\!\right]_{\Gamma} }\nolimits}

\newtheorem{problem}{Problem}

\newcommand{\OPTXFEM}{\labelstyle{opt-xfem}}
\newcommand{\OPTFEM}{\labelstyle{opt-fem}}

\newcommand{\TPFA}{\labelstyle{tpfa}{}}
\newcommand{\FEM}{\labelstyle{fem}{}}
\newcommand{\MFEM}{\labelstyle{mfem}{}}
\newcommand{\MVEM}{\labelstyle{mvem}{}}

\newcommand{\UPWIND}{\labelstyle{upwind}{}}
\newcommand{\SUPG}{\labelstyle{supg}{}}

\newcommand{\TPFAUP}{\labelstyle{tpfaup}}
\newcommand{\MFEMUP}{\labelstyle{mfemup}}
\newcommand{\MVEMUP}{\labelstyle{mvemup}}
\newcommand{\MFEMSUPG}{\labelstyle{mfemsupg}}
\newcommand{\FEMSUPG}{\labelstyle{femsupg}}
\newcommand{\XFEMSUPG}{\labelstyle{xfemsupg}}

\begin{document}


\journal{Advances in Water Resources}


\begin{frontmatter}

\title{Analysis of conforming, non-matching, and polygonal methods for Darcy and advection-diffusion-reaction simulations in
discrete fracture networks}

\author[polito,indam]{Andrea Borio}
\author[polito,indam]{Alessio Fumagalli}
\author[polito,indam]{Stefano Scial{\`o}}

\address[polito]{Dipartimento di Scienze Matematiche, Politecnico di Torino,
Corso Duca degli Abruzzi 24, 10129 Torino, Italy}
\address[indam]{Member of the INdAM research group GNCS}


\begin{abstract}

  The aim of this study is to compare numerical methods for the simulation of
  single-phase flow and transport in fractured media, described here by means of the Discrete Fracture Network (DFN) model. A Darcy problem is solved to
  compute the advective field, then used in a subsequent time dependent transport-diffusion-reaction problem. The numerical schemes are benchmarked in terms of flexibility in handling geometrical complexity, mass conservation and stability issues for advection dominated flow regimes. To this end, three benchmark cases have been specifically designed and are here proposed, representing some of the most critical issues encountered in DFN simulations.

\end{abstract}




\end{frontmatter}


\section{Introduction}

The movement of liquids in the underground is heavily influenced by
the presence of fractures and their relative intersections
\cite{Dowd2009,Hardebol2015,Flemisch2016a}. Fractures are
discontinuities (here assumed planar) along which a rock has been
broken, mainly due to geological movements or to artificial
stimulation \cite{McClure2015}.  In this work, we are considering only
open structures, characterized by a geometrical aperture, that allow a
liquid to flow through \cite{Erhel2009}.  Possibilities are actual
fractures, faults, and joints. We are thus excluding low permeable
(closed/impervious) objects such as veins or dykes.  For particular
underground compositions (e.g., granite, shale, or sandstone) the rock
permeability is several orders of magnitude smaller than the fracture
permeability. It is a common choice and a reasonable approximation to
ignore the rock matrix effect in the simulations and rely only on the
fractures. The framework is the discrete fracture network model (DFN),
where the aperture is not a geometrical constraint but a parameter in
the bidimensional representation of the fractures by reduced models,
see \cite{Martin2005} and the forthcoming references.

The geometrical complexity of natural fracture networks may impose
difficulties in the numerical simulations, due to the presence of
small intersections between fractures, different intersection
configurations (e.g., Y-type or L-type), small angles of intersection,
and small distances between intersection lines. In the literature,
three main approaches are developed to overcome these difficulties.

The first approach considers rather standard numerical scheme for the
discretization of the physical equations, and relies on robust grid
generation software. The coupling conditions among fractures are
imposed via Lagrange multipliers on a representation of the interfaces
conforming with the computational mesh on all sides. To be specific,
the edges of the computational elements have to match at the
intersections.  This approach, commonly known as conforming
discretization, may suffer, for example, when intersections lines are
very close to each other, see
\cite{Hyman2014a,Hyman2015a,Mustapha2007a,NGO2017}.  To partially
overcome these problems, a possibility is to consider the class of
virtual element methods that allows grid elements of general
shapes. Two different numerical methods of this type are introduced in
\cite{Benedetto2014,Fumagalli2016a,Benedetto2017,Fumagalli2017a,BBorth}.

A second possibility is to keep the intersections explicitly
represented, but relaxing the conformity of the edges. This approach,
named non-conforming, requires more advanced numerical schemes based
on the mortar technique. In this case, we relax the actual generation
of the fracture meshes which usually gives less discrete elements than
the conforming approach \cite{Benedetto2016,Nordbotten2018}. This
method may still suffer in presence of severe geometries. Also in this
case the virtual element methods are an interesting option to further
decrease the computational cost, see
\cite{BeiraodaVeiga2014b,BeiraoVeiga2016}.

A third family of schemes comprises the so-called non-matching
discretizations. In terms of mesh generation, in this case the
intersections do not place any constraints, as the fractures are
meshed independently and the coupling conditions are imposed by an
optimization procedure. A functional that measures the mismatch in the
coupling conditions is minimized iteratively, where only the degrees
of freedom involved in the intersections (the cut region) are
considered \cite{BPSc,BPSd,BSV,BBBPSsupg,BPSf}. This procedure is
independent on the actual numerical schemes, might take advantage of
ad-hoc strategies to enrich the solution in the cut region.  Extended
finite elements are a successful example, see, for example,
\cite{DAngelo2011,Formaggia2012,Fumagalli2012d,Fumagalli2012bb}.

The present work extends and enriches the one proposed in
\cite{Scialo2017} to more complex physical phenomena. The concepts
previously discussed are applied to Darcy and heat
transport-and-diffusion models. The Darcy velocity is computed first
and then used as an advective field for the heat equation. High
quality computation or reconstruction of the Darcy velocity may
significantly impact the temperature distribution. Moreover, in the
heat equation the transport part usually dominates the diffusion
(P\'eclet greater than 1) and stable or stabilized schemes are needed
to avoid or limit spurious oscillations that might compromise the
accuracy of the solution. In the numerical tests, we are considering
several numerical schemes for the comparison to cover most of the
combinations discussed and try to assess their performances. The aim
of this work is thus twofold: establish a set of benchmark cases and
give guidance in the development of more advanced numerical schemes to
solve this problem.

The paper is organized as follow. In Section \ref{sec:model} the Darcy
and heat models are introduced and discussed, with particular focus on
the coupling conditions. Section \ref{sec:discretization} is devoted
to the description of the proposed numerical schemes. Three numerical
examples are presented and discussed in Section \ref{sec:examples},
comparing the performances of the considered numerical
schemes. Finally, in Section \ref{sec:conclusion} we draw some
conclusions and suggestions for future developments.



\section{Mathematical model}
\label{sec:model}

In this section we introduce the mathematical model used to describe
the hydraulic head and Darcy velocity profiles in a discrete fracture
network. Once this problem is solved, the Darcy velocity is considered
as advective field to simulate the transport and diffusion of heat in
 DFNs.

Fractures are considered as non-overlapping planar objects, which can
be connected to other fractures through intersection segments, also called traces.  We consider
$N_\Omega$ fractures $\Omega_i \subset \mathbb{R}^3$ with boundary
$\partial \Omega_i$, which compose the discrete fracture network
$\Omega = \cup_{i=1}^{N_\Omega} \Omega_i$, and we denote
its boundary as $\partial \Omega_{\rm}$ with outward unit normal
$\bm{n}_{\rm ext}$, defined on each fracture plane as the unit vector
normal to the fracture boundary pointing outward from the fracture
polygon.  To keep the presentation clear, we make a distinction
between external boundaries, where we will impose data, and internal
interfaces, i.e. the traces, where we will couple the mathematical models on the
fractures.

\begin{figure}
  \centering 
\begingroup%
  \makeatletter%
  \providecommand\color[2][]{%
    \errmessage{(Inkscape) Color is used for the text in Inkscape, but the package 'color.sty' is not loaded}%
    \renewcommand\color[2][]{}%
  }%
  \providecommand\transparent[1]{%
    \errmessage{(Inkscape) Transparency is used (non-zero) for the text in Inkscape, but the package 'transparent.sty' is not loaded}%
    \renewcommand\transparent[1]{}%
  }%
  \providecommand\rotatebox[2]{#2}%
  \newcommand*\fsize{\dimexpr\f@size pt\relax}%
  \newcommand*\lineheight[1]{\fontsize{\fsize}{#1\fsize}\selectfont}%
  \ifx\svgwidth\undefined%
    \setlength{\unitlength}{158.96681538bp}%
    \ifx\svgscale\undefined%
      \relax%
    \else%
      \setlength{\unitlength}{\unitlength * \real{\svgscale}}%
    \fi%
  \else%
    \setlength{\unitlength}{\svgwidth}%
  \fi%
  \global\let\svgwidth\undefined%
  \global\let\svgscale\undefined%
  \makeatother%
  \begin{picture}(1,0.66288874)%
    \lineheight{1}%
    \setlength\tabcolsep{0pt}%
    \put(0,0){\includegraphics[width=\unitlength,page=1]{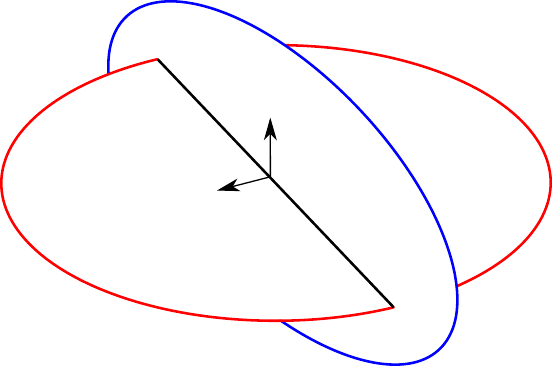}}%
    \put(0.05782312,0.32252549){\color[rgb]{0,0,0}\makebox(0,0)[lt]{\lineheight{0}\smash{\begin{tabular}[t]{l}$\Omega_i$\end{tabular}}}}%
    \put(0.72400433,0.05569427){\color[rgb]{0,0,0}\makebox(0,0)[lt]{\lineheight{0}\smash{\begin{tabular}[t]{l}$\Omega_j$\end{tabular}}}}%
    \put(0.33114717,0.51657863){\color[rgb]{0,0,0}\makebox(0,0)[lt]{\lineheight{0}\smash{\begin{tabular}[t]{l}$\Gamma_k$\end{tabular}}}}%
    \put(0.41769339,0.272502){\color[rgb]{0,0,0}\makebox(0,0)[lt]{\lineheight{0}\smash{\begin{tabular}[t]{l}$\bm{n}_{i,-}$\end{tabular}}}}%
    \put(0.5084644,0.39153042){\color[rgb]{0,0,0}\makebox(0,0)[lt]{\lineheight{0}\smash{\begin{tabular}[t]{l}$\bm{n}_{j,+}$\end{tabular}}}}%
  \end{picture}%
\endgroup%

  \caption{Representation of two fractures $\Omega_i$ and $\Omega_j$
    intersecting in $\Gamma_k$.}%
  \label{fig:dfn}
\end{figure}
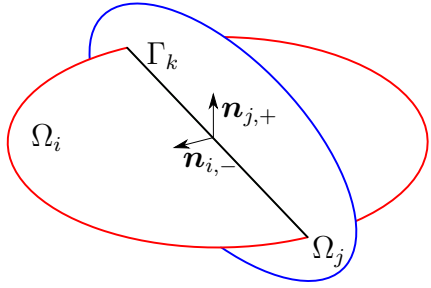

Given two distinct and intersecting fractures $\Omega_i$ and
$\Omega_j$, with $i \neq j$, we indicate their intersection (trace) as
$\overline{\Gamma_k} = \overline{\Omega_i} \cap
\overline{\Omega_j}$. For simplicity, we assume that a trace is formed
only by two distinct fractures, however this assumption can be
relaxed. A natural order of indexes can be introduced to numerate the
traces $\Gamma_k$ from 1 to $N_\Gamma$, being the latter their
cardinality. We consider also the function
$t: \{1, \ldots, N_\Omega\} \times \{1, \ldots, N_\Omega\} \rightarrow
\{1, \ldots, N_\Gamma\}$ such that $k = t(i, j)$ with
$\overline{\Gamma_k} = \overline{\Omega_i} \cap
\overline{\Omega_j}$. We have $t(i, j) = t(j, i)$ and its inverse
$t^{-1}$ is well defined such that $(i, j) = t^{-1}(k)$ where
$i<j$. We indicate with $\Gamma = \cup_{k=1}^{N_\Gamma} \Gamma_k$ the
union of all the traces and by $\Gamma_{\Omega_i}$ the set of traces
belonging to the fracture $\Omega_{i}$. Moreover, consider a fracture
$\Omega_i$ and a trace $\Gamma_k$, with
$k\in\Gamma_{\Omega_i}$. $\Gamma_k$ naturally subdivides $\Omega_i$ in
two sub-regions, indicated by $\Omega_{i,+}^{k}$ and
$\Omega_{i,-}^{k}$, such that
$\Gamma_k\subset(\partial \Omega_{i,+}^k\cap \Omega_{i,-}^k)$. To each
of these sub-regions we associate an outward unit normal perpendicular
to $\Gamma_k$, denoted by $\bm{n}_{i,+}$ and $\bm{n}_{i,-}$ (with
$\bm{n}_{i,+} = -\bm{n}_{i,-}$) and a trace operator, respectively
denoted by $\tr_{i,+}^k$ and $\tr_{i,-}^k$.  Generic trace operators
on fracture $\Omega_i$ are denoted by $\tr_i$. An example of a simple DFN with the
introduced nomenclature is given in Figure \ref{fig:dfn}.

We present the Darcy problem in
strong and weak form in Subsection \ref{subsec:darcy}, to compute the hydraulic head and Darcy
velocity, whereas the heat equation is
introduced along with its weak formulation and functional setting in Subsection \ref{subsec:heat}.

\subsection{The Darcy model}\label{subsec:darcy}

This section is devoted to the presentation of the mathematical models used to describe the
hydraulic head field $h$ and Darcy velocity $\bm{u}$ in a discrete fracture
network, further details being available in \cite{BPSc,BPSd,Fumagalli2016a,Scialo2017} and references therein.
Unknowns and parameters restricted to a fracture
$\Omega_i$ are denoted by a subscript $i$.

For clarity in the exposition, we start considering a single fracture
$\Omega_i$.  The Darcy model on $\Omega_i$ reads: find $(\bm{u}_i, h_i)$ such that
\begin{subequations}\label{eq:darcy_strong}
  \begin{gather}\label{eq:darcy_strong_eq}
    \begin{aligned}
      &\bm{u}_i + K_i \nabla h_i = \bm{0}\\
      &\nabla \cdot \bm{u}_i = f_i
    \end{aligned}
    \quad \text{in } \Omega_i \setminus \Gamma
  \end{gather}
  Variables, data and differential operators are defined on the
  tangent plane of the fracture. $K_i$ is the hydraulic conductivity
  tensor, which is symmetric and positive definite, and $f_i$ is a
  scalar source/sink term. In our applications, following lubrication
  theory \cite{Bear1972,Wangen2010,Schon2011}, we consider a
  scalar permeability, obtained by the cubic law:
  $k_i = \epsilon_i^2/12$,
  where $\epsilon_i$ is the fracture aperture. Moreover, the hydraulic
  conductivity in \eqref{eq:darcy_strong_eq} is isotropic and defined as
  \begin{gather*}
    K_i = \frac{\epsilon_i k_i \rho_w g}{\mu} I\,,
  \end{gather*}
  with $\rho_w$ the fluid density, $g$ the gravity
  acceleration, $\mu$ the dynamic viscosity of the fluid.

  The boundary conditions on $\Omega_i$ are
  \begin{gather}\label{eq:darcy_strong_bc}
    \begin{aligned}
      & \tr_i^\partial h_i = \overline{h}_i && \text{on } \partial \Omega_{i, D}\\
      & \tr_i^\partial \bm{u}_i \cdot \bm{n}_{\rm ext} =0 && \text{on }
      \partial \Omega_{i, N}
    \end{aligned}
  \end{gather}
  where $\partial \Omega_{i, D}$ and $\partial \Omega_{i, N}$ are
  disjoint portions of the boundary of $\Omega_i$ such that
  $\overline{\partial \Omega_i} = \overline{\partial \Omega_{i, D}}
  \cup \overline{\partial \Omega_{i, N}}$, and $\tr_i^\partial$ gives the trace on the boundary of $\Omega_i$.
  The data $\overline{h}_i$ is the given pressure head at the boundary $\partial \Omega_{i, D}$
  and we model $\partial \Omega_{i, N}$ as impervious. Other boundary
  conditions are possible, however, in order to keep the presentation as simple as possible and
  coherent with the examples proposed in Section \ref{sec:examples}, we
  consider only these, being the generalization to other boundary conditions straightforward.
  We assume, for a single fracture, that $\partial \Omega_{i, D}$ is not empty; in presence of intersecting
  fractures, instead, we can have $\partial \Omega_{i, D} = \emptyset$ for all
  the fractures but one.  Equations \eqref{eq:darcy_strong_eq} and
  \eqref{eq:darcy_strong_bc} are well studied in the literature, refer
  to the aforementioned references.

  Let us now consider two distinct fractures $\Omega_i$ and $\Omega_j$
  forming an intersection $\Gamma_k$. On both fractures we apply
  Equations \eqref{eq:darcy_strong_eq} and \eqref{eq:darcy_strong_bc}
  and we assume continuous coupling conditions for the hydraulic head
  and the normal component of the Darcy velocity at the trace. The
  coupling conditions between two fractures $\Omega_i$ and $\Omega_j$
  such that $t(i,j) = k$, read
  \begin{gather}
    \label{eq:darcy_strong_coupling}
    \begin{aligned}
      &\sum_{l \in \{i, j\}} \tr_{l,+}^k \bm{u}_l \cdot \bm{n}_{l,+} +
      \tr_{l,-}^k \bm{u}_l \cdot
      \bm{n}_{l,-} = 0\\
      &\tr_{i,+}^k h_i = \tr_{i,-}^k h_i = \tr_{j,+}^k h_j =
      \tr_{j,-}^k h_j
    \end{aligned}
    \quad \text{on $\Gamma_k$}\,.
  \end{gather}
  The case of multiple fractures follows immediately from our
  assumptions.
\end{subequations}
The Darcy equation can be formalised in the following problem.
\begin{problem}[Darcy equation on a DFN - mixed
  formulation]\label{pb:model_darcy}
  Given the set of fractures $\Omega$ and traces $\Gamma$, find
  $(\bm{u}, h)$ such that
  \eqref{eq:darcy_strong_eq}-\eqref{eq:darcy_strong_bc}-\eqref{eq:darcy_strong_coupling}
  are satisfied.
\end{problem}
The weak formulation of Problem \ref{pb:model_darcy} requires the
introduction of the functional spaces on each $\Omega_i$, namely
\begin{gather*}
  {V}_i = H_{\nabla \cdot}(\Omega_i) = \left\{ \bm{v} \in
    \left[L^2(\Omega_i)\right]^2:\, \nabla \cdot \bm{v} \in L^2(\Omega_i)\right \}
  \quad \text{and} \quad Q_i = L^2(\Omega_i)\,,
\end{gather*}
equipped with the natural norms that make them complete. The spaces
for the global network are
\begin{equation*}
  \begin{split}
    V &= \left\{ \bm{v}\colon \bm{v}_i \in V_i \,,\, \tr_i^\partial \bm{v}_i
      \cdot \bm{n}_{\rm ext} =0 \text{ on $\partial \Omega_{i,N}$}
      \quad \forall i\,, \vphantom{\sum_{l \in \{i, j\}}}\right.
    \\
    &\quad \left.\sum_{l \in \{i, j\}} \tr_{l,+}^k \bm{u}_l \cdot
      \bm{n}_{l,+} + \tr_{l,-}^k \bm{u}_l \cdot \bm{n}_{l,-} = 0 \quad
      \forall k \,,\, t^{-1}(k) = (i,j) \right\} \,,
  \end{split}
\end{equation*}
and $Q = \bigoplus_i Q_i$, both endowed with their natural norms. The following linear forms are then introduced:
$a_i : V_i \times V_i \rightarrow \mathbb{R}$,
$b_i: V_i \times Q_i \rightarrow \mathbb{R}$,
$G_i: \left(H^{\frac{1}{2}}(\partial \Omega_{i, D})\right)^\prime
\rightarrow \mathbb{R}$ and $F_i: Q_i \rightarrow \mathbb{R}$ such
that, $\forall \bm{u},\bm{v} \in V_i$ and $\forall q \in Q_i$
\begin{align*}
  a_i(\bm{u}, \bm{v}) &= \left(K_i^{-\frac{1}{2}} \bm{u},
                        K_i^{-\frac{1}{2}} \bm{v}\right)_{\Omega_i}\,,
  &
    b_i(\bm{u}, q) &= - (\nabla \cdot \bm{u}, q)_{\Omega_i}\,,
  \\
  G_i(\bm{v}) &= -\langle \overline{h}_i, \tr_i \bm{v} \cdot \bm{n}
                \rangle_{\partial \Omega_{i,D}}\,, & F_i(q) &= - (f, q)_{\Omega_i} \,.
\end{align*}
We also require that $f \in Q_i$,
$\overline{h}_i \in H^{\frac{1}{2}}(\partial \Omega_{i, D})$, and
$K_i \in L^\infty(\Omega_i)$.  The less standard functional spaces
associated with $G_i$ are due to the low regularity of a
$H_{\nabla \cdot}$-function at a boundary, as well as the duality
pairing used instead of a scalar product.  The global forms
$a : V \times V \rightarrow \mathbb{R}$,
$b: V \times Q \rightarrow \mathbb{R}$,
$G: H^* \rightarrow \mathbb{R}$, and $F: Q \rightarrow \mathbb{R}$ are
the sums of the respective local ones, where we denoted by
$H^{*} = \bigoplus_{i} \left(H^{\frac{1}{2}}(\partial \Omega_{i,
    D})\right)^\prime$, the space associated to the boundary term.
\begin{problem}[Darcy equation on a DFN - mixed weak
  formulation]\label{pb:model_darcy_weak}
  The weak formulation of Problem \ref{pb:model_darcy} is: find
  $( \bm{u}, h) \in V \times Q$ such that
  \begin{align*}
    &a(\bm{u}, \bm{v}) + b(\bm{v}, h) = G(\bm{v}) && \forall \bm{v} \in V\\
    &b(\bm{u}, q) = F(q) &&  \forall q \in Q
  \end{align*}
\end{problem}

Some numerical formulations require the problem written in term of
pressure head alone, thus we consider also the primal formulation of
Problem \ref{pb:model_darcy}. Given a fracture $\Omega_i$ the model
is: find $h_i$ such that
\begin{subequations}\label{eq:darcy_strong_primal}
  \begin{align}\label{eq:darcy_strong_primal_eq}
    \begin{aligned}
      &- \nabla \cdot \left(K_i \nabla h_i\right) = f_i && \text{in }
      \Omega_i \setminus
      \Gamma\\
      & \tr_i h_i = \overline{h}_i && \text{on } \partial \Omega_{i, D}\\
      & \tr_i K_i \nabla {h}_i \cdot \bm{n}_{\rm ext} =0 && \text{on }
      \partial \Omega_{i, N}
    \end{aligned}
  \end{align}
  Given a second distinct fracture $\Omega_j$ which intersects
  $\Omega_i$ in $\Gamma_k$, the following coupling conditions are added:
  \begin{gather}
    \label{eq:darcy_strong_primal_coupling}
    \begin{aligned}
      &\sum_{l \in \{i, j\}} \tr_{l,+} K_l \nabla h_l \cdot
      \bm{n}_{l,+} + \tr_{l,-} K_l \nabla h_l \cdot
      \bm{n}_{l,-} = 0\\
      &\tr_{i,+}^k h_i = \tr_{i,-}^k h_i = \tr_{j,+}^k h_j =
      \tr_{j,-}^k h_j
    \end{aligned}
    \quad \text{on } \Gamma_k.
  \end{gather}
\end{subequations}

Once the pressure head is computed, the Darcy velocity is given
by the first equation of \eqref{eq:darcy_strong_eq}.  We can
formalize also the primal formulation of the Darcy equation in the
following problem.
\begin{problem}[Darcy equation on a DFN - primal
  formulation]\label{pb:model_darcy_primal}
  Given the set of fractures $\Omega$ and traces $\Gamma$, find $h$
  such that
  \eqref{eq:darcy_strong_primal_eq}-\eqref{eq:darcy_strong_primal_coupling}
  are satisfied $\forall i,k$.
\end{problem}
In this case, the weak formulation is simpler than the previous
one. To keep the notation standard and since we are introducing a
separate formulation of the Darcy problem, in the following we commit
abuses in notation. We introduce the functional spaces
\begin{gather*}
  V_i = \{ v \in H^1(\Omega_i):\, \tr_i^\partial v = \overline{h}_i \text{ on }
  \partial \Omega_{i, D} \} \,,
\end{gather*}
and $V_{i,0}$ in the special case of homogeneous head condition. The linear
forms now read: $a_i: V_i \times V_{i,0} \rightarrow \mathbb{R}$ and
$F: V_{i, 0} \rightarrow \mathbb{R}$ such that,
$\forall v \in V_{i,0}$,
\begin{gather*}
  a_i(h, v) = (K_i^{\frac{1}{2}} \nabla h, K_i^{\frac{1}{2}} \nabla
  v)_{\Omega_i} \quad\text{and}\quad F_i(v) = (f, v)_{\Omega_i} \,.
\end{gather*}
In this case, we ask $f \in \left(V_{i,0}\right)^\prime$,
$\overline{h}_i \in H^{\frac{1}{2}}(\partial \Omega_{i, D})$, and
$K_i \in L^\infty(\Omega_i)$.  The functional spaces on the whole network are
\begin{gather*}
  \begin{split}
    V = \left\{v\colon v_i \in V_i \quad\forall i \,,\quad \tr_{i,+}^k
      v_i = \tr_{i,-}^k v_i = \tr_{j,+}^k v_j = \tr_{j,-}^k v_j
      \quad\right.&
    \\
    \quad\left.\forall k\,,\, t^{-1}(k) = (i,j) \right\}&\,,
  \end{split}
  \\
  \begin{split}
    V_0 = \left\{v\colon v_i \in V_{i,0} \quad\forall i \,,\quad
      \tr_{i,+}^k v_i = \tr_{i,-}^k v_i = \tr_{j,+}^k v_j =
      \tr_{j,-}^k v_j \quad\right.&
    \\
    \quad\left.\forall k\,,\, t^{-1}(k) = (i,j) \right\}&\,,
  \end{split}
\end{gather*}
with their natural norms.  The global forms
$a : V \times V_0 \rightarrow \mathbb{R}$ and
$F: V_0 \rightarrow \mathbb{R}$ are the sums of the respective local
ones.
\begin{problem}[Darcy equation on a DFN - primal weak
  formulation]\label{pb:model_darcy_primal_weak}
  The weak formulation of Problem \ref{pb:model_darcy_primal} is: find
  $h \in V$ such that
  \begin{align*}
    a(h, v) = F(v) \qquad \forall v \in V_0 \,.
  \end{align*}
\end{problem}
We note that to obtain symmetric forms in Problem
\ref{pb:model_darcy_primal_weak} a lifting technique of the Dirichlet datum can be
considered.

\subsection{The heat equation}\label{subsec:heat}

In this part, we present the heat equation on the DFN. Once Problem
\ref{pb:model_darcy} is solved, the Darcy velocity can be used as
advective field in the transport problem. We denote the temperature in a DFN as
$\theta$, and its restriction to fracture $\Omega_i$ as $\theta_i$.  The heat equation on $\Omega_i$ reads: given
$\bm{u}_i$ find $\theta_i$ such that
\begin{subequations}\label{eq:heat_strong}
  \begin{gather}\label{eq:heat_strong_eq}
    \zeta_i \partial_t \theta_i + \nabla \cdot (\bm{u}_i \theta_i -
    D_i \nabla \theta_i) + \iota_i (\theta_i - \hat{\theta}_i) = 0 \quad
    \text{in } \Omega_i \setminus \Gamma \times (0, T] \,,
  \end{gather}
  where $T\in\mathbb{R}$ is the end time of the simulation. Also in
  this case, the variables, data, and differential operators are
  defined on the tangent plane of the fracture. The relations
  to compute the physical parameters are
  \cite{Bear1972,Wangen2010,Schon2011}:
  \begin{align*}
    \zeta_i &= \frac{\epsilon_i c_{e,i}}{\rho_w c_w} \,,
    &
      D_i &= \frac{\epsilon_i \lambda_{e,i}}{\rho_w c_w} \,,
    &
      \iota_i &= \frac{\gamma_{e,i}}{\rho_w c_w} \,,
  \end{align*}
  where $\epsilon_i$ is the fracture aperture, $\rho_w$ the fluid
  density, $c_w$ the fluid specific thermal capacity,
  $c_{e,i}=\phi_i\rho_w c_w+(1-\phi_i)\rho_m c_m$ the fracture effective
  thermal capacity, $\phi_i$ the fracture porosity, $\rho_m$
  the density of the rock matrix, $c_m$ the specific thermal capacity
  of the rock matrix,
  $\lambda_{e,i} = \lambda_w^{\phi_i}\lambda_m^{1-\phi_i}$ the
  effective thermal conductivity, $\gamma_{e,i}$ the effective heat
  transfer coefficient between fluid and rock, and, finally, $\hat{\theta}_i$ is the
  temperature of the rock matrix, acting as external heat
  source/sink. When $\phi_i = 1$ the fracture is completely open and no
  sediments are present, whereas values of $\phi_i < 1$ indicate that fracture is partially filled with
  crystals/rocks. Equation \eqref{eq:heat_strong_eq} is in
  conservative form, being the total flux the conserved quantity.

  For the sake of brevity and clarity of explanation we do not present boundary conditions in the most general context, but conforming to the numerical tests of Section~\ref{sec:examples}; generalization is however again straightforward. Recalling that $\partial \Omega_{i,D}$ is the Dirichlet portion of the boundary of fracture $\Omega_i$ for the Darcy problem in Section~\ref{subsec:darcy}, let us split $\partial \Omega_{i,D}$ into two parts, namely $\partial \Omega_{i,D}^{\rm inflow}$ and  $\partial \Omega_{i,D}^{\rm outflow}$. Thus, according to the computed Darcy velocity $\bm{u}_i$, on each fracture $\Omega_i$, the
  inflow boundary $\partial \Omega_{i,D}^{\rm inflow}$ is the portion of $\partial \Omega_{i,D}$ where
  $\bm{u}_i \cdot \bm{n}_{\rm ext} < 0$ and conversely the outflow boundary is the portion of $\partial \Omega_{i,D}$ where $\bm{u}_i \cdot \bm{n}_{\rm ext} > 0$, thus linking the nature of the boundary to the solution of the Darcy problem. Please note that $\partial \Omega_{i,D}^{\rm inflow}$ and $\partial \Omega_{i,D}^{\rm outflow}$ might be both empty for most of the fractures. Then  boundary conditions on fracture $\Omega_i$ are:
  \begin{gather}\label{eq:heat_strong_bc}
    \begin{aligned}
      & \tr_i \theta_i = \overline{\theta}_i && \text{on } \partial
      \Omega_{i, D}^{\rm inflow}\times (0, T] \,,
      \\
      & \tr_i D_i \nabla \theta_i \cdot \bm{n}_{\rm ext} = 0 &&
      \text{on } \partial \Omega_{i, D}^{\rm outflow}\times (0, T] \,,
      \\
      & \tr_i D_i \nabla \theta_i \cdot \bm{n}_{\rm ext} =0 &&
      \text{on } \partial \Omega_{i, N}\times (0, T] \,.
    \end{aligned}
  \end{gather}
  By \eqref{eq:darcy_strong_bc}, on the portion $\partial \Omega_{i, N}$
  both the diffusive and advective terms are zero for $\theta_i$, so
  $\partial \Omega_{i, N}$ represents an impervious boundary also for
  the heat equation. Finally, the initial condition on $\Omega_i$ is
  given by
  \begin{gather}\label{eq:heat_strong_ic}
    \theta_i = \overline{\overline{\theta_i}} \quad \text{in }
    \Omega_i \setminus \Gamma \times \{0\}.
  \end{gather}

  We consider now two distinct fractures $\Omega_i$ and $\Omega_j$
  with an intersection $\Gamma_k$. On $\Omega_i$ the temperature
  $\theta_i$ fulfils \eqref{eq:heat_strong_eq},
  \eqref{eq:heat_strong_bc}, and \eqref{eq:heat_strong_ic}. The same
  for $\theta_j$ in the fracture $\Omega_j$. The coupling conditions
  on $\Gamma_k \times (0,T]$ are
  \begin{equation}
    \label{eq:heat_strong_cc}
    \begin{split}
      &\sum_{l \in \{i, j\}} \!\!\tr_{l,+} (\bm{u}_l \theta_l -
      D_l \nabla \theta_l) \cdot \bm{n}_{l,+} + \!\tr_{l,-} (\bm{u}_l \theta_l
      - D_l \nabla \theta_l) \cdot \bm{n}_{l,-} \!=\! 0 \,,
      \\
      &\tr_{i,+} \theta_i = \tr_{i,-} \theta_i = \tr_{j,+} \theta_j =
      \tr_{j,-} \theta_j \,.
    \end{split}
  \end{equation}
  Also in this case, in presence of multiple fractures the
  generalization of the model is immediate.
\end{subequations}
The heat equation can be formalized in the following problem.
\begin{problem}[Heat equation on a DFN]\label{pb:model_heat}
  Given the set of fractures $\Omega$ and traces $\Gamma$ and the
  Darcy velocity $\bm{u}$, find $\theta$ such that
  \eqref{eq:heat_strong_eq}-\eqref{eq:heat_strong_bc}-\eqref{eq:heat_strong_ic}-\eqref{eq:heat_strong_cc}
  are satisfied.
\end{problem}

To derive the weak formulation of Problem \ref{pb:model_heat} we need
to introduce the Bochner space
\begin{gather*}
  W_i^0 = L^2\left(0, T; \left\{ v \in H^1(\Omega_i):\, \tr_i v \, (t)
      = {\overline{\theta_i}}\, (t) \text{ on } \partial \Omega_{i,
        D}^{\rm inflow} \,\,\, \forall t \in (0, T] \right\}\right),
\end{gather*}
while for the test functions we define
\begin{gather*}
  W_{i}^1 = H^1\left(0, T; \left\{ v \in H^1(\Omega_i):\, \tr_i v \,
      (t) = 0 \text{ on } \partial \Omega_{i, D}^{\rm inflow} \,\,\,
      \forall t \in (0, T] \right\}\right) \,.
\end{gather*}
We need to introduce space-time forms associated to each fracture, in
particular we have: $m_i$, $d_i$, $g_i$, $r_i$ which are defined on
$W_i^0 \times W_{i}^1 \rightarrow \mathbb{R}$, and
$S_i : W_{i}^1 \rightarrow \mathbb{R}$ such that
\begin{gather*}
  m_i(\theta, v) = -\int_0^T (\theta, \partial_t \zeta_i
  v)_{\Omega_i}\qquad d_i(\theta, v) = \int_0^T (D_i^{\frac{1}{2}}
  \nabla \theta, D_i^{\frac{1}{2}}
  \nabla v)_{\Omega_i}\\
  g_i(\theta, v) = \int_0^T (\nabla \cdot (\bm{u}_i \theta),
  v)_{\Omega_i}\quad
  r_i(\theta, v) = \int_0^T (\iota_i \theta, v)_{\Omega_i}\\
  S_i(v) = \int_0^T (\iota_i \hat{\theta}_i, v)_{\Omega_i} +
  (\overline{\overline{\theta_i}}, v(0))_{\Omega_i}
\end{gather*}
The global functional spaces are
\begin{gather*}
  \begin{split}
    W^0 = \left\{v\colon v_i \in W_i^0 \quad\forall i \,,\,
      \tr_{i,+}^k v_i(t) = \tr_{i,-}^k v_i(t) = \tr_{j,+}^k v_j(t) =
      \tr_{j,-}^k v_j(t) \quad\right.&
    \\
    \quad\left.\forall t \in (0,T]\,,\, \forall k\,,\, t^{-1}(k) =
      (i,j) \right\}&\,,
  \end{split}
  \\
  \begin{split}
    W^1 = \left\{v\colon v_i \in W_{i}^1 \quad\forall i \,,\,
      \tr_{i,+}^k v_i(t) = \tr_{i,-}^k v_i(t) = \tr_{j,+}^k v_j(t) =
      \tr_{j,-}^k v_j(t) \quad\right.&
    \\
    \quad\left.\forall t \in (0,T]\,,\, \forall k\,,\, t^{-1}(k) =
      (i,j) \right\}&\,,
  \end{split}
\end{gather*}
and the global bilinear forms $m$, $d$, $g$, $r$ and $S$ are the sum
over all fractures of the local ones. We finally have the weak
formulation of the problem.
\begin{problem}[Heat equation on a DFN - weak
  formulation]\label{pb:model_heat_weak}
  The weak formulation of Problem \ref{pb:model_heat} is: find
  $\theta \in W^0$ such that
  \begin{align*}
    m(\theta, v) + d(\theta, v) + g(\theta, v) + r(\theta, v) = S(v)\qquad
    \forall v \in W^1 \,.
  \end{align*}
\end{problem}
Also in this case to obtain a symmetric problem we can use a lifting technique of the Dirichlet datum.



\section{Numerical discretization}\label{sec:discretization}

In this section, we present various discretization strategies, both
well established, both unconventional, that can be used to solve the
models described in Section~\ref{sec:model}. These strategies have
similarities and differences that can be used to categorize them.  A
first point concerns the computational mesh and, in particular, how
the meshing is performed at fracture intersections: it is possible to
have a matching or non-matching grids. In the former case fracture
grids are conforming to the intersections among fractures, while in
the latter, grid elements arbitrarily cross intersections. A second
issue is related to mass conservation: in computing the Darcy
velocity, some schemes are locally mass conservative and some other
are only globally conservative, and this property may impact the
subsequent solution of the heat problem.  Also, some numerical schemes
are characterized by high numerical diffusivity which might impact the
solution but also avoids nonphysical spurious oscillations in
advection dominated flow regimes, whereas other schemes need to adopt
stabilitazion techniques.

Six different approaches are considered in the present work, given as
the combination of a numerical scheme for the computation of the Darcy
velocity and a numerical scheme for the spatial semi-discretization of
the subsequent non-stationary advection-diffusion-reaction problem
(shortly denoted as Heat equation). The implicit Euler method is used,
instead, in all cases, for time evolution. The approaches are listed
in Table~\ref{tab:tag}. The scheme tagged \TPFAUP{} is given by the
combination of the Two Point Flux Approximation (\TPFA{}) method for
the Darcy problem and the TPFA with upwinding for the advection term
(\TPFA{}+\UPWIND{}) for the Heat equation. Scheme \MFEMUP{} uses,
instead mixed finite elements (\MFEM{}) for the Darcy equation and
again \TPFA{}+\UPWIND{} for the Heat equation. The method \MVEMUP{},
uses instead the Virtual Element Method in mixed formulation (\MVEM{})
for the Darcy problem, on matching polygonal meshes. These schemes are
implemented in PorePy, a simulation tool written in Python for
fractured and deformable porous media, see
\cite{Keilegavlen2017a,Keilegavlen2017b}. PorePy is freely available
on GitHub along with the numerical tests proposed in Section
\ref{sec:examples}.  The method labeled \MFEMSUPG{} is based on mixed
standard finite elements for the numerical resolution of the Darcy
problem and on standard finite elements (\FEM{}) with Streamline
Upwind Petrov-Galerkin (\SUPG{}) stabilization
\cite{Franca-Frey-Hughes:1992} for the Heat equation. Finally methods
\FEMSUPG{} and \XFEMSUPG{} use a non-conventional optimization based
approach for the Darcy equation and, with \SUPG{} stabilization also
for the Heat equation. The optimization approach can adopt different
baseline discretization methods: here we consider the variants using
standard finite elements (\OPTFEM{}) and extended finite elements
(\OPTXFEM). Methods \MFEMSUPG{} \FEMSUPG{} and \XFEMSUPG{} are
implemented in C++ and Matlab$^\circledR$.

The forthcoming parts describe in more details the previous
approaches, grouping them according to the coupling at the traces:
matching coupling at traces in Subsection~\ref{subsec:matching} and
non-matching coupling in Subsection~\ref{subsec:nonmatching}.

\newcommand{\tagtablecell}[1]{#1}
\begin{table}
  \centering
  \begin{tabular}{c|ccc}
    Tag & Darcy eq.  & Heat eq.  & Grid
    \\\hline
    {\small \TPFAUP} & {\small \TPFA} & \tagtablecell{\small \TPFA{} + \UPWIND}
                                 & \tagtablecell{\small Matching triangles}
    \\\hline
    {\small \MFEMUP} & {\small \MFEM}         & \tagtablecell{\small \TPFA{} + \UPWIND}
                                 & \tagtablecell{\small Matching triangles}
    \\\hline
    {\small \MVEMUP} & {\small \MVEM}         & \tagtablecell{\small \TPFA{} + \UPWIND}
                                 & \tagtablecell{\small Matching polygons}
    \\\hline
    {\small \MFEMSUPG} & {\small \MFEM}  &  \tagtablecell{\small \FEM{} + \SUPG}
                                 & \tagtablecell{\small Matching triangles}
    \\\hline
    {\small \FEMSUPG} & {\small \OPTFEM}      & \tagtablecell{\small \OPTFEM{} + \SUPG}
                                 & \tagtablecell{\small Non-match. triangles}
    \\\hline
    {\small \XFEMSUPG} & {\small \OPTXFEM}     & \tagtablecell{\small \OPTXFEM{} + \SUPG}
                                 & \tagtablecell{\small Non-match. triangles}
  \end{tabular}
  \caption{List of numerical schemes proposed to solve the Darcy and
    heat problems, with the type of meshes used.}%
  \label{tab:tag}
\end{table}

\subsection{Matching discretization at traces}\label{subsec:matching}

\begin{figure}
  \centering
  \begin{subfigure}{.32\linewidth}
    \includegraphics[width=\linewidth]{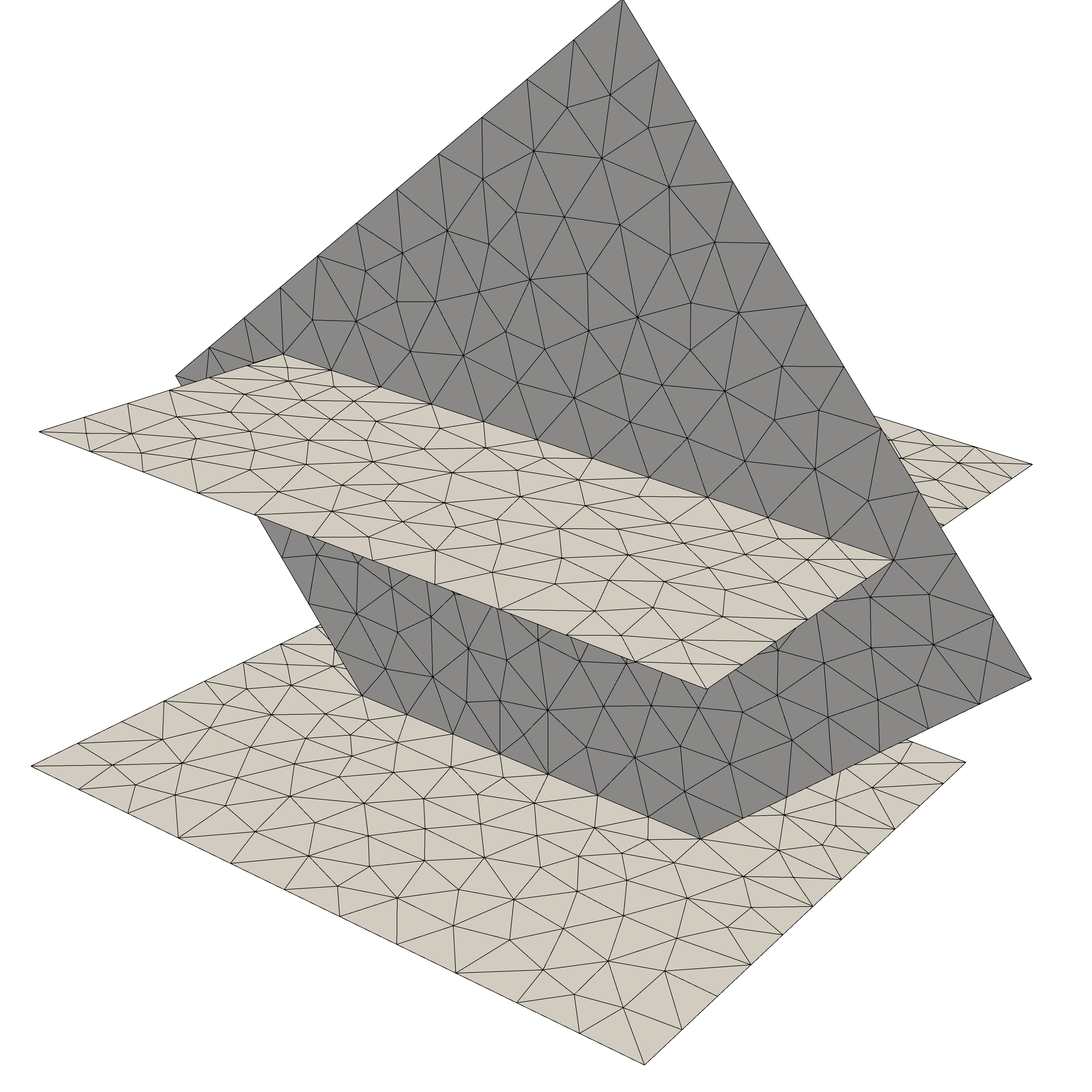}
    \caption{matching triangles}
    \label{fig:example-tri-conf-mesh}
  \end{subfigure}
  \hfill
  \begin{subfigure}{.32\linewidth}
    \includegraphics[width=\linewidth]{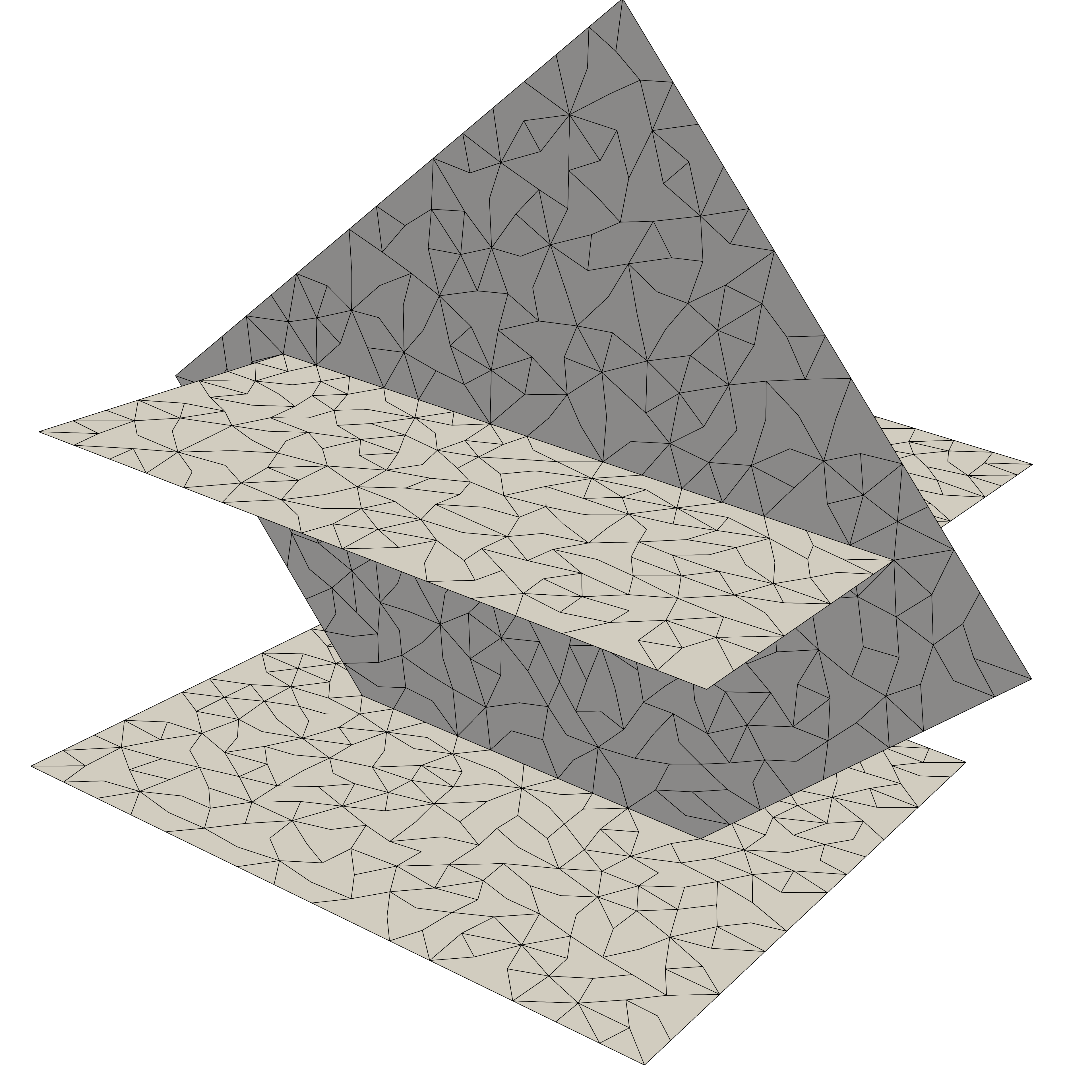}
    \caption{matching polygons obtained by coarsening}
    \label{fig:example-poly-conf-mesh}
  \end{subfigure}
  \hfill
  \begin{subfigure}{.32\linewidth}
    \includegraphics[width=\linewidth]{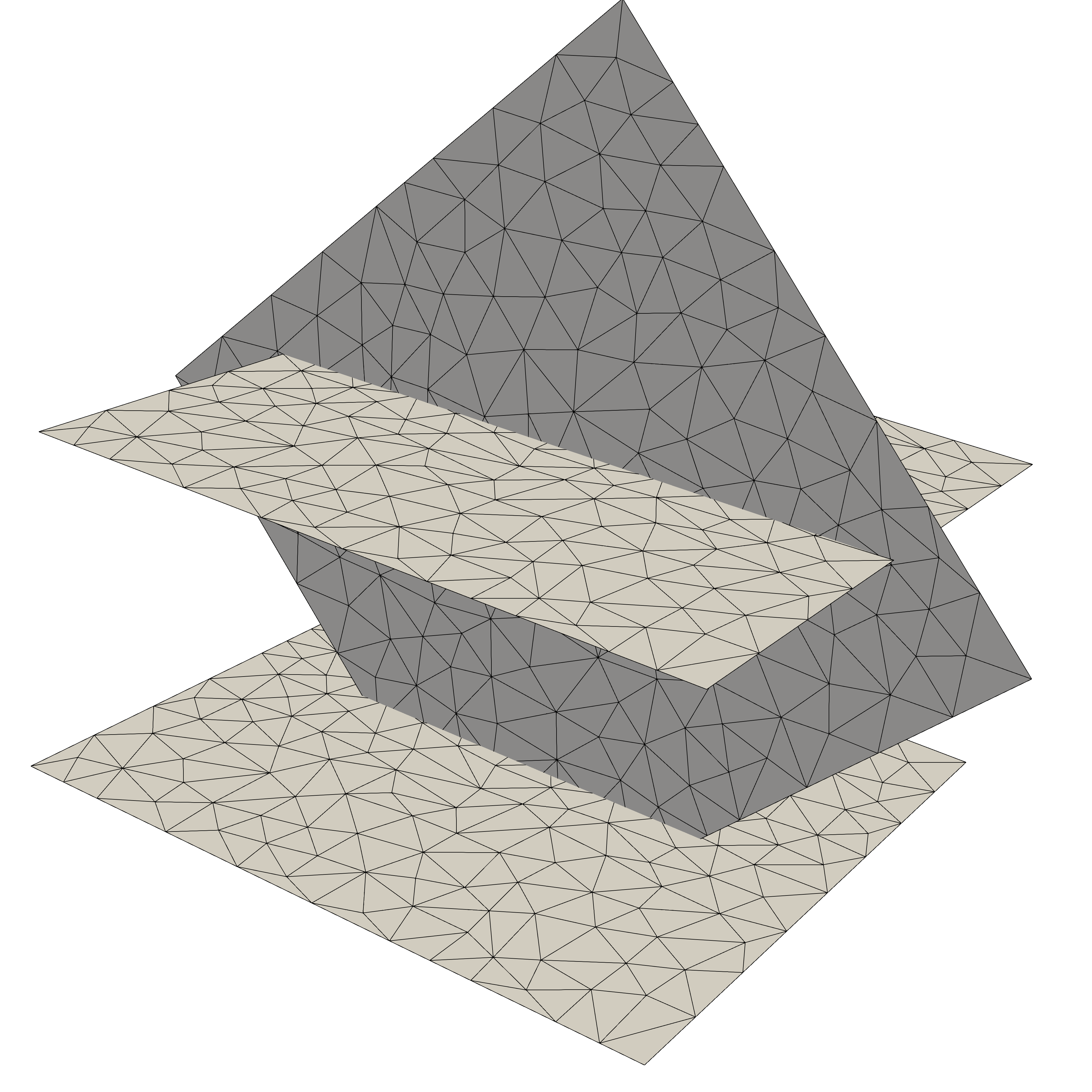}
    \caption{non-matching triangles}
    \label{fig:example-tri-nonconf-mesh}
  \end{subfigure}
  \caption{Examples of matching meshes.}
  \label{fig:example-mesh}
\end{figure}

Here advantages and drawbacks of a conforming discretization at the
traces are discussed. As mentioned before in a conforming grid, the
meshes of both the intersecting fractures match the trace with their
geometry. The trace is thus entirely covered by contiguous cell edges
of the two fractures, see
Figure~\ref{fig:example-tri-conf-mesh}-\ref{fig:example-poly-conf-mesh}
as an example. This approach has the clear advantage of an easy
applicability to most of the existing and well established numerical
methods (finite volumes, finite elements). However, in the case of
complex geometries the computational cost might increase and become a
severe constraint in complex fracture networks.

To solve the Darcy problem, we rely on different classes of numerical
schemes: finite volumes, finite elements in primal and mixed
formulation, and virtual elements in mixed formulation. For the finite
volume class, we choose the two-point flux approximation \TPFA{} on
simplicial grids, applied to Problem~\ref{pb:model_darcy_primal_weak},
see \cite{Aavatsmark2007a} for details.  This scheme is well known in
the industry field and widely used for its velocity in assembling the
discrete problem and for having a narrow matrix stencil. The scheme is
locally conservative and robust with respect to strong variations of
the hydraulic conductivity coefficient, however it is consistent only
for $K$-orthogonal grids. Regarding the class of finite elements, we
consider mixed finite elements with the pair of spaces
$\mathbb{RT}_0-\mathbb{P}_0$ for the Darcy velocity and piecewise
constants for the hydraulic head, denoted by \MFEM{}, defined on a
simplicial grid, for further references see
\cite{Raviart1977,Roberts1991}. It is well known that
$\mathbb{RT}_0-\mathbb{P}_0$ is locally conservative and more robust
than primal $\mathbb{P}_1$ finite elements with respect to strong
variations on the hydraulic conductivity coefficient. On the other
hand, it gives a larger linear system, with a saddle-point
structure. The \MFEM{} scheme solves the Darcy problem in the form
presented in Problem \ref{pb:model_darcy_weak}. In some particular
scenarios the regularity requirements on meshes formed by triangles
are too restrictive and schemes able to handle generally shaped cells
are more suitable. In these cases we rely on the new class of virtual
element methods, which are variational methods where the basis
functions of the discrete spaces are not prescribed a-priori, and are
defined implicitly on general star-shaped elements as solutions of
suitable local PDEs. See the seminal works
\cite{BeiraodaVeiga2013a,Brezzi2014,BeiraodaVeiga2014a,BeiraoVeiga2016}
and those related to DFN
\cite{Benedetto2014,Benedetto2016,Benedetto2016b,Benedetto2017,Fumagalli2016a}. In
our analysis we consider only the lowest order mixed (\MVEM)
formulation, which can be viewed as a generalization
$\mathbb{RT}_0-\mathbb{P}_0$ mixed finite elements on generally shaped
cells, solving Problem \ref{pb:model_darcy_weak}. Virtual element
methods share many properties with their finite element counterpart:
indeed, \MVEM{} are locally conservative, robust with respect to the
hydraulic conductivity variation and have the same grid stencil of the
$\mathbb{RT}_0-\mathbb{P}_0$. Here we use \MVEM{} on polygonal grids,
obtained by coarsening a mesh originally made of triangular elements
in order to reduce the number of cells required in the simulation for
complex geometries, as described in
\cite{Fumagalli2016a,Fumagalli2017}.

To solve the heat equation (Problem \ref{pb:model_heat_weak}), we use
a \TPFA{} scheme for the diffusive term and a weighted upwind scheme
for the advective term. The advantage of this choice is that we obtain
a stable scheme which respects the maximum and minimum principle,
without oscillation due to high grid-P\'eclet numbers. However, for
some applications the obtained scheme might be too diffusive. Finally,
we consider also standard $\mathbb{P}_1$ finite elements, with a
\SUPG{} stabilized discrete variational formulation, which is a
globally conservative numerical scheme.

Since in all the above cases we have matching meshes, and thus the
degrees of freedom on the trace of one fracture correspond to the
degrees of freedom of the intersecting one, coupling conditions can be
imposed strongly.

Methods \TPFAUP{}, \MFEMUP{} and \MFEMSUPG{} are considered here
standard reference approaches of different discretization strategies
(finite volume and finite element based), and are used to benchmark
the behavior of the other less conventional approaches, based e.g. on
polygonal or non-matching discretizations.


\subsection{Non-matching discretization at traces}
\label{subsec:nonmatching}

When dealing with huge networks, the generation of conforming meshes
may require a large computational cost. Then, it is worth considering
a class of methods that do not require any kind of conformity of the
fracture meshes to traces, see, for example, Figure
\ref{fig:example-tri-nonconf-mesh}.

In \cite{Berrone2013,Berrone2013a,BPSc,BPSd,BBoS,BBV2019} a
PDE-constrained optimization approach is proposed, based on
non-conforming meshes, that can be applied both to Problem
\ref{pb:model_darcy_primal_weak} and Problem
\ref{pb:model_heat_weak}. In this framework, the problem is rewritten
as a minimization problem for a functional measuring the error in
fulfilling matching conditions, constrained by local PDEs on each
fracture. This approach provides not only a numerical approximation of
the solution but also a directly computed approximation of the flux
exchanged at traces, which is of interest for many applications.  The
discretization can be based on different methods: standard
$\mathbb{P}_1$ finite elements are the simplest choice, and the
resulting scheme is denoted as \OPTFEM{}. However, as mesh elements
arbitrarily cross the traces, the jump of the co-normal derivative of
the solution at fracture intersections, still directly
computed by the method, can not be correctly represented by non
conforming $\mathbb{P}_1$ finite elements. Thus the use of local
extended finite elements is also considered, being at the basis of the
method denoted \OPTXFEM{}. When used for advection dominated flow
regimes, the \SUPG{}-stabilized versions of the method are used
(\FEMSUPG{}, \XFEMSUPG{}).

The resulting numerical schemes inherit the mass-conservation
properties of the local discrete formulations, thus they are globally
but not locally conservative. A huge advantage of this method is that
the matrices resulting from the local discrete problems can be
computed in parallel, and also the solution can be computed strongly
relying on parallel computing. In \cite{BSV} a MPI-based parallel
algorithm is proposed for assembling and solving the discrete problem
using the conjugate gradient method on huge networks of fractures,
whereas in \cite{BDV} an implementation suitable for GPGPUs is
presented.




\section{Examples}\label{sec:examples}

In this section, we present three test cases with the aim of validating and
comparing the previously introduced models. Extending the work proposed in
\cite{Scialo2017}, here we analyse the various approaches for time dependent problems
of advection-diffusion-reaction, where the advection velocity is computed by
means of the same approach solving a diffusion problem. The key aspects of the
various schemes will be highlighted and investigated, along with the impact of
the lack of conservation of fluxes, both locally and through a trace, that
characterizes some of the proposed approaches. For the considered test cases,
both local and global quantities will be computed at different time steps, and
used to assess and compare the behaviour of the various approaches, such as
\textit{i)} the integral mean in space of the temperature on a fracture
$\Omega$, denoted as $\langle {\theta} \rangle_{\Omega}$; \textit{ii)} the total flux
mismatch on a trace, $\delta \Phi_\Gamma$ defined as the integral of the sum of
the net total fluxes $\Phi_{\Gamma,i}$ and $\Phi_{\Gamma,j}$ entering or leaving
the two fractures $\Omega_i$ and $\Omega_j$ meeting at trace $\Gamma$,
respectively, i.e. $\delta \Phi=\vert{\Phi_{\Gamma,i} + \Phi_{\Gamma,j}}\vert$,
$$
    \Phi_{\Gamma,i}=-\int_{\Gamma}\left(\jumpS{\tr_i D_i\nabla\theta_i \cdot
    \bm{n}_{i\Gamma}}
    + \jumpS{\tr_i \theta_i \tr_i \bm{u_i}\cdot \bm{n}_{i\Gamma}}\right)
$$
with $\bm{n}_{i\Gamma}$ the unit normal vectors to $\Gamma$, with fixed orientation
on $\Omega_i$; and \textit{iii)} the averaged $\theta$ on the outflow boundary, $\partial
\Omega^{\mathrm{outflow}}_D$, denoted as
\begin{gather*}
    \langle {\theta} \rangle_{\mathrm{outflow}}
    :=\dfrac{1}{\vert{\partial \Omega^{\mathrm{outflow}}_D}\vert}
    \int_{\partial \Omega^{\mathrm{outflow}}_D} \tr \theta.
\end{gather*}

The considered test cases are designed in order to challenge the proposed
approaches with complex geometries and/or realistic models. In
Subsection~\ref{subsec:example1} we consider the effect of vanishing trace from
a simple network, by a sequence of simulations. In the second example presented
in Subsection~\ref{subsec:example2} a synthetic small network of $10$ fractures
is considered, to analyse the behaviour of the methods on a more general, yet
simple configuration. We finally conclude with a realistic example in
Subsection~\ref{subsec:example3}, where a network of 89 fractures is generated
from an extrusion of an interpreted natural outcrop and physically sound values
for the various parameters are used.

With the purpose of establishing a standard reference for the analysis of
numerical schemes for flow in fractures, test cases in Subsections
\ref{subsec:example1} and \ref{subsec:example3} are borrowed from
\cite{Scialo2017} and adapted to the present context, all the geometrical data
being available in \cite{FKSdata}.


\subsection{Vanishing trace between intersecting fractures}\label{subsec:example1}

As a first test case, named Test Case 1, the same setting of the problem
proposed in the example of \cite[Subsection 4.3.1]{Scialo2017} is considered. In
the reference, the Darcy problem was tackled, whereas here a non-stationary
advection-diffusion problem for the passive scalar $\theta$ is solved.

\begin{figure}[htbp]
    \centering
    \resizebox{0.33\textwidth}{!}{\fontsize{15pt}{8}\selectfont%
    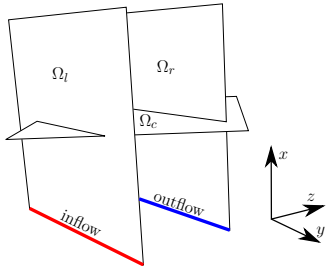}
    \caption{Geometry of Test Case 1. The
    inflow and outflow portions of the boundary are represented by red and blue
    lines, respectively.}%
    \label{fig:geometry_example_1}
\end{figure}

\begin{figure}
    \centering
    \includegraphics[width=0.25\textwidth]{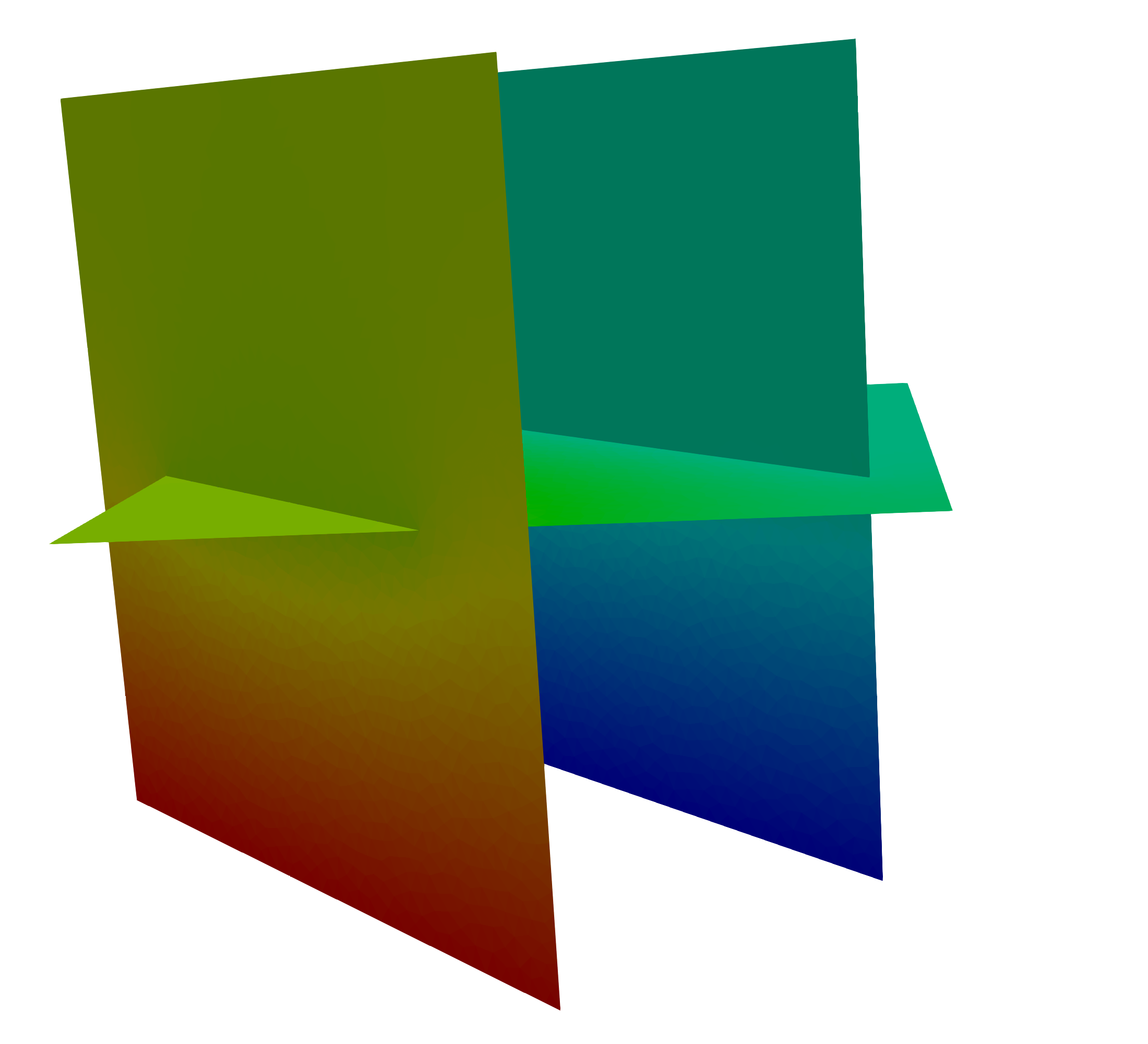}%
    \includegraphics[width=0.25\textwidth]{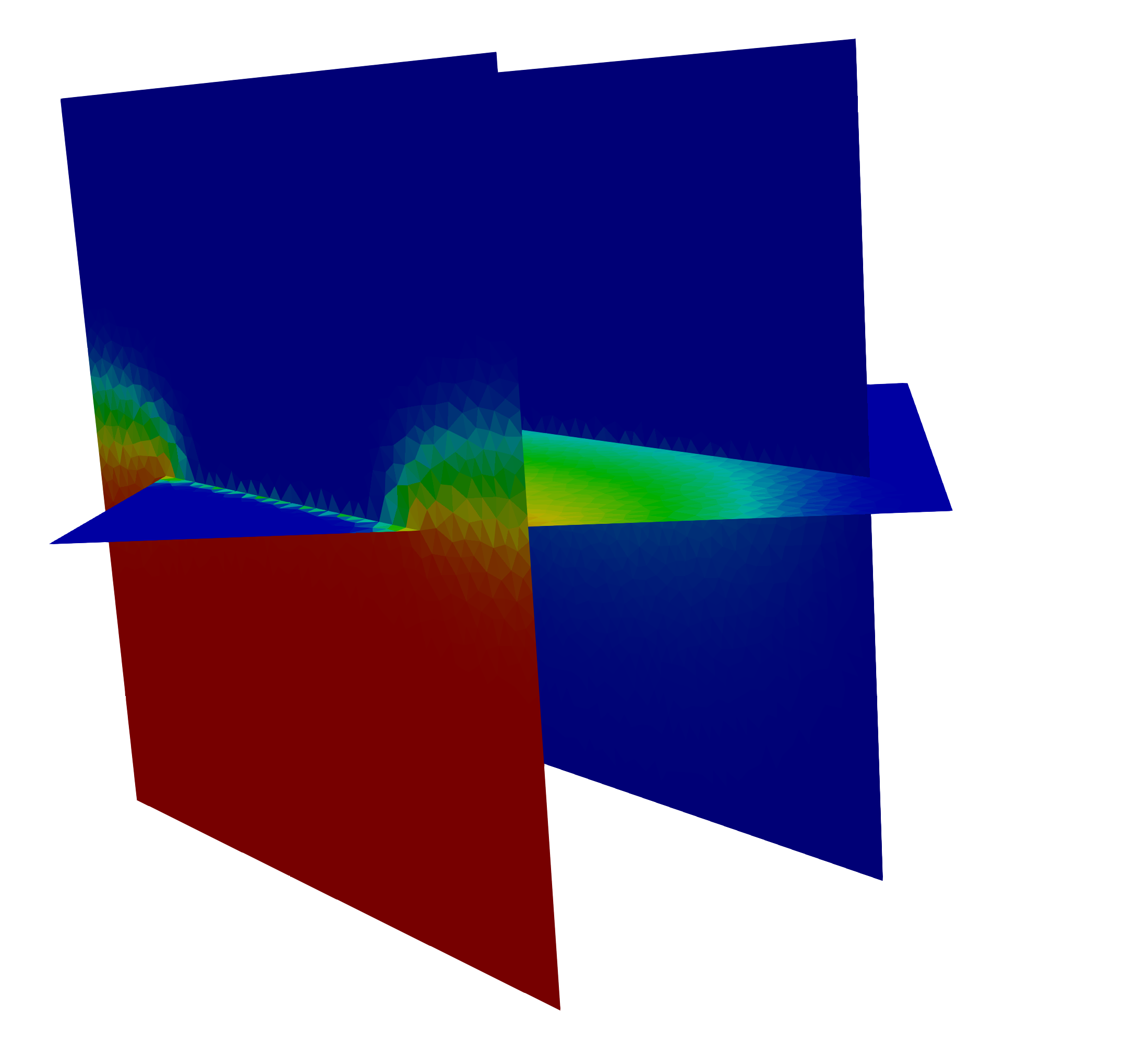}%
    \includegraphics[width=0.25\textwidth]{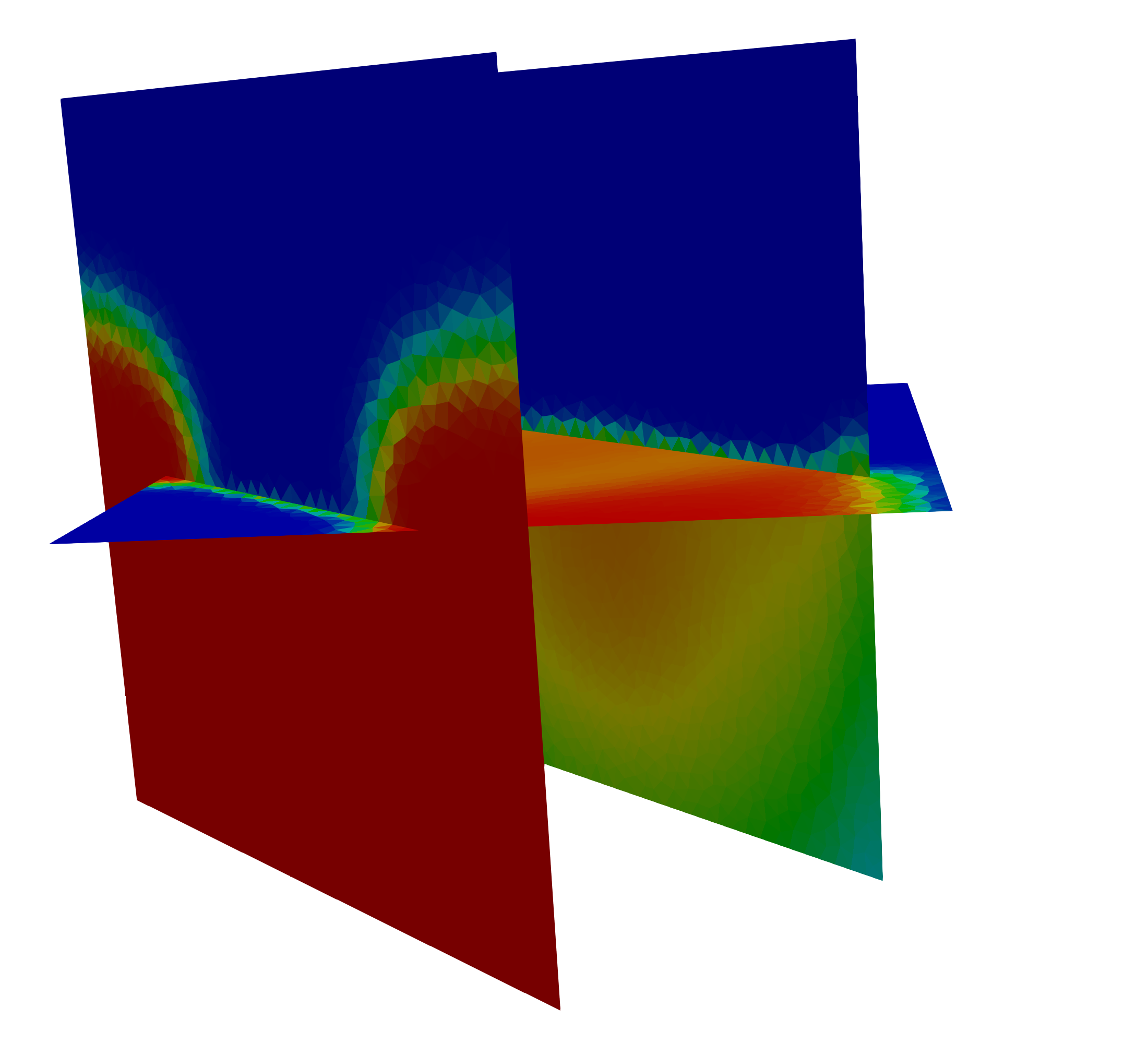}%
    \includegraphics[width=0.25\textwidth]{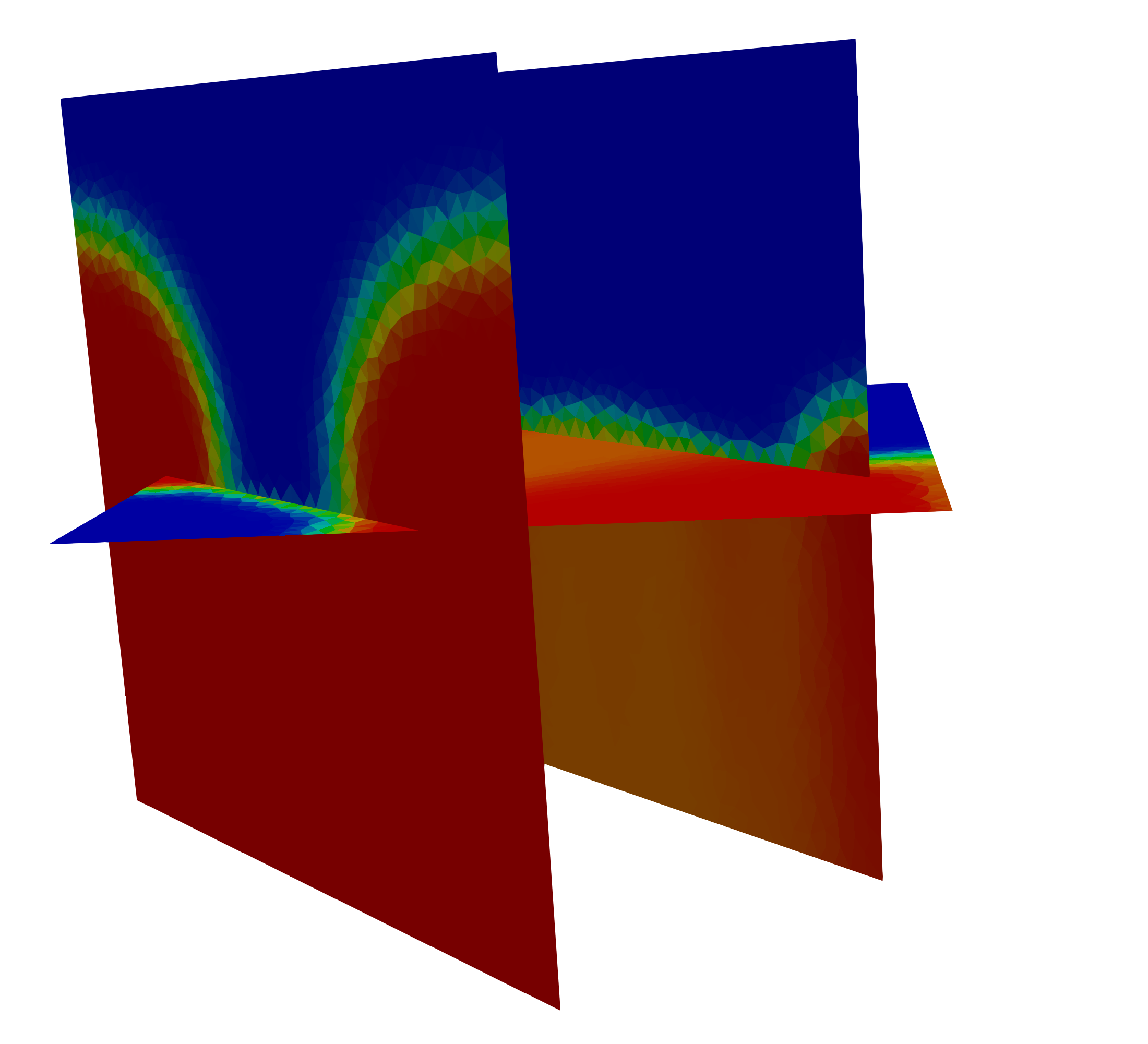}\\
    \includegraphics[width=0.25\textwidth]{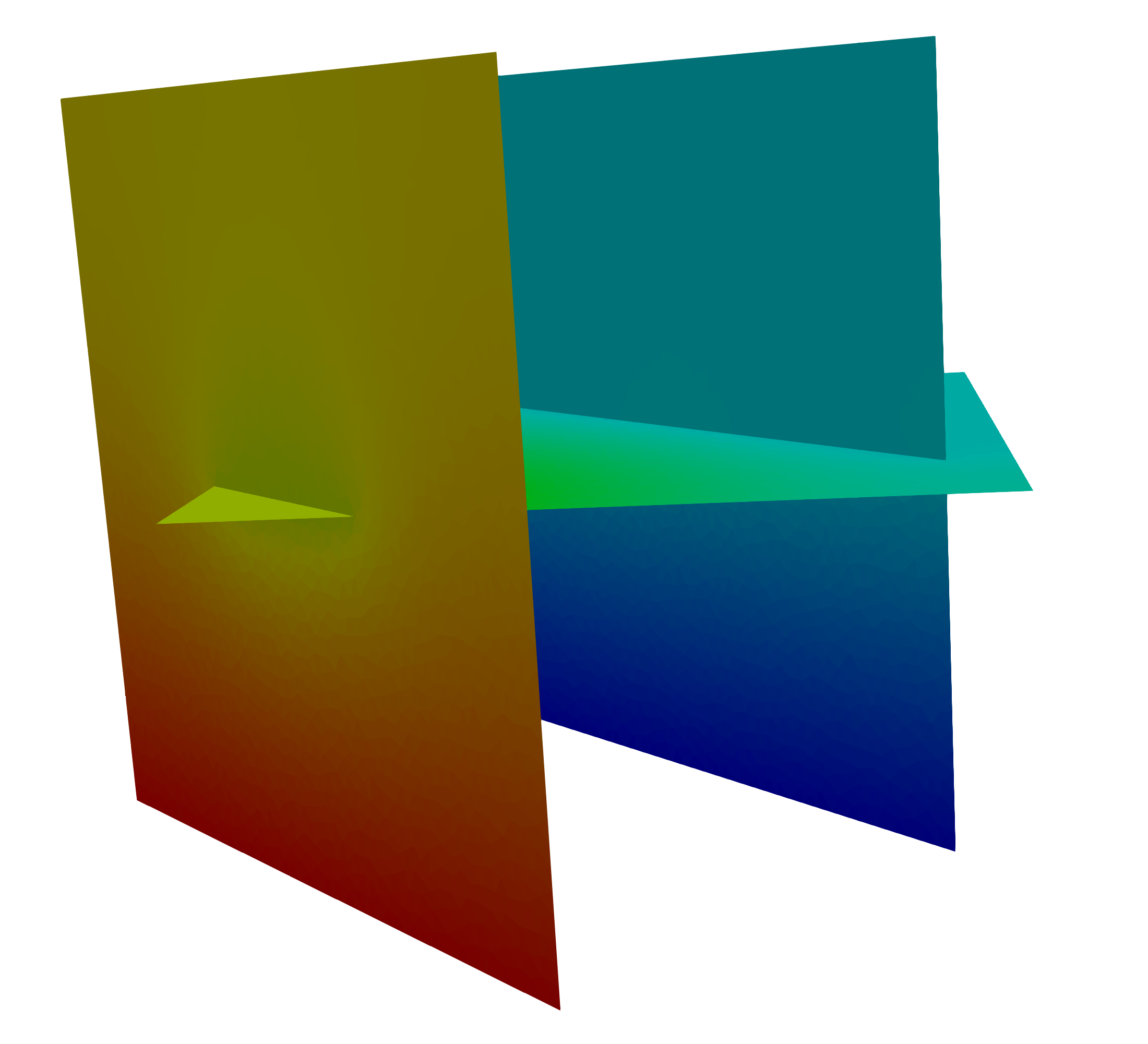}%
    \includegraphics[width=0.25\textwidth]{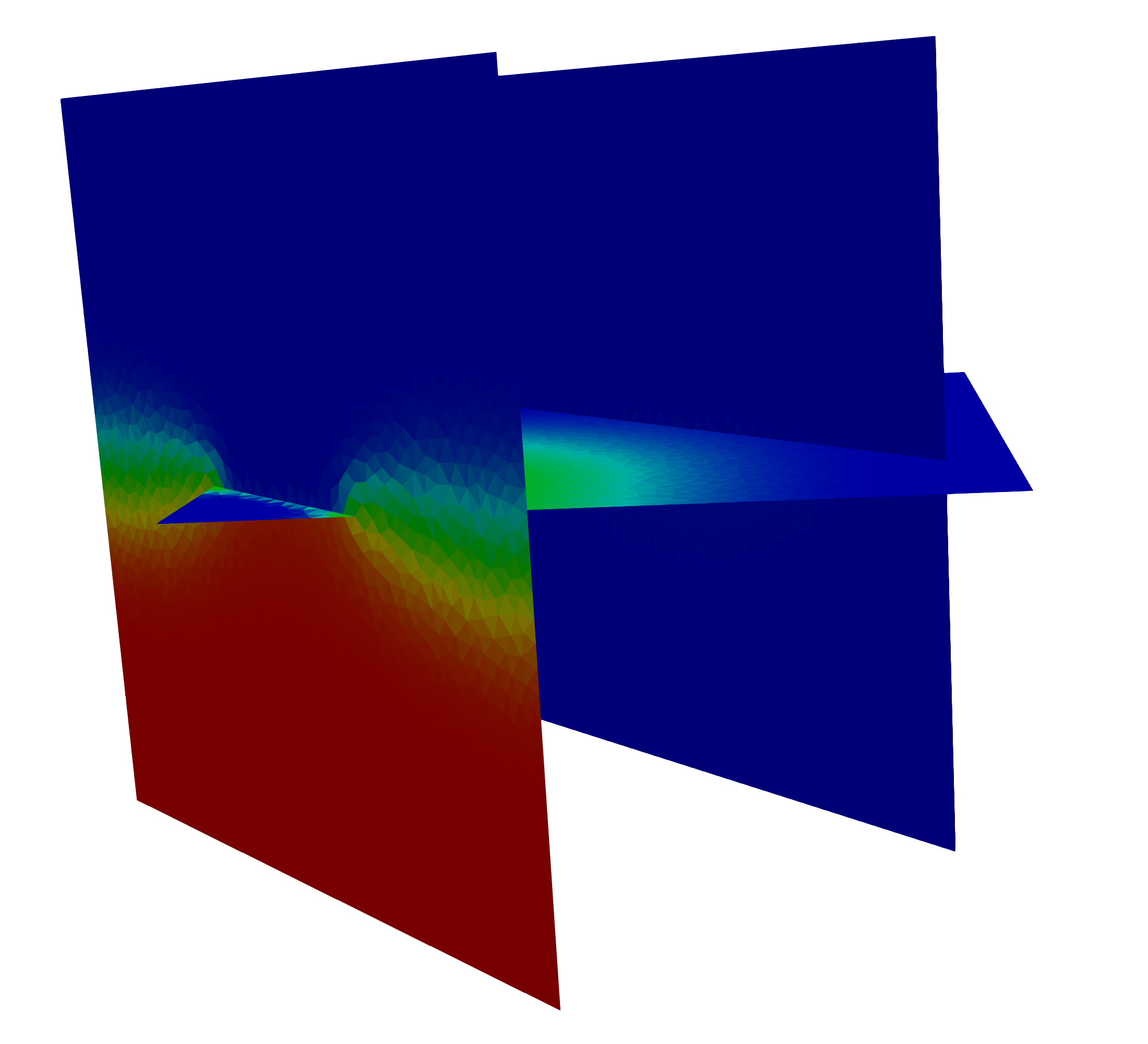}%
    \includegraphics[width=0.25\textwidth]{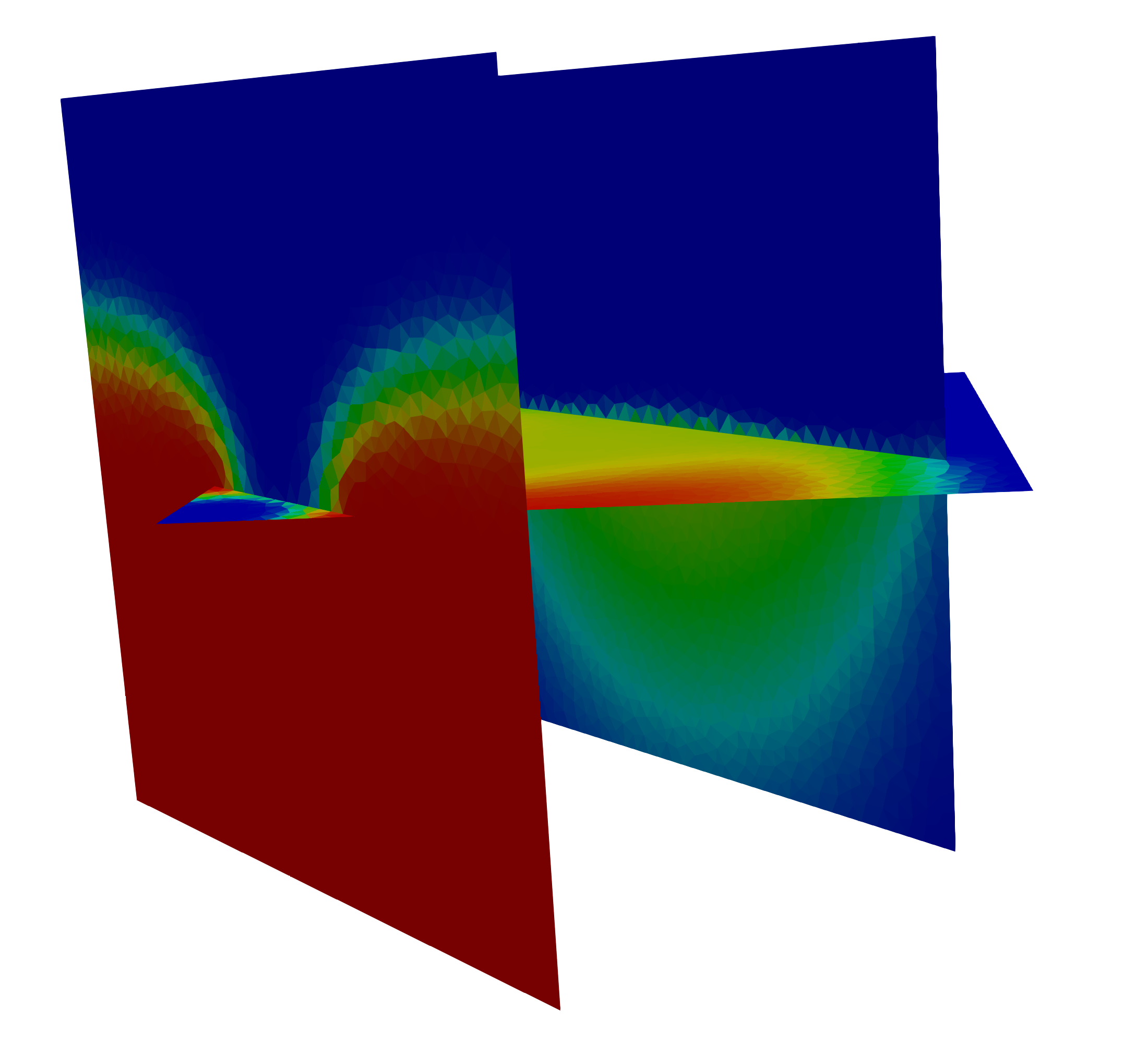}%
    \includegraphics[width=0.25\textwidth]{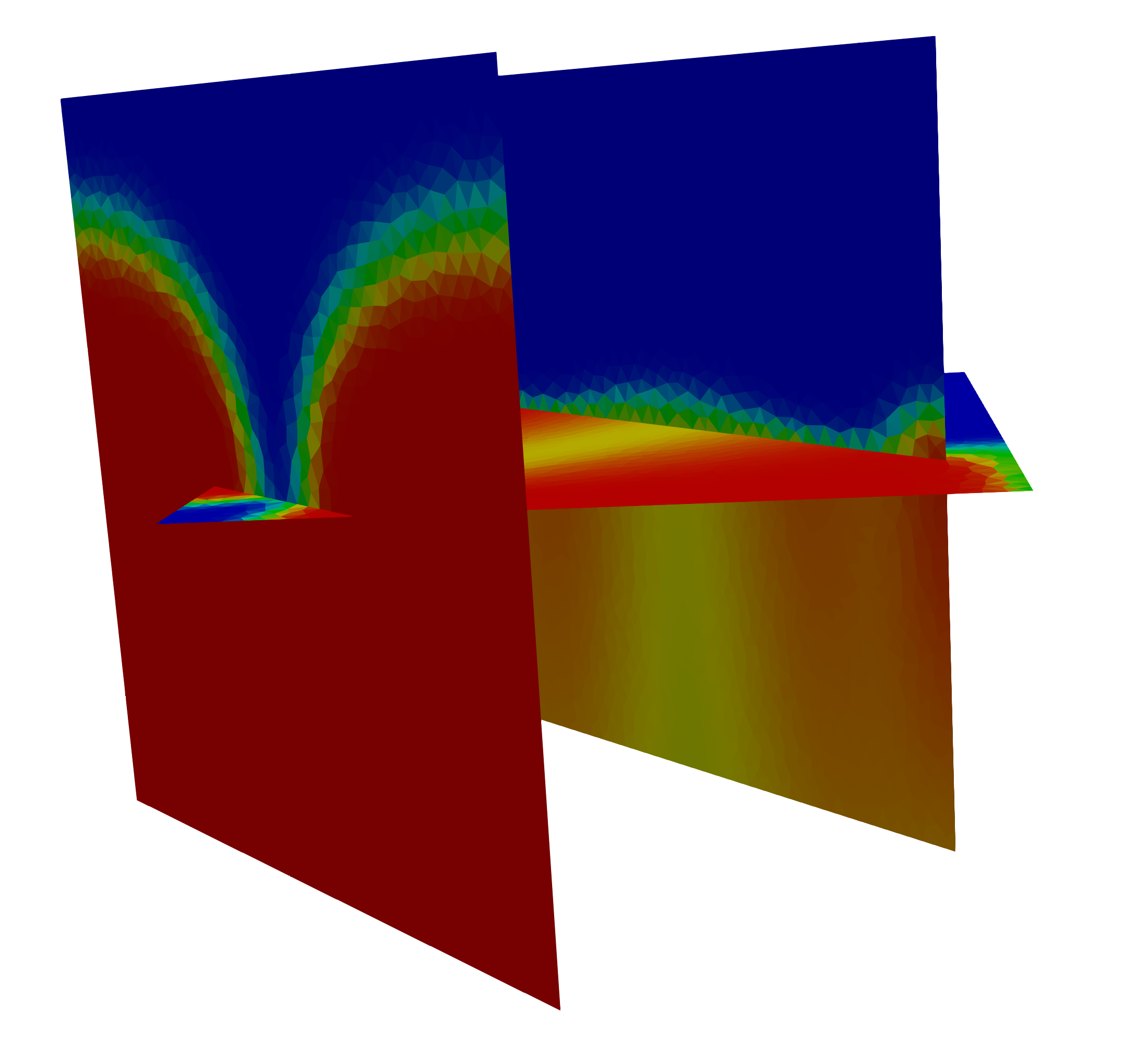}\\
    \includegraphics[width=0.25\textwidth]{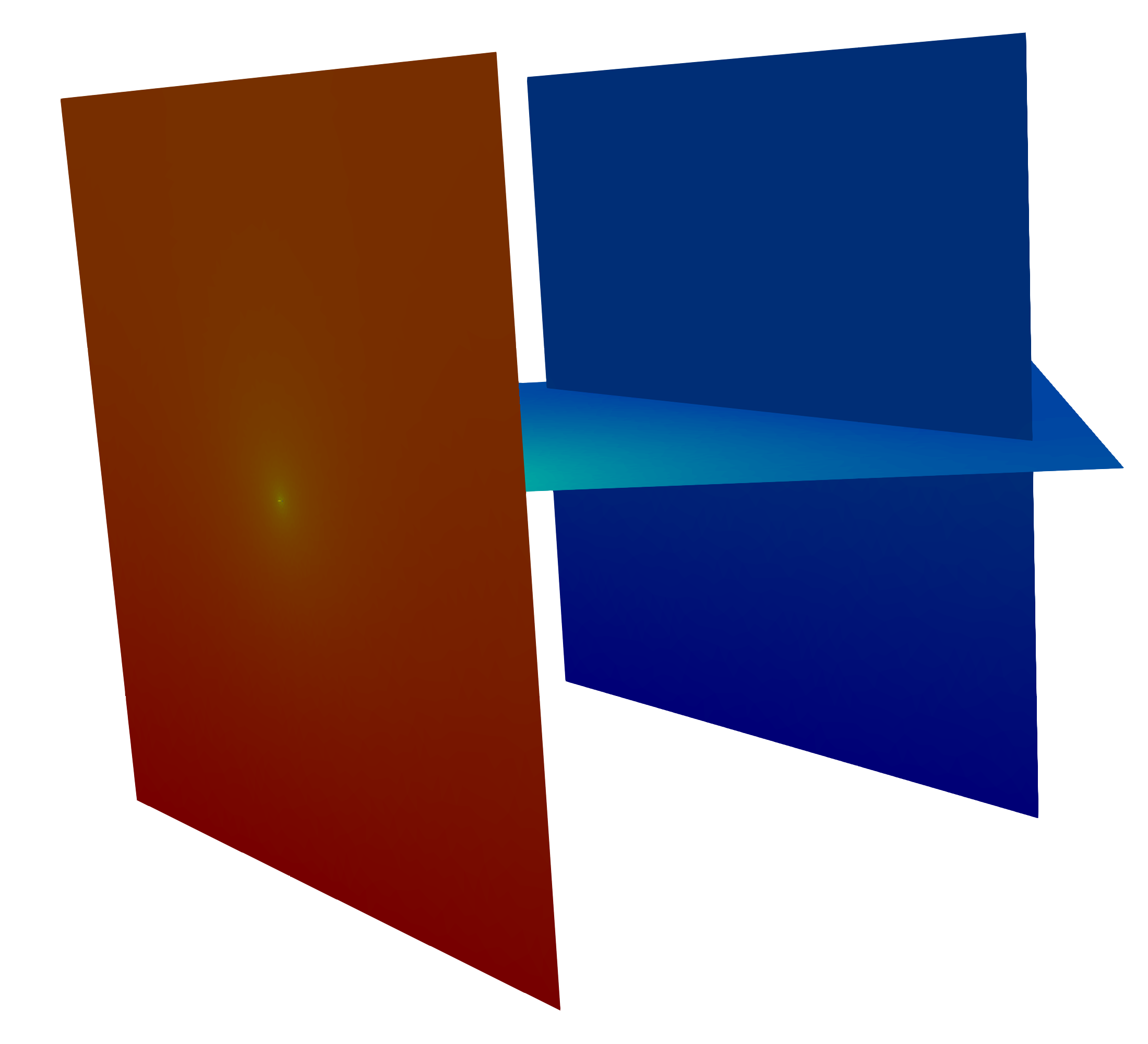}%
    \includegraphics[width=0.25\textwidth]{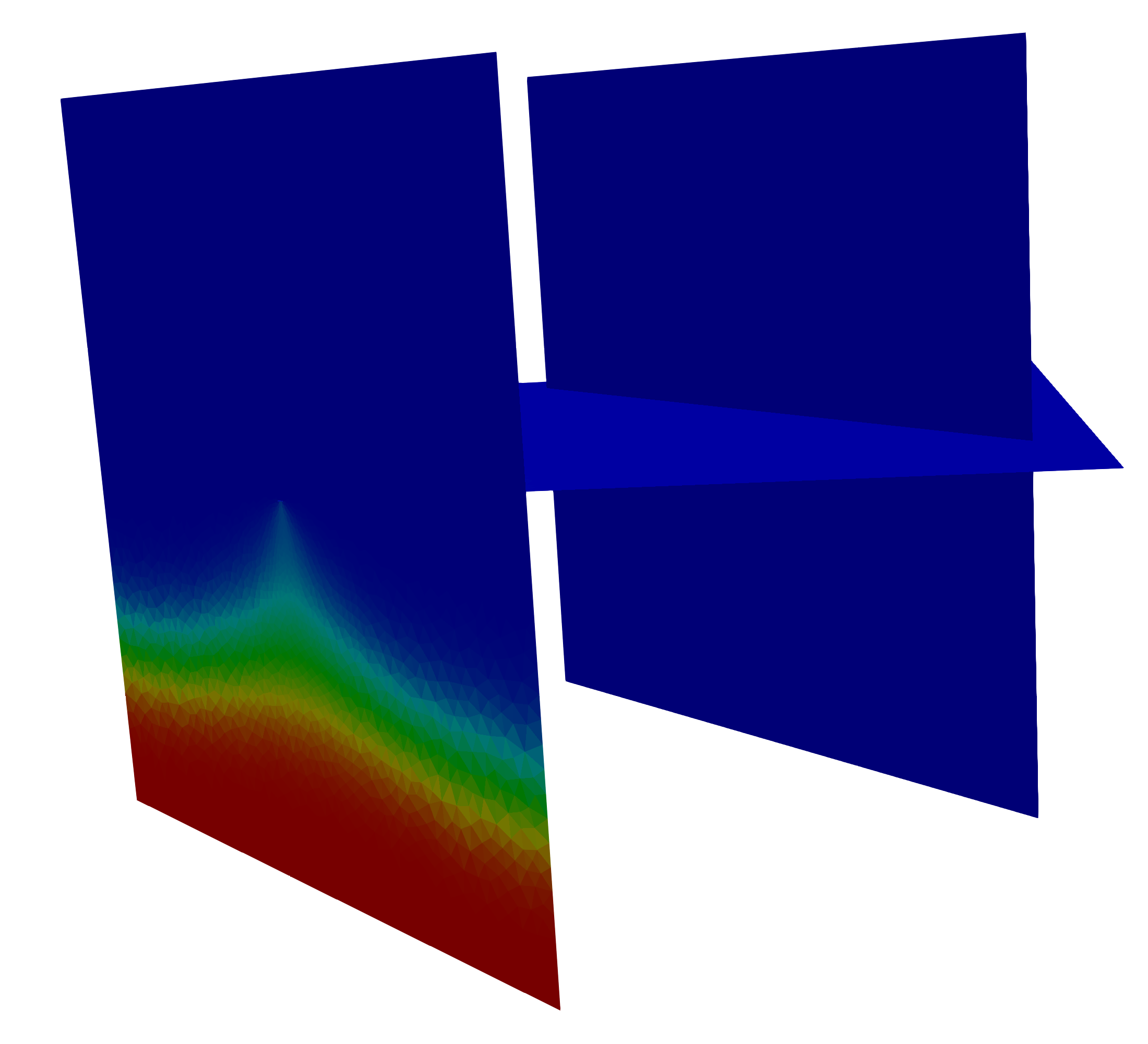}%
    \includegraphics[width=0.25\textwidth]{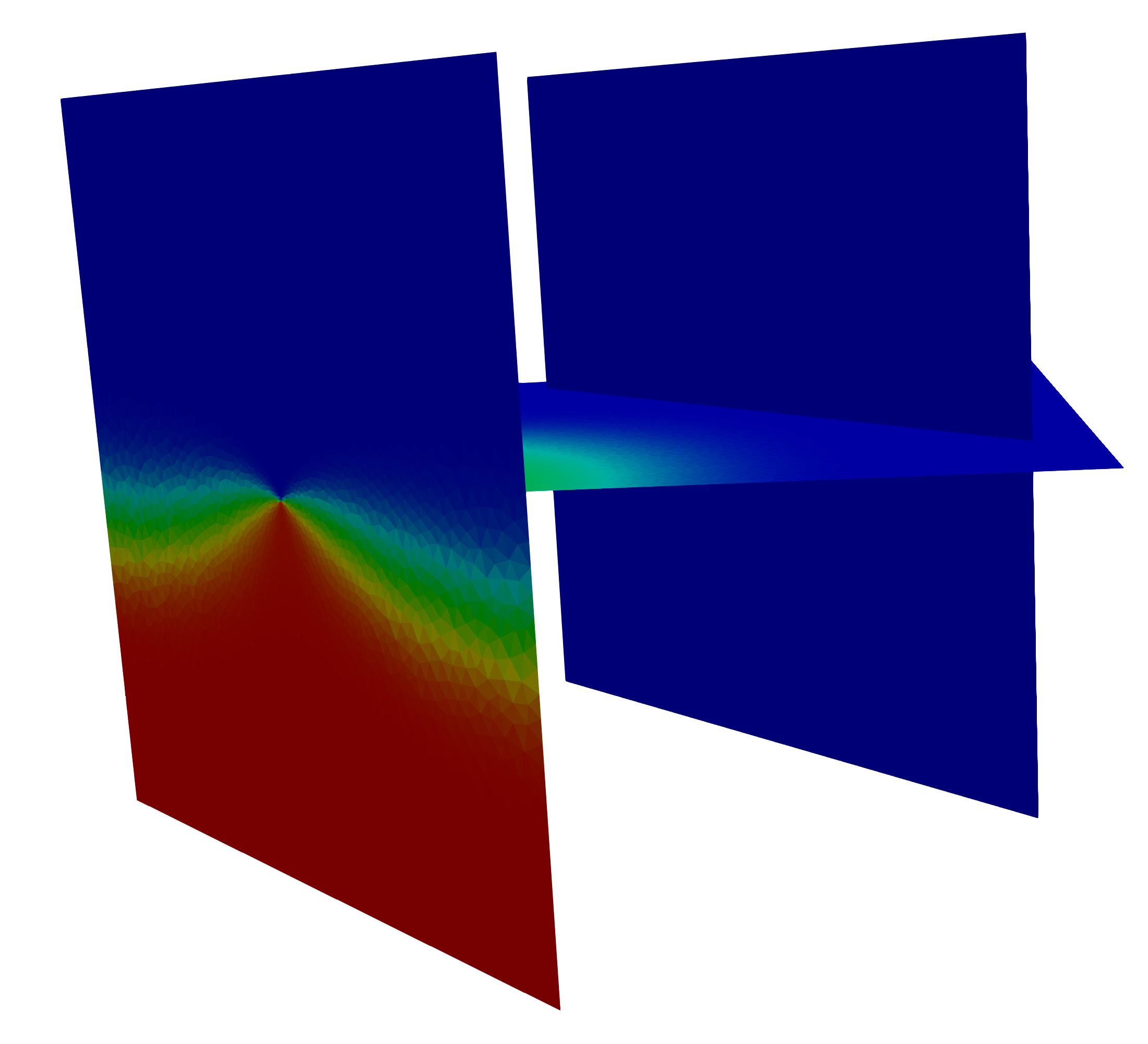}%
    \includegraphics[width=0.25\textwidth]{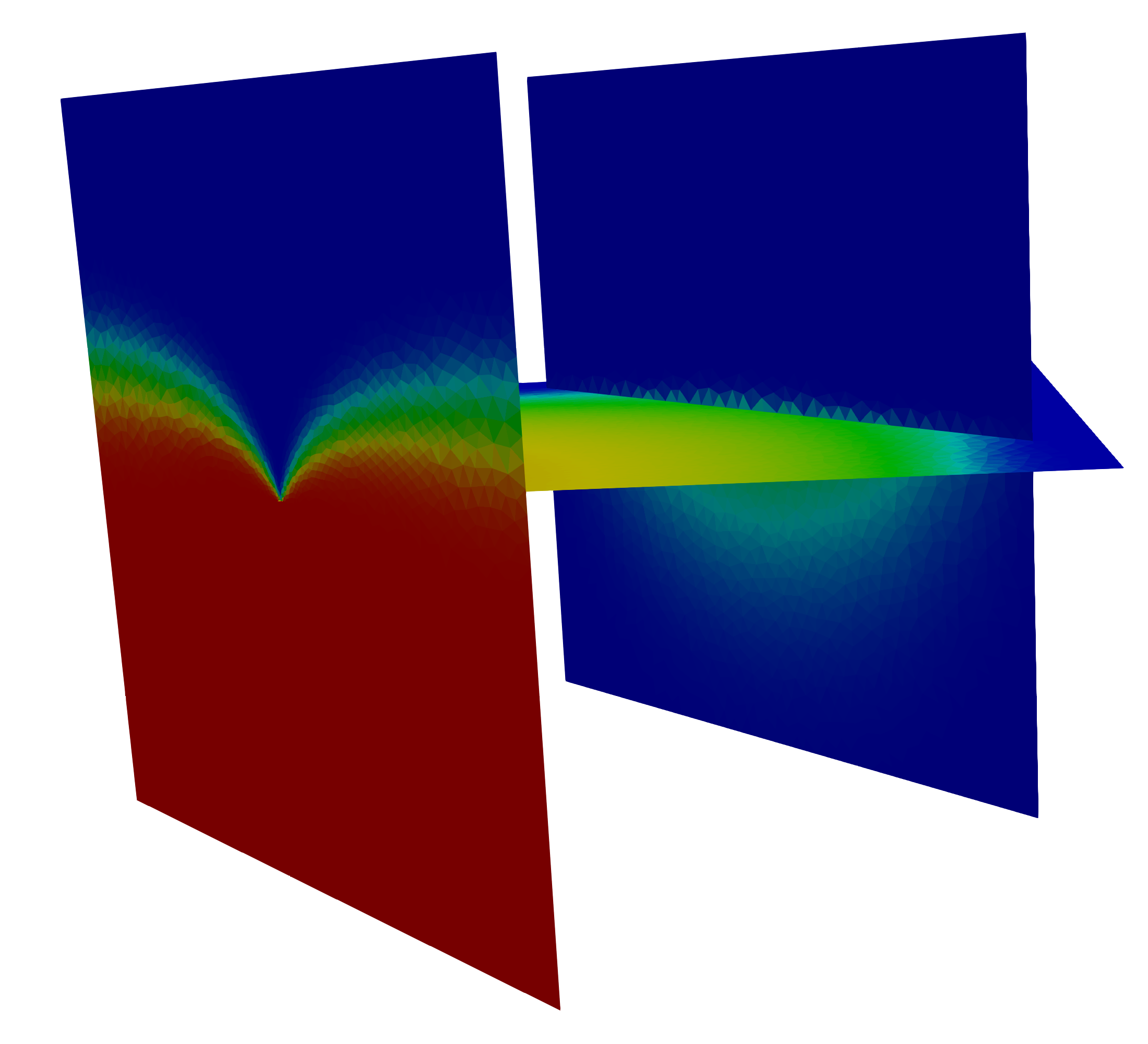}\\
    \caption{Solution of Test Case 1. On the pressure head solution is shown on
    the leftmost column at configurations $C0$, $C10$ and $C20$. On the
    remaining columns the solution $\theta$ is represented at $t=\{1.25, 2.5,
    5\}$, respectively. The colour scale spans from $0$ to $1$.}%
    \label{fig:solution_example_1}
\end{figure}

In this test case the same problem is solved on different geometries. Let us
consider a network composed of three fractures, named $\Omega_l$, $\Omega_r$,
and $\Omega_c$, as shown in Figure~\ref{fig:geometry_example_1}. Fracture
$\Omega_l$ has a fixed position, whereas fractures $\Omega_r$, and $\Omega_c$
are displaced, for each different geometry of a same distance along the
$z$-direction, such that the length of the intersection line between fractures
$\Omega_l$ and $\Omega_r$, denoted as $\Gamma_0$, progressively reduces from the
configuration shown in Figure~\ref{fig:solution_example_1}, being instead fixed
the intersection between $\Omega_r$, and $\Omega_c$, denoted as $\Gamma_1$. In
such a way 21 different configurations are obtained, with the length of the
vanishing trace $\Gamma_0$ ranging from $0.6$ at configuration $C0$, as shown in
Figure~\ref{fig:geometry_example_1}, to $0.01$, at configuration $C20$. The
geometries corresponding to configurations $C0$, $C10$ and $C20$ are shown in
Figure~\ref{fig:solution_example_1}.  For each configuration, the Darcy problem,
formulated as in \ref{pb:model_darcy_weak} or
\ref{pb:model_darcy_primal_weak}, depending on the method, is first solved in
order to compute the Darcy velocity $\bm{u}$, with a null source term. A unitary
hydraulic conductivity $K$ is used for all the fractures and pressure head boundary
conditions are prescribed on the bottom part of $\partial \Omega_l$ (inflow) and
$\partial \Omega_r$ (outflow) equal to $1$ and $0$, respectively (see
Figure~\ref{fig:geometry_example_1}). On the other portions of the boundary a
no-flux boundary condition is imposed. Subsequently, an advection-diffusion
problem is solved, obtained from Problem~\ref{pb:model_heat_weak} setting the
reaction operator $r(\cdot,\cdot)\equiv 0$, and with null source. We assume a
unitary coefficient $\zeta$, a diffusion coefficient $D$ equal to $10^{-4}$ and
a constant in time unitary Dirichlet boundary condition on the inflow part of
the domain boundary, whereas homogeneous Neumann boundary condition are set on
the rest of the boundary.

\begin{figure}
    \centering
    \includegraphics[width=0.33\textwidth]{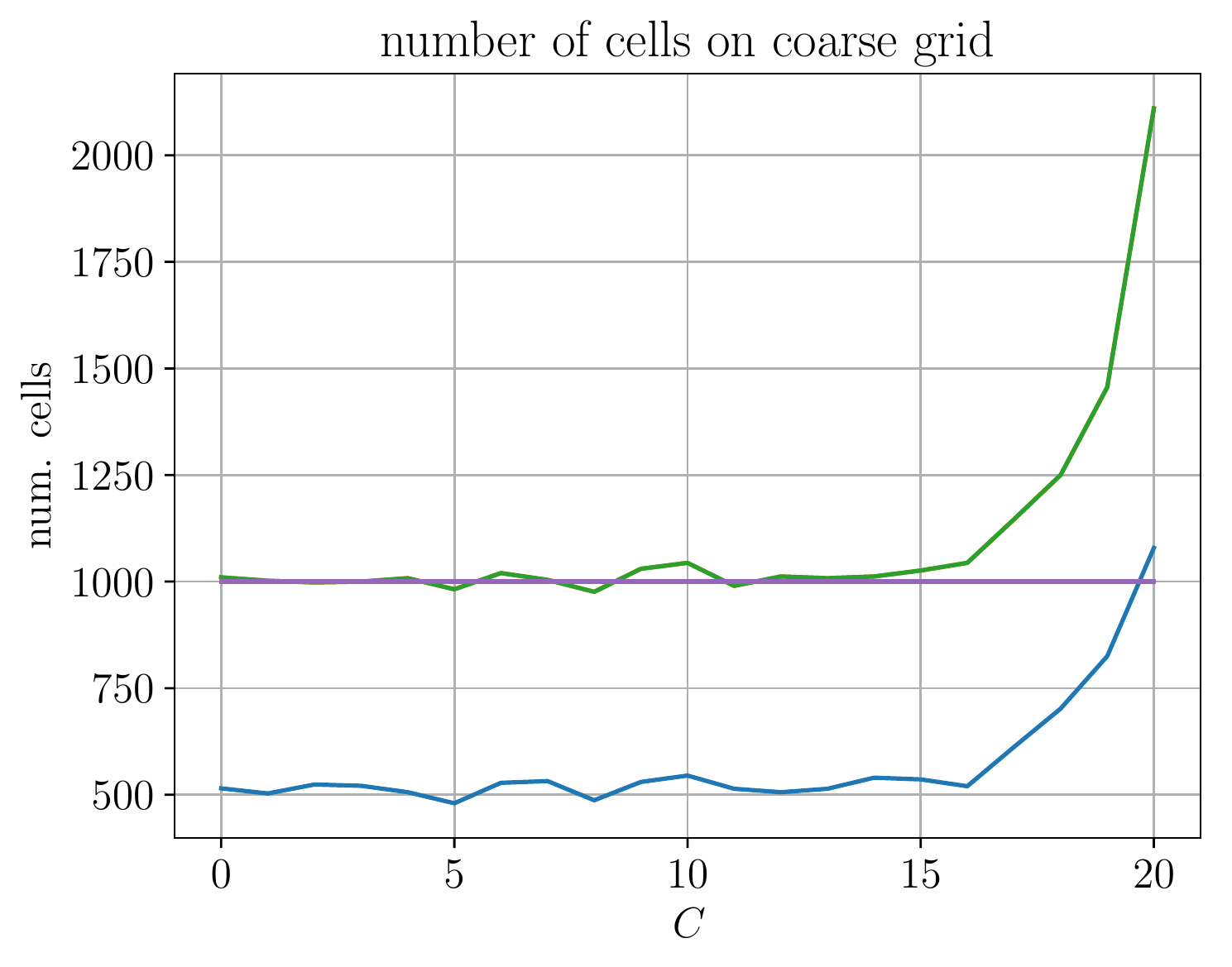}%
    \includegraphics[width=0.33\textwidth]{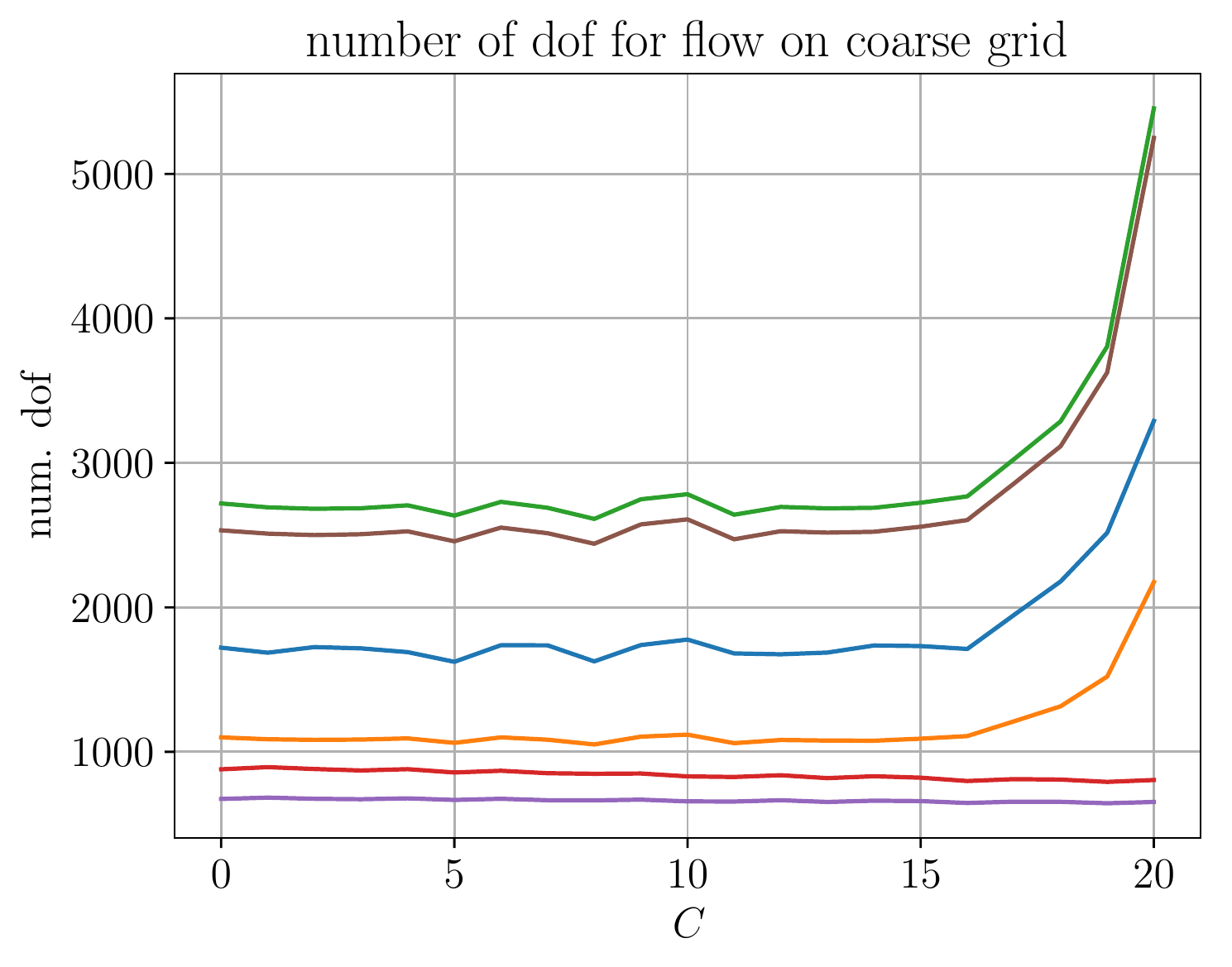}%
    \includegraphics[width=0.33\textwidth]{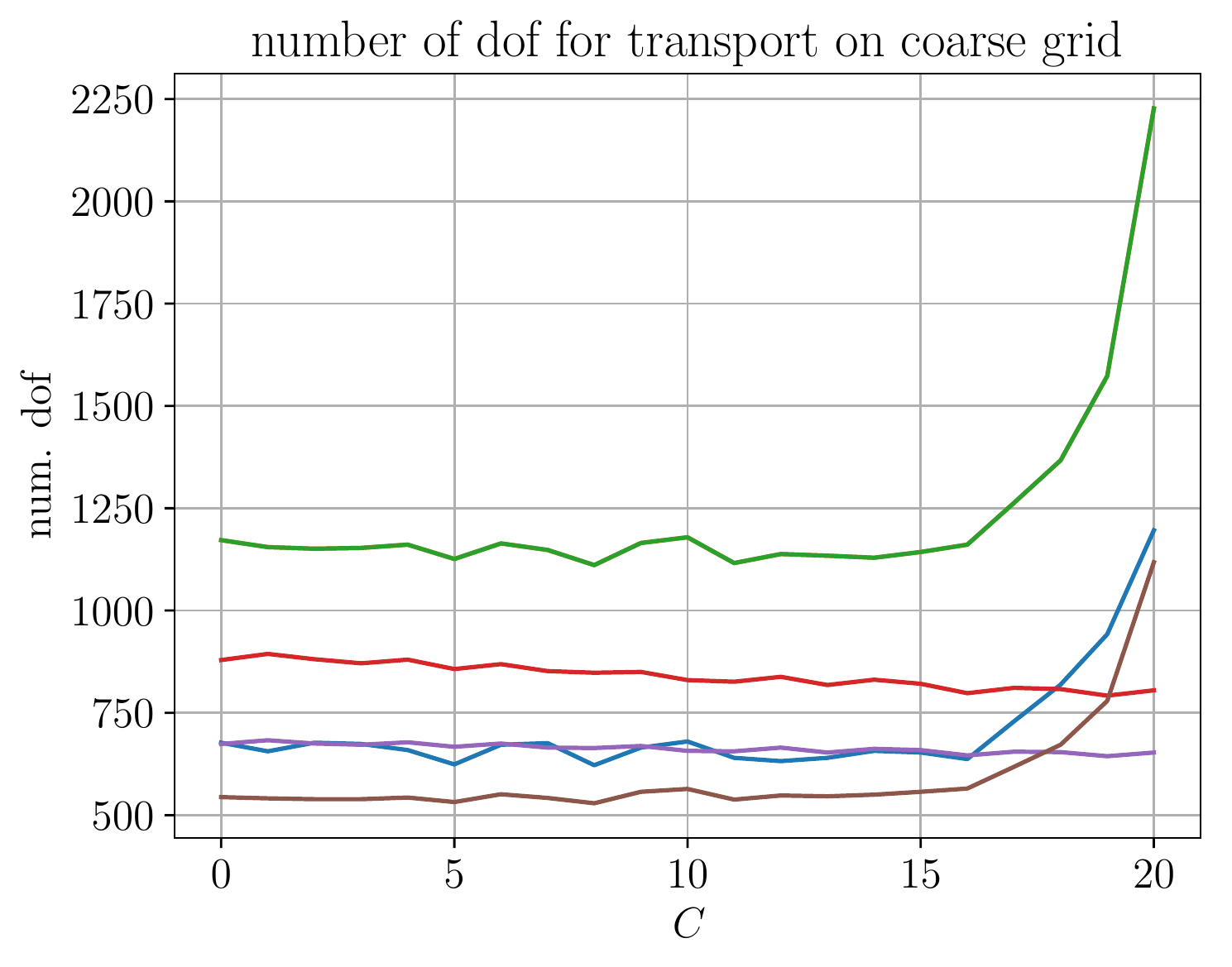}\\
    \includegraphics[width=0.32\textwidth]{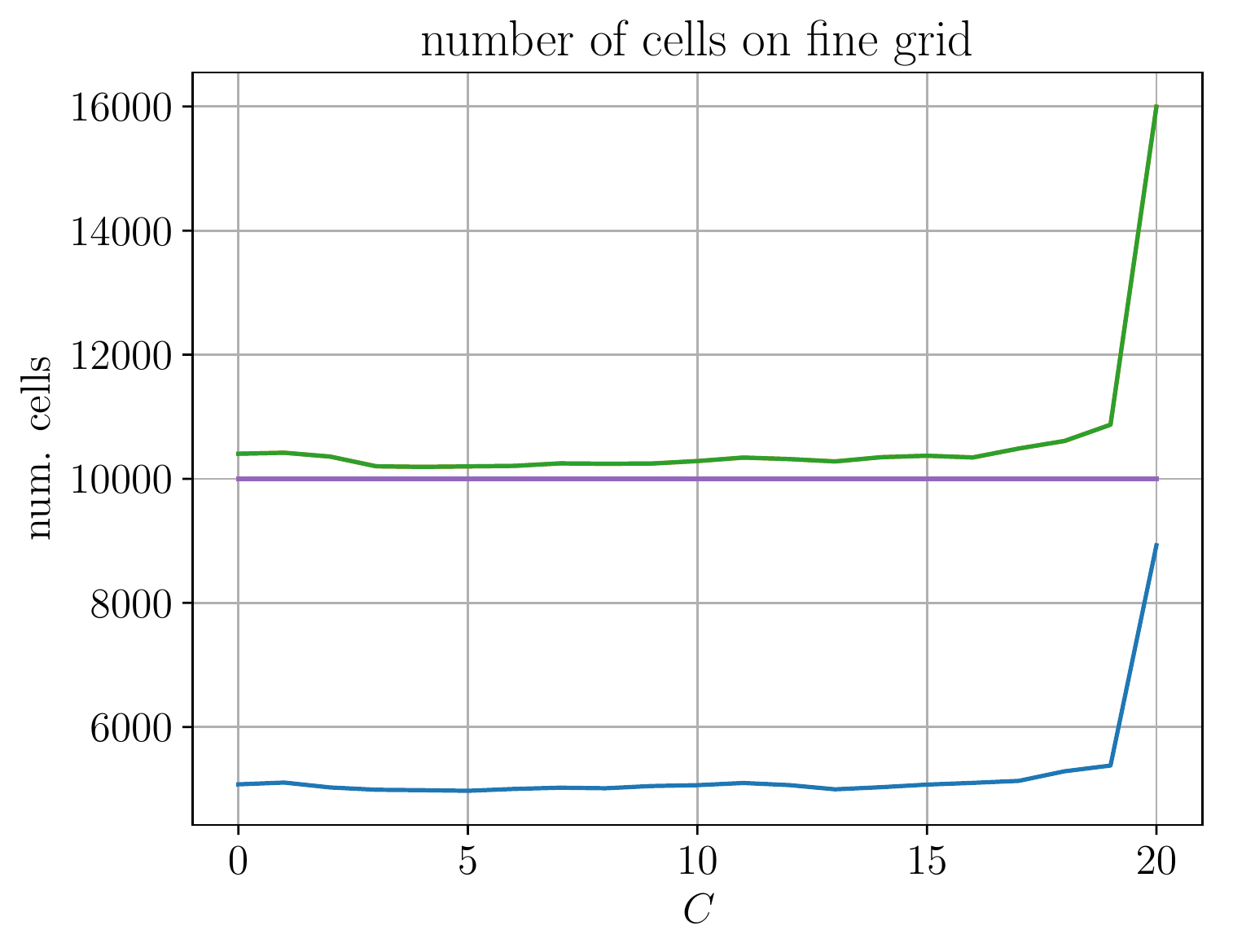}
    \includegraphics[width=0.33\textwidth]{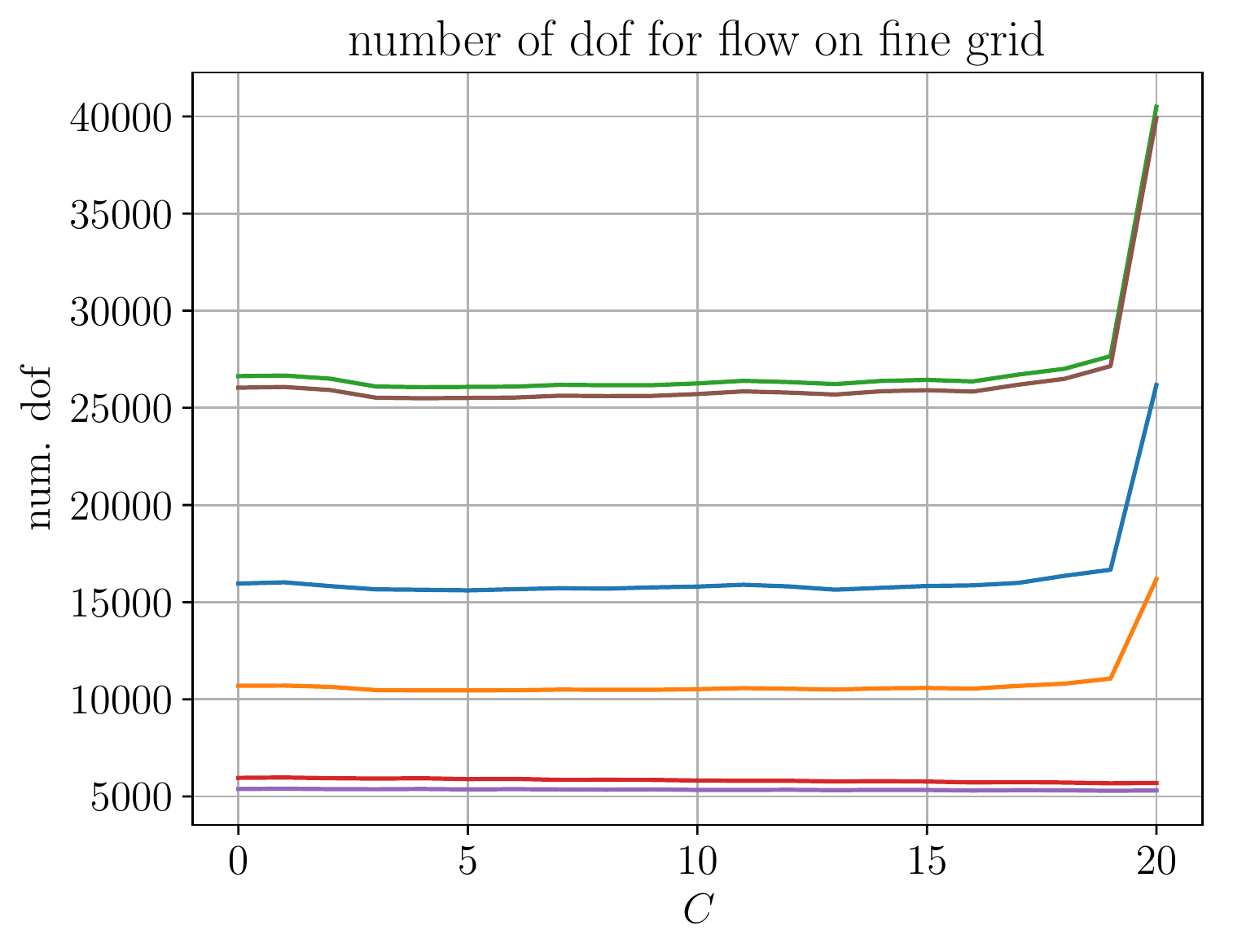}
    \includegraphics[width=0.33\textwidth]{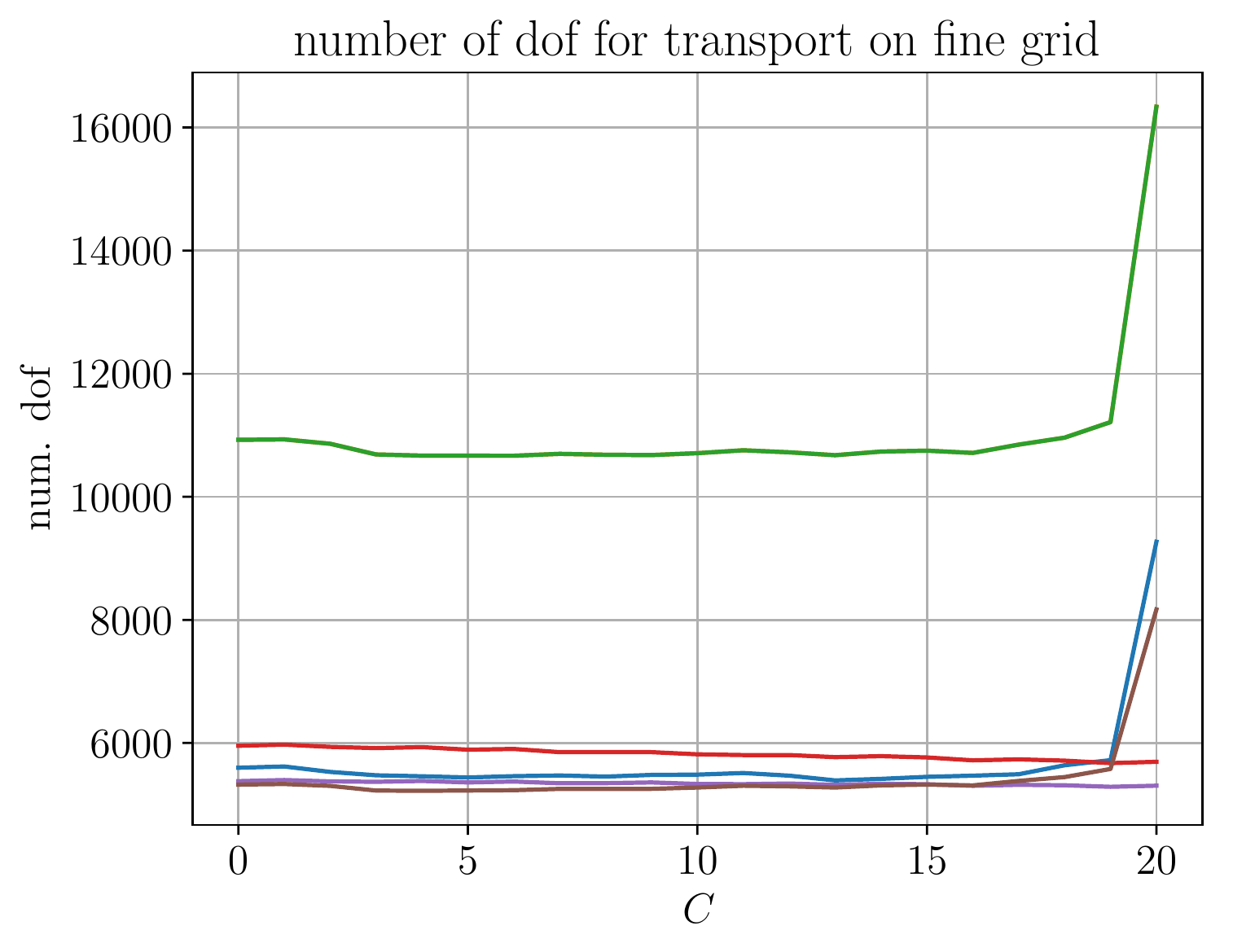}\\
    \includegraphics[width=\textwidth]{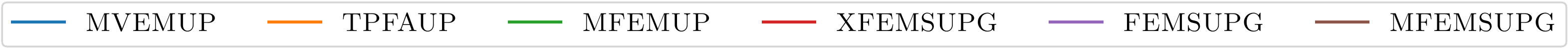}
    \caption{Number of cells (left), number of DOFs for the Darcy problem
    (middle) and number of DOFs for the dispersion problem (right) against
    configuration id for the two mesh parameters of Test Case 1, coarse mesh on
    the top and fine mesh on the bottom grids.}%
    \label{fig:example1_num_cells_dofs}
\end{figure}

The meshes used for the space discretization are different for conforming-triangular, conforming-polygonal and non-matching strategies, as discussed above. The
mesh for the conforming schemes are generated using the Gmsh
\cite{Geuzaine2009} library,
and two different mesh parameters are used, corresponding to about $10^3$ and
$10^4$ elements for the configuration $C0$. Clearly, as the length of trace
$\Gamma_0$ progressively reduces from configuration $C0$ to $C20$, the mesh
generation tool tends to increase the number of elements in order to meet the
conformity requirement without compromising mesh quality. This process
inevitably leads to a large increment of the number of elements from
configuration $C0$ to $C20$, for each of the three initial refinement levels.
The mesh for the non-matching schemes is obtained through the Triangle
\cite{Shewchuk1996} software, using two mesh parameters, again corresponding to
about $10^3$, and $10^4$ elements. In this case, instead, since the mesh is
independent from the traces, the number of elements is practically unaffected by
the change of the geometry from configurations $C0$ to $C20$ (small oscillations
of few elements are observed as a consequence of the change of the coordinates
of fracture $\Omega_l$). The polygonal mesh for the \MVEM{} approach is finally
built gluing together the triangular elements of the conforming mesh, thus
aiming at mitigating the increase of mesh elements and degrees of freedom. The
number of elements for the various schemes, for each geometrical configuration
is reported in Figure~\ref{fig:example1_num_cells_dofs}, on the leftmost column,
for the coarse (top) and fine (bottom) meshes. The centre and right columns of
Figure~\ref{fig:example1_num_cells_dofs} report instead the number of degrees of
freedom corresponding to each method, which can be used as an indication of the
corresponding computational cost. We remark, however, that the
computational cost is also largely affected by the approach used to solve linear
systems: some of the proposed methods, e.g., have efficient parallel solvers,
other take advantage of standard domain decomposition strategies to achieve
computational efficiency. Also, the availability of effective preconditioners
should be taken into account. Here, however, the main focus is on the analysis
of the response of the methods in terms of prediction accuracy versus
geometrical and model complexities typical of DFN simulations, thus we refer to
the literature of each method for further details on computational efficiency
issues.

In the following of this test case, for brevity, the mesh with $10^3$ cells on configuration $C0$
will be denoted as coarse mesh and the mesh with $10^4$ cells on $C0$ as the
fine. For time discretization an equally spaced mesh with time-step equal to
$0.05$ is considered and $300$ time-steps are performed, starting from an all zero initial condition.

\begin{figure}
\centering
\includegraphics[width=0.99\textwidth]{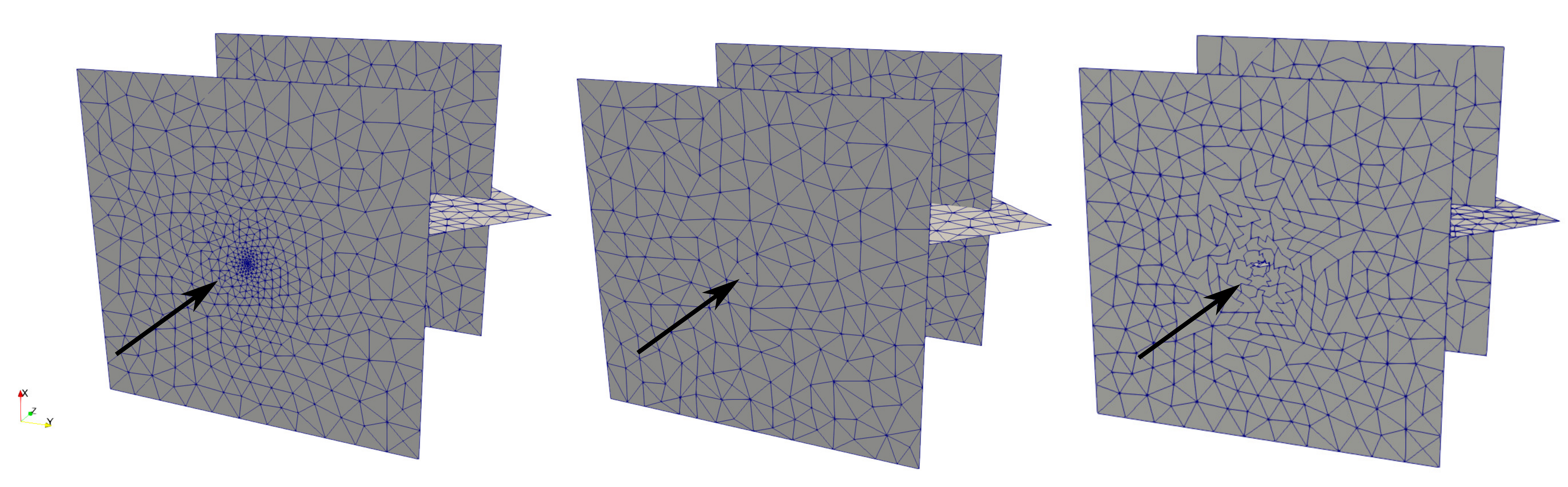}
\caption{Mesh for Test Case 1 on configuration $C20$ for conforming (left) non-conforming (middle) and coarsened (right) schemes}
\label{fig:meshes_example_1}
\end{figure}

The solution obtained with the scheme \MFEM{} is reported, as an example, in
Figure~\ref{fig:solution_example_1} on configurations $C0$, $C10$ and $C20$. On
the leftmost column the pressure head distribution in the network is shown for the
three geometries, whereas the remaining columns depict the temperature
distribution $\theta$ at three time-steps corresponding to a time $t=1.25$,
$t=2.5$ and $t=5$, respectively. We can notice that, as the trace
between fractures $\Omega_l$ and $\Omega_r$ vanishes, the pressure head distribution
displays a steeper gradient around the intersection, and the effective
permeability of the network reduces, thus also reducing the penetration depth of
the higher temperature zone in the network at fixed time frames. From a
computational point of view, simulations become more and more challenging as the
solution starts to display steep gradients, especially for methods built on
non-conforming meshes, that are not adapted to the geometry, as shown in
Figure~\ref{fig:meshes_example_1}.


For all the proposed numerical schemes, at each time-step, the following
quantities are computed: the average temperature $\langle{\theta} \rangle_\Omega$, the
average temperature $\langle {\theta}\rangle_{\mathrm{outflow}}$ on the outflow portion of the boundary, and,
for non-locally conservative schemes, as the optimization-based methods, the
total flux mismatch $\delta \Phi_\Gamma$ at each trace.

A reference solution is computed using the \MFEM{} method on a mesh much finer
than the three considered here for the simulations, counting about
$6.3\times10^5$ cells, almost independently from the configuration id, as cell
size is capable of resolving the smallest length of the vanishing trace.

\begin{figure}
    \centering
    \includegraphics[width=0.33\textwidth]{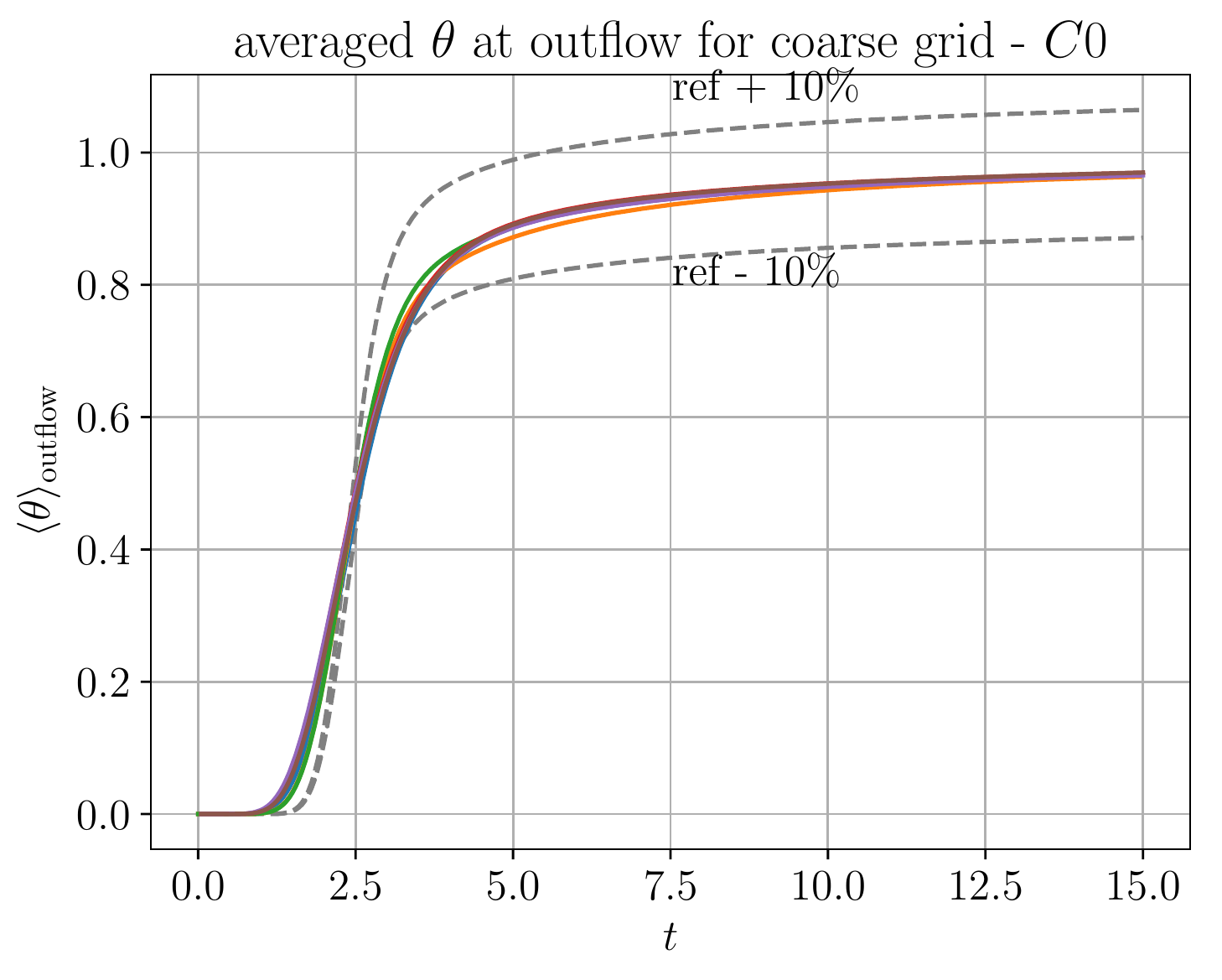}%
    \includegraphics[width=0.33\textwidth]{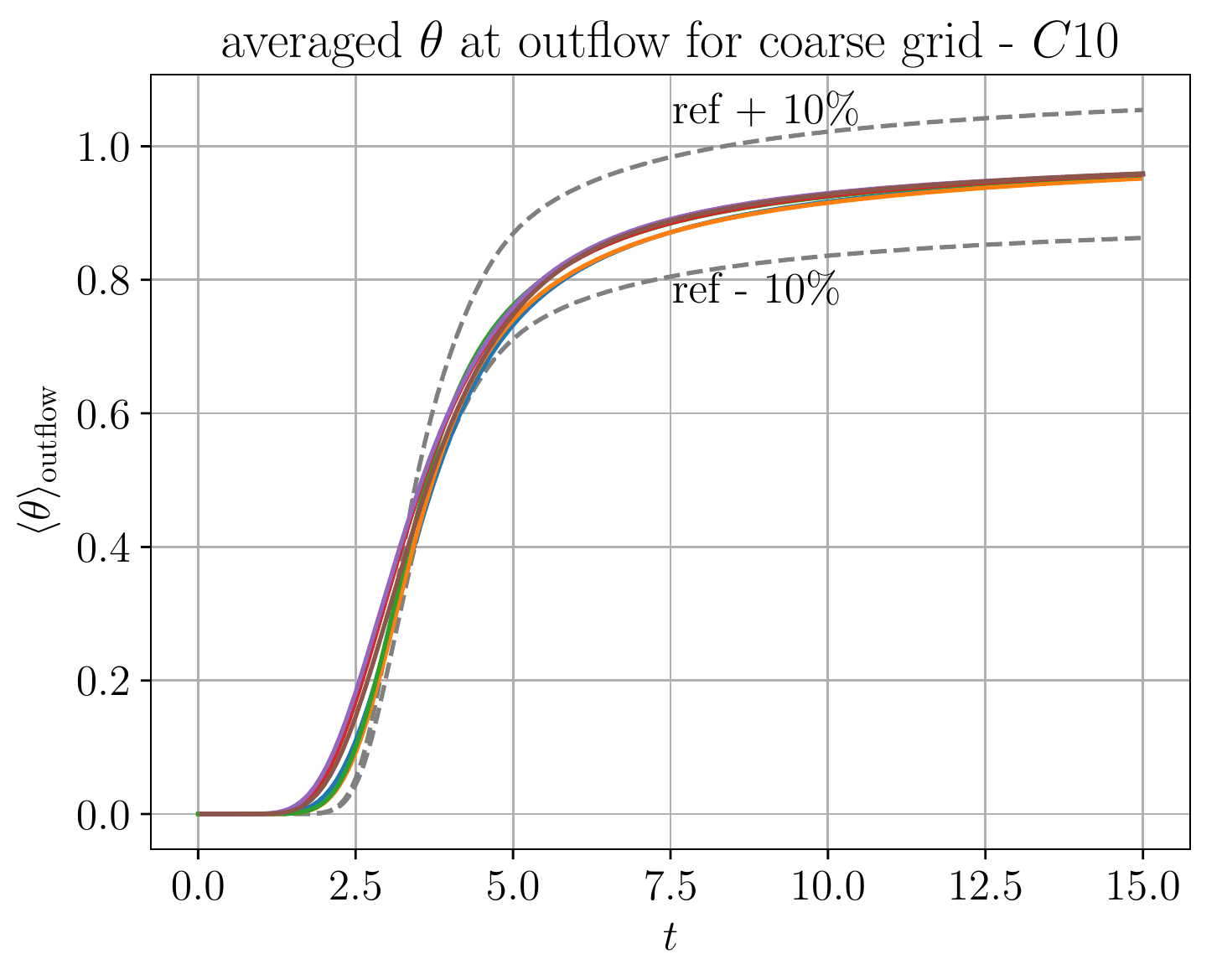}%
    \includegraphics[width=0.33\textwidth]{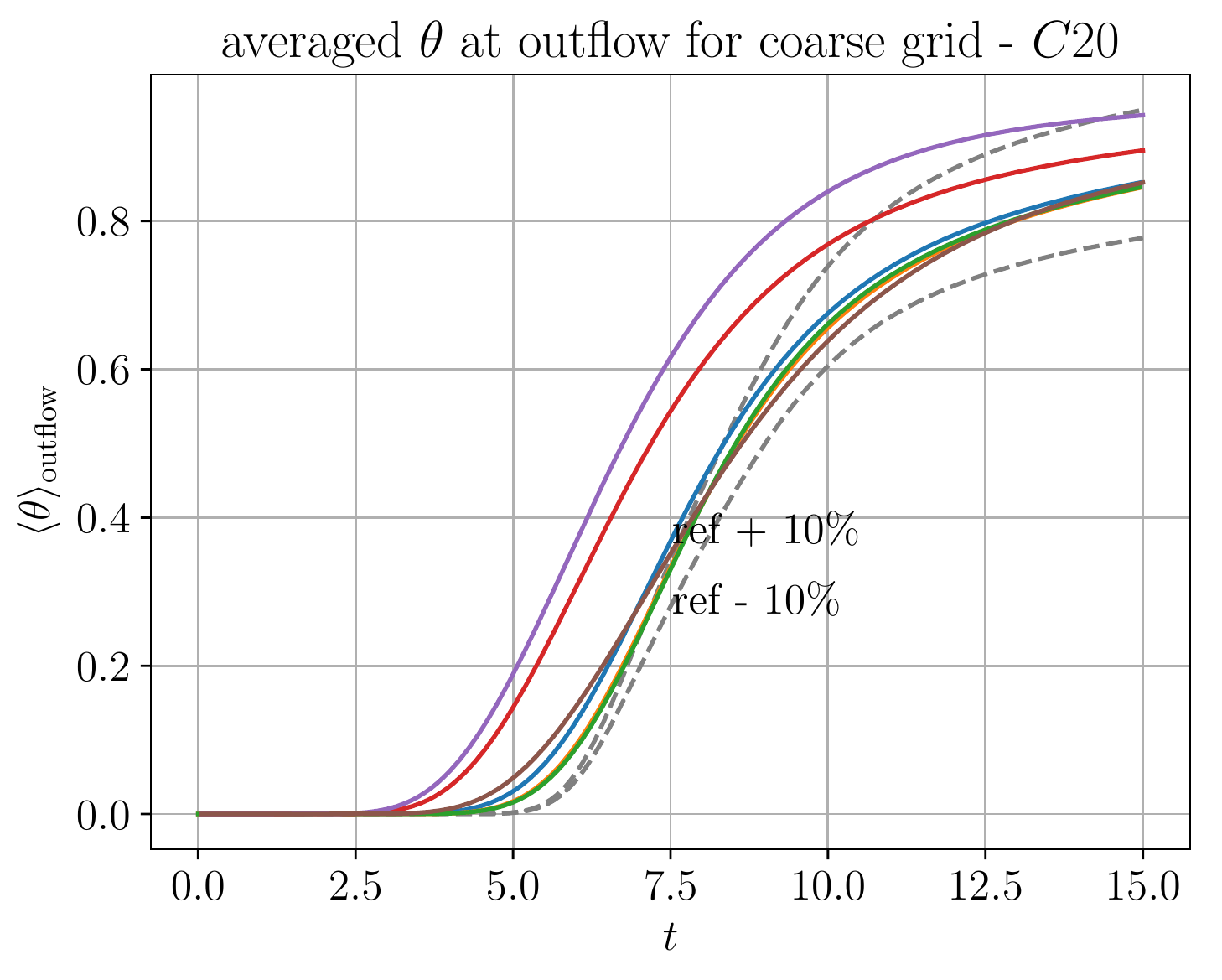}\\
    \includegraphics[width=0.33\textwidth]{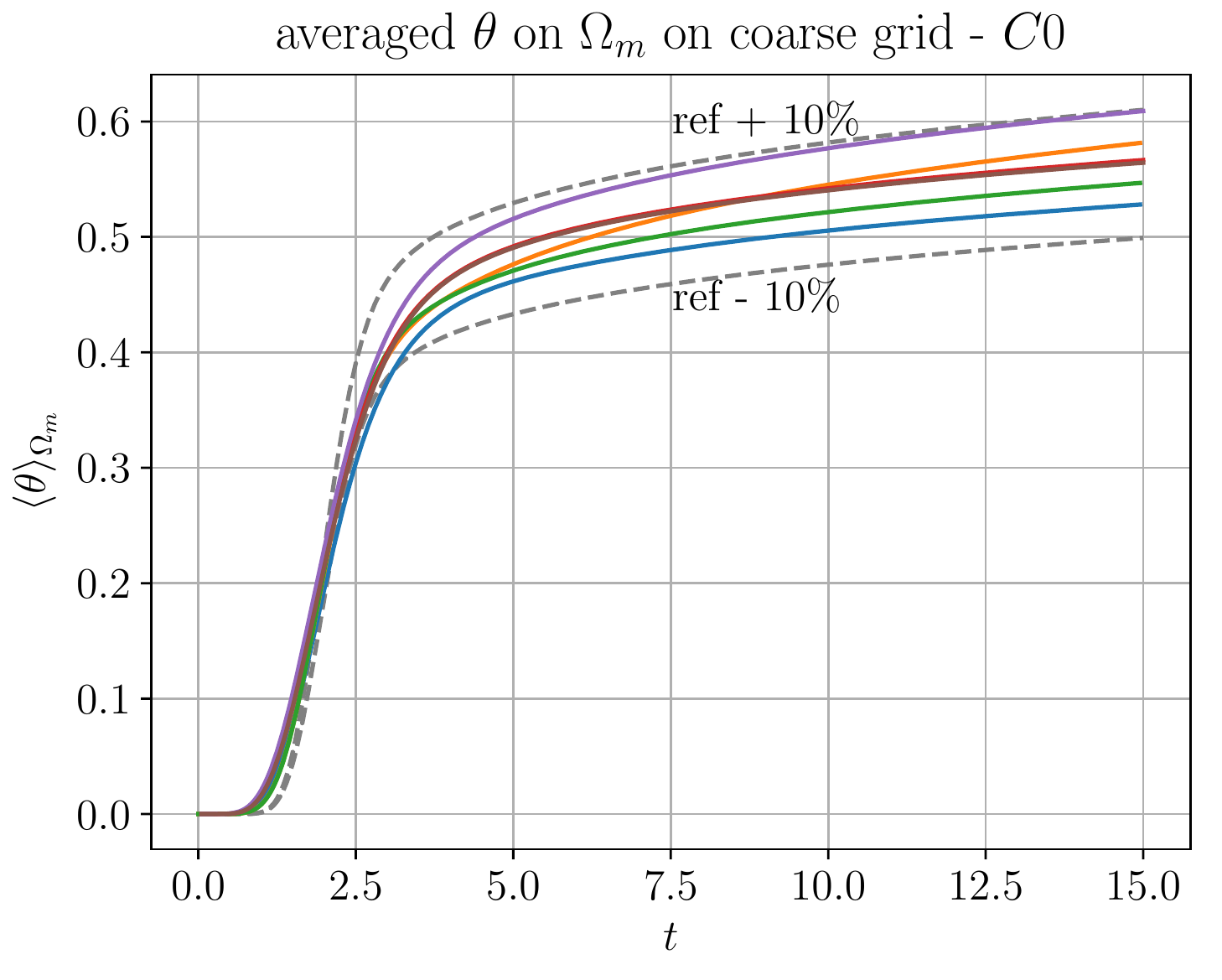}%
    \includegraphics[width=0.33\textwidth]{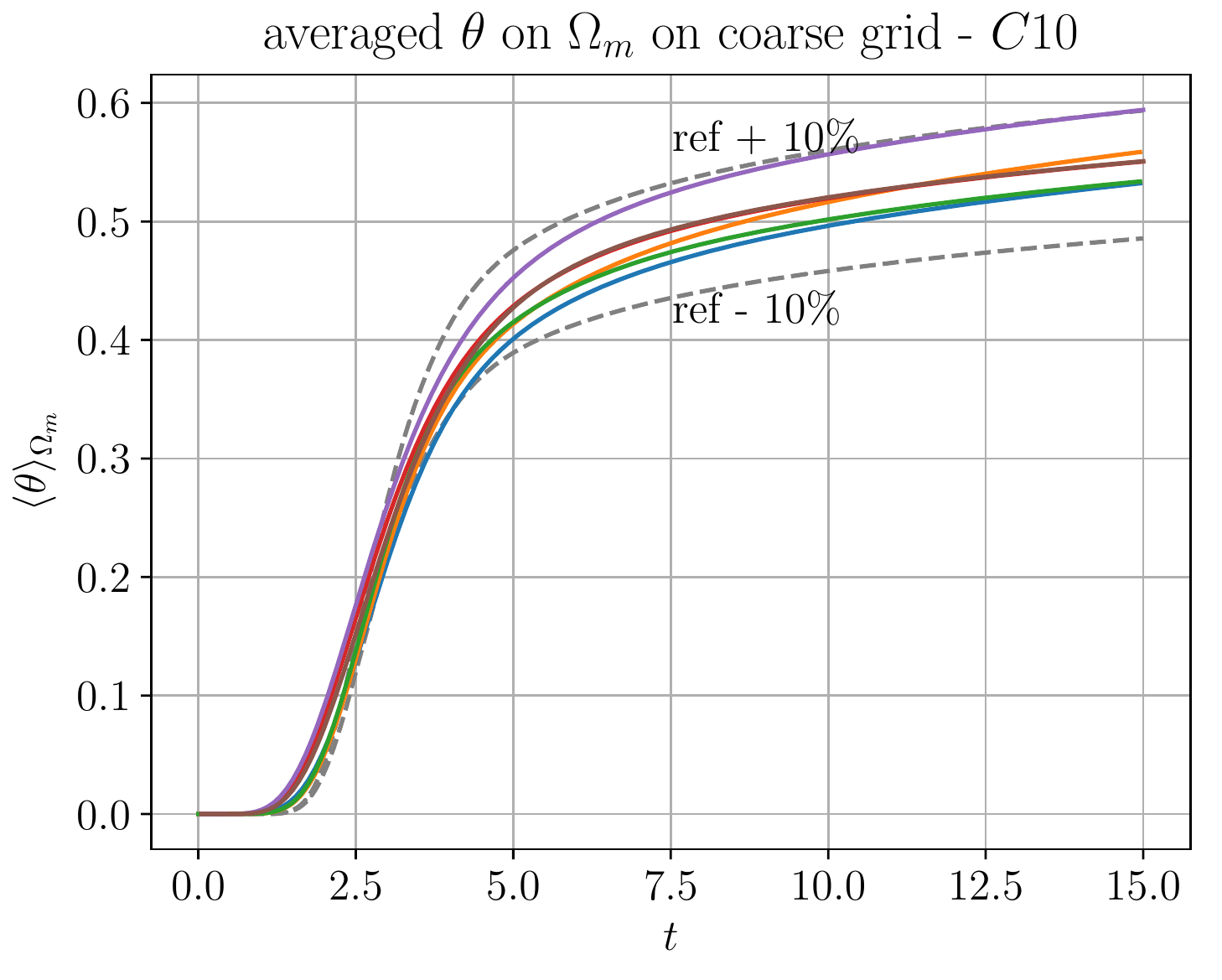}%
    \includegraphics[width=0.33\textwidth]{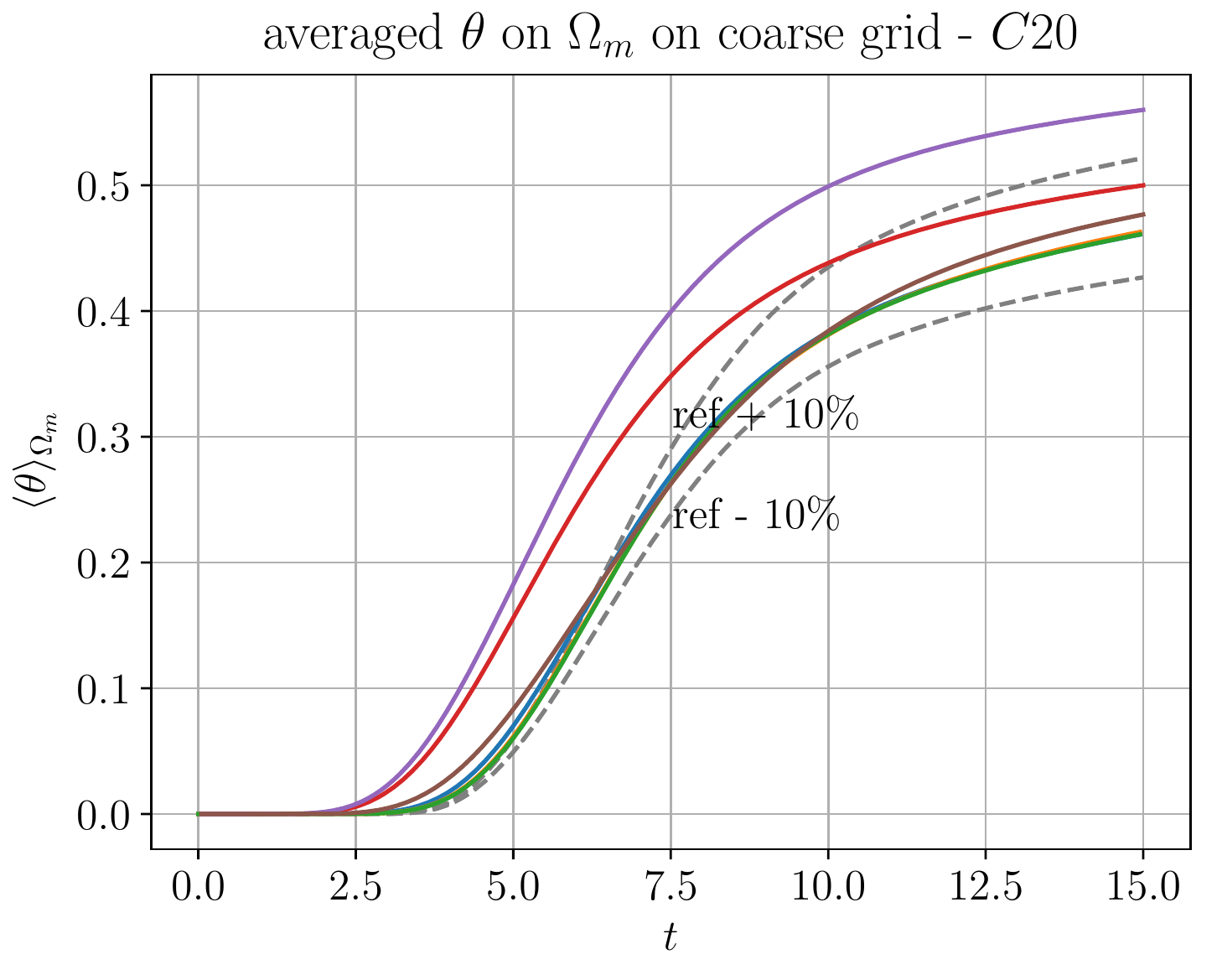}\\
    \includegraphics[width=0.33\textwidth]{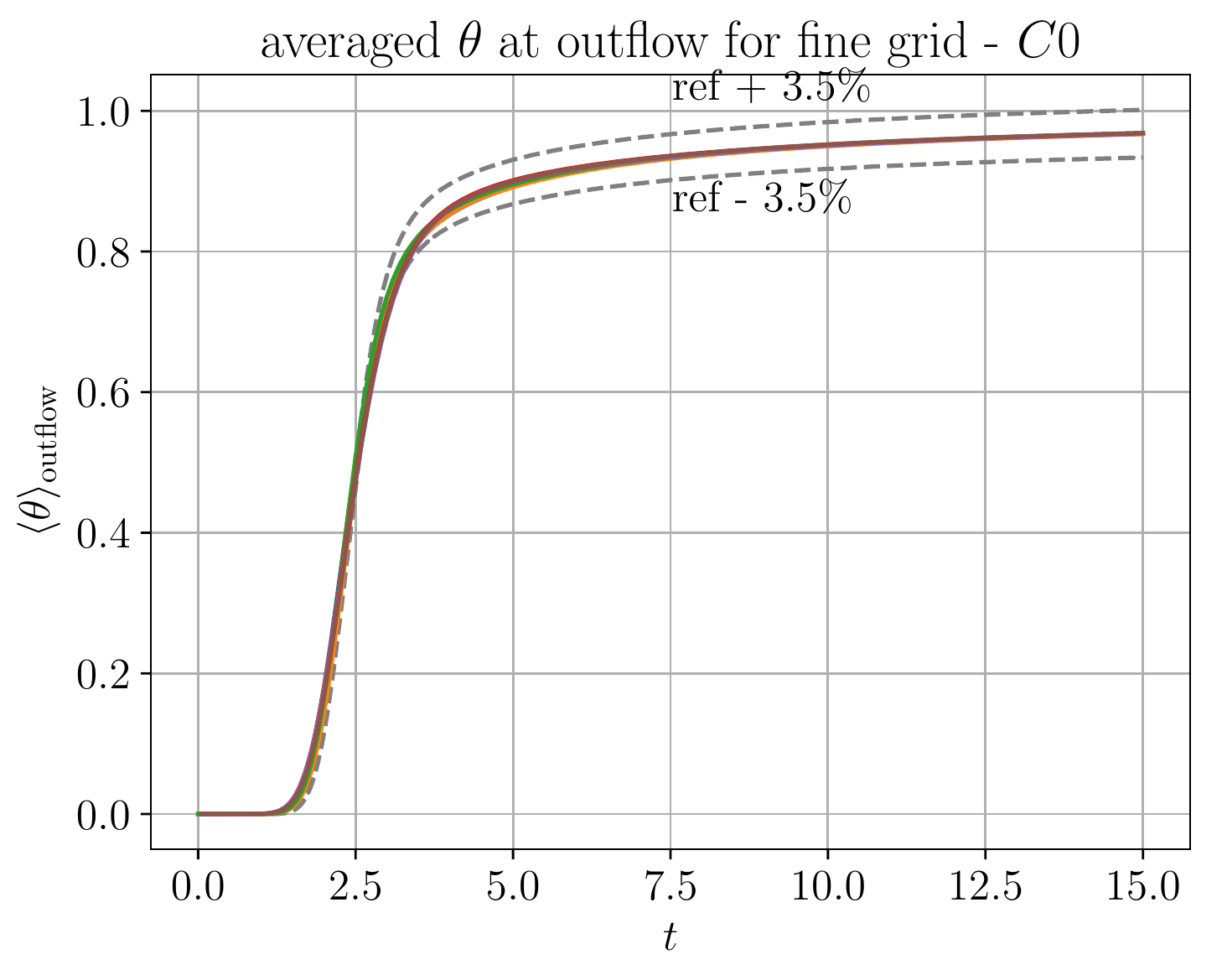}%
    \includegraphics[width=0.33\textwidth]{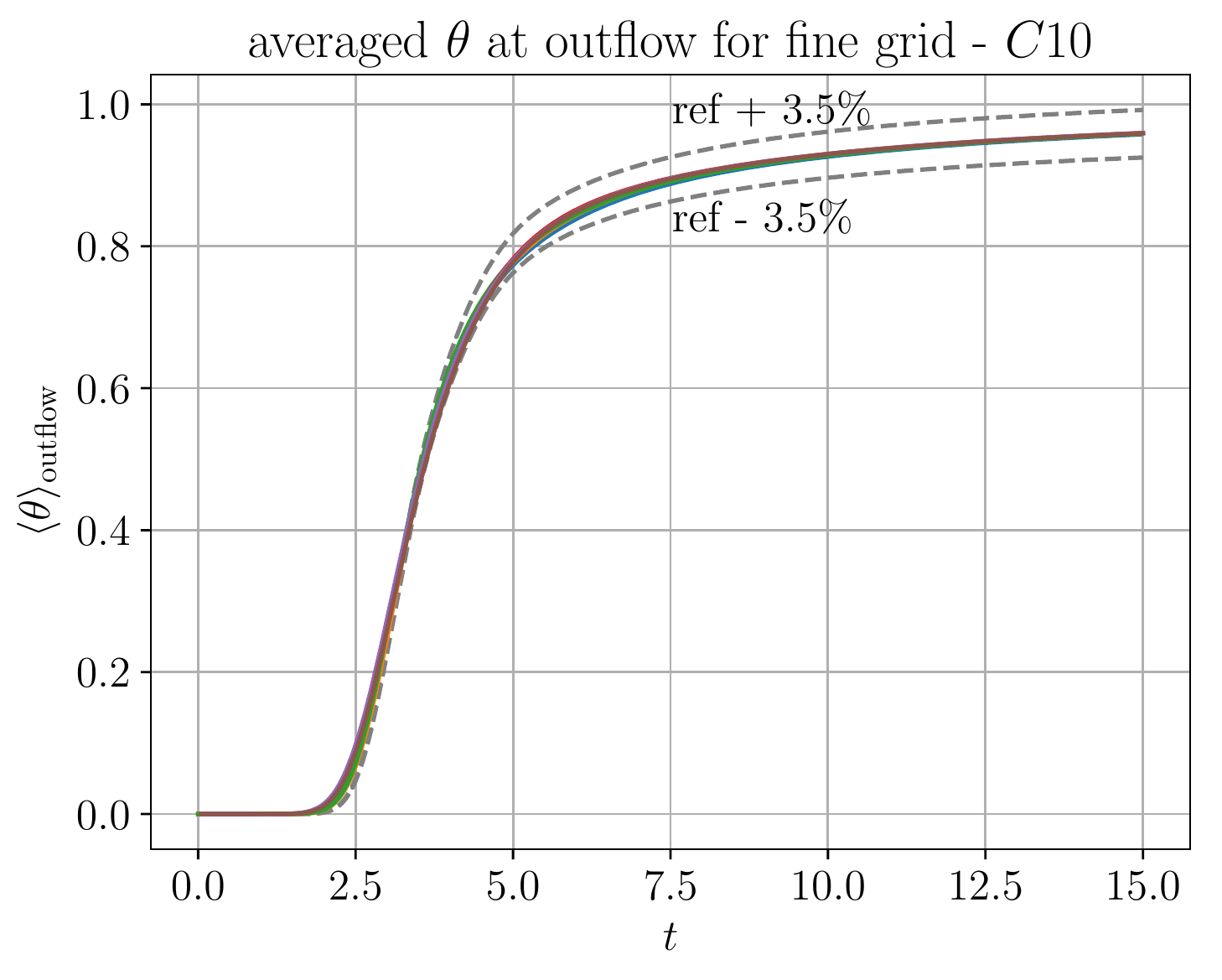}%
    \includegraphics[width=0.33\textwidth]{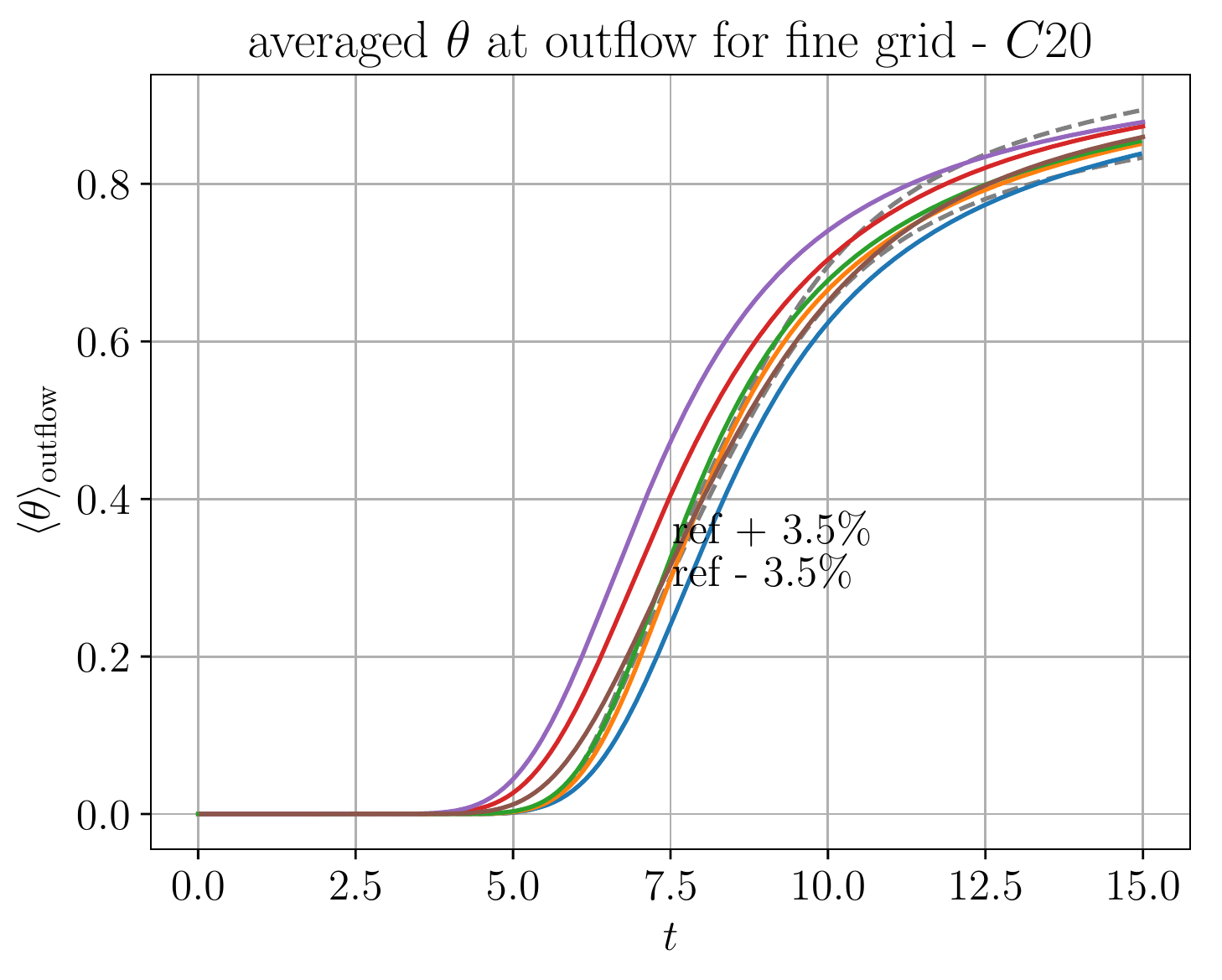}\\
    \includegraphics[width=0.33\textwidth]{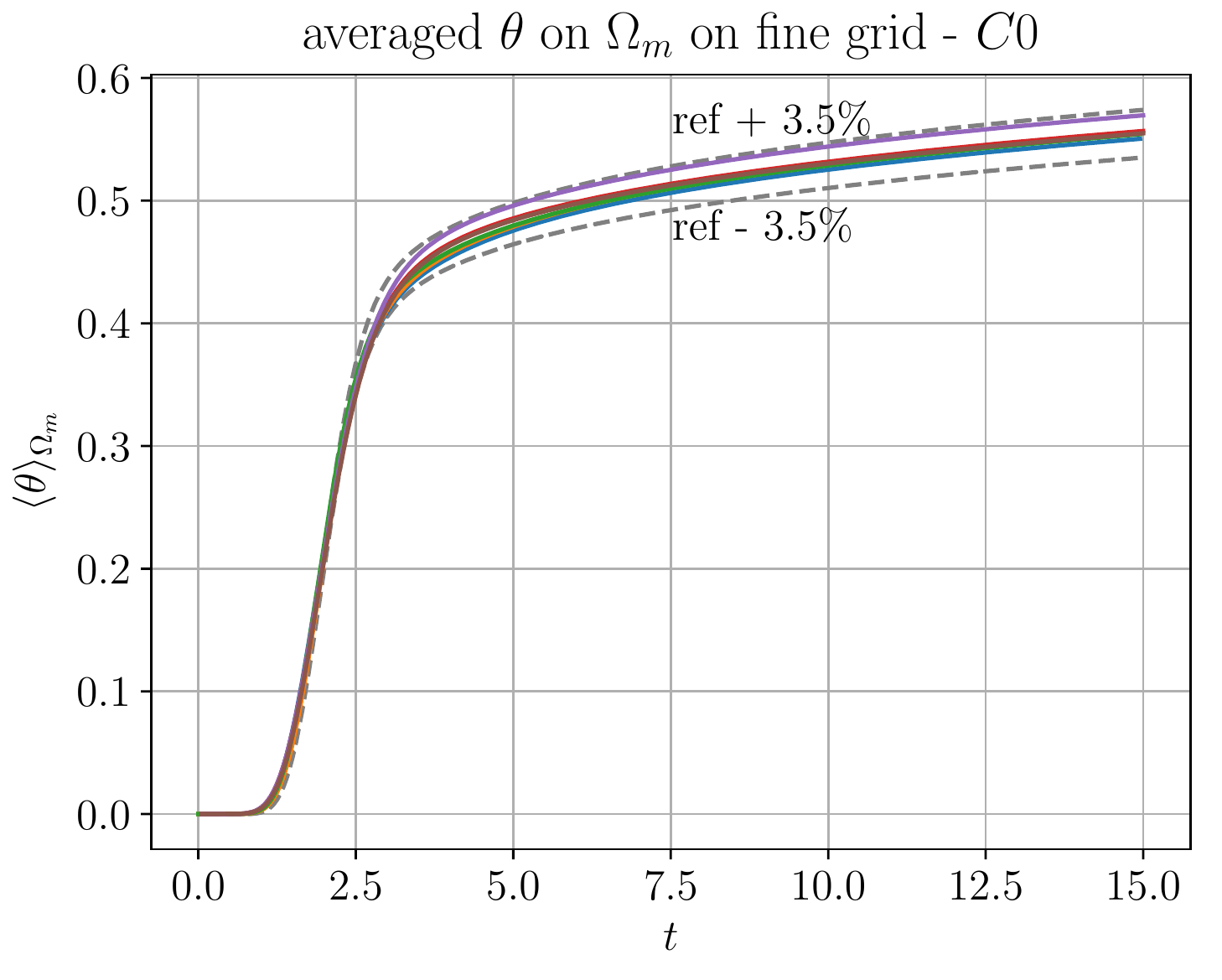}%
    \includegraphics[width=0.33\textwidth]{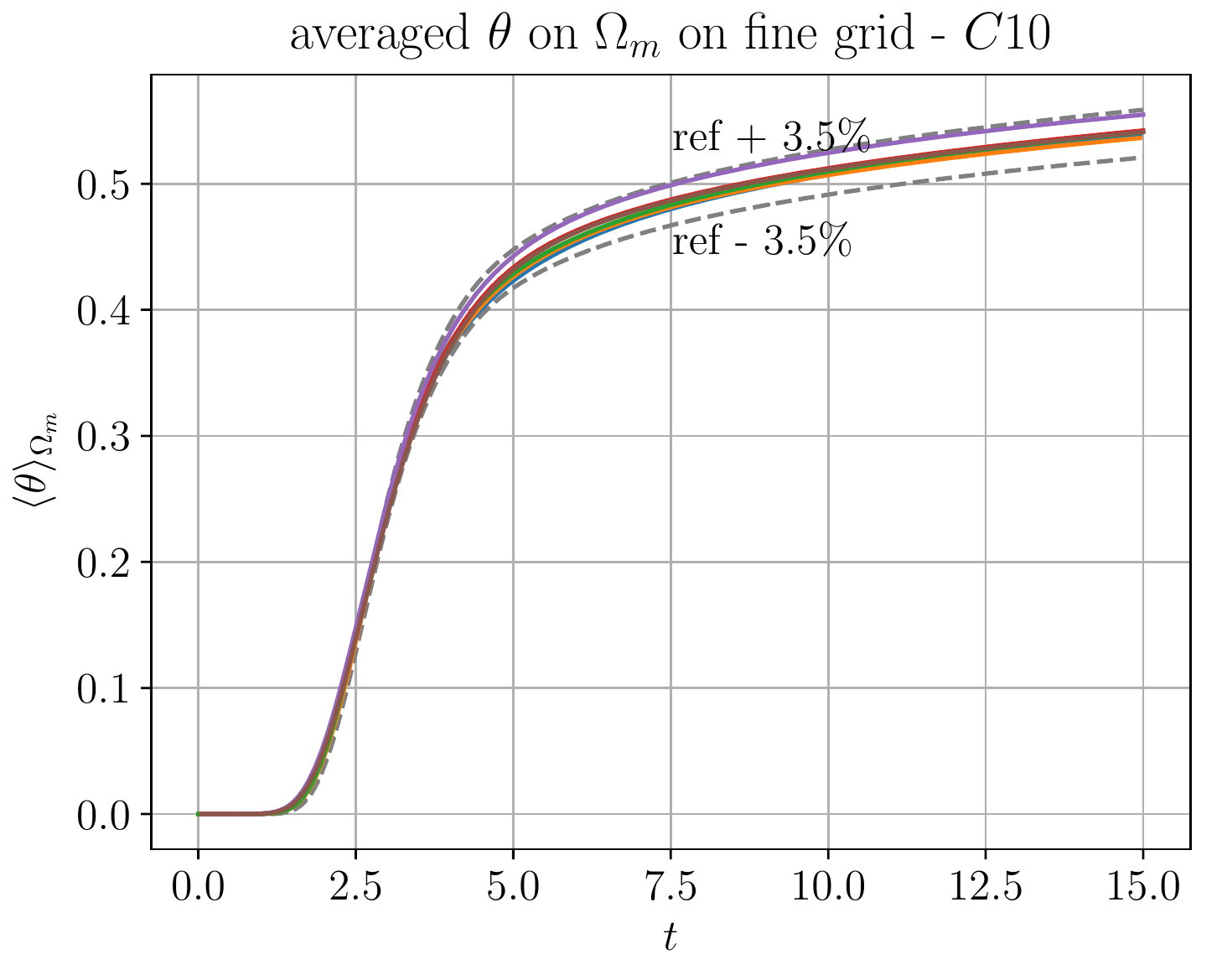}%
    \includegraphics[width=0.33\textwidth]{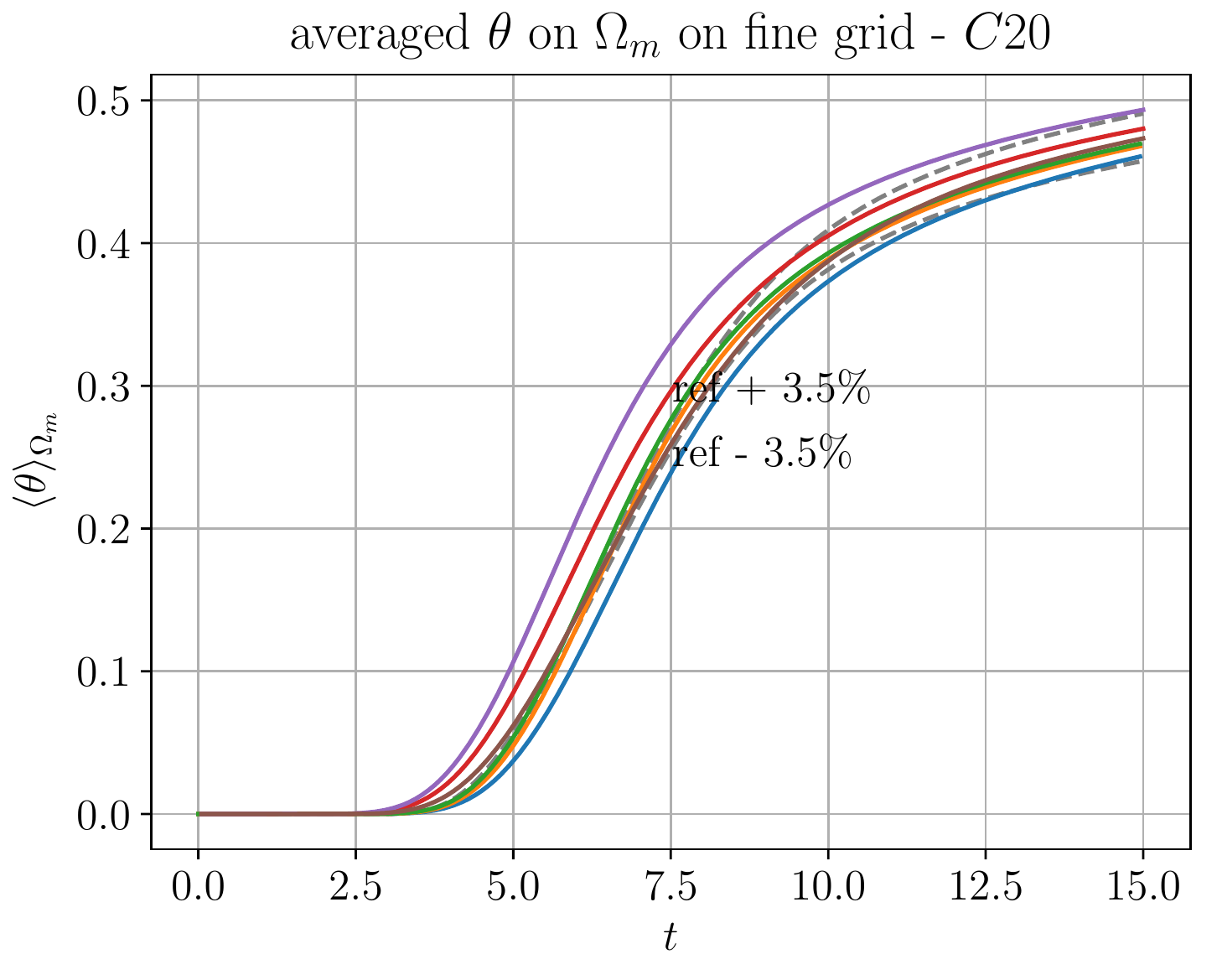}\\
    \includegraphics[width=\textwidth]{label}
    \caption{Average temperature on curves against time for Test
    Case 1 on the coarsest (first six pictures) and finest (last six pictures)
    for selected fracture and configurations. Columns refer to the same
    geometrical configuration: $C0$ on column 1, $C10$ on column 2 and $C20$ on
    column 3. The average temperature curves refer to fracture $\Omega_m$.}%
    \label{fig:example1_plots}
\end{figure}

The plots in Figure~\ref{fig:example1_plots} propose a comparison of average
temperature for all the proposed methods against time, on
three selected configuration, for all the considered meshes. In the picture,
dashed black lines represent relative errors with respect to the reference
solution. As expected, all the curves are in very good agreement, also on the
coarsest mesh, for configuration $C0$, whereas small discrepancies appear for
configuration $C10$ on the coarse mesh that however disappear as the mesh is
refined. Larger differences appear, for methods \FEMSUPG{} and \XFEMSUPG{},
instead for the simulations on configuration $C20$. This is expected,
since, as mentioned, when the varying trace becomes very small, methods built on
non-matching meshes can not rely on the effect of mesh refinement around
the vanishing trace which clearly improves representation capabilities of conforming methods.
We observe that method \MVEMUP{} retains good approximation capabilities despite the coarsening.

\begin{figure}
\centering
\includegraphics[width=0.33\textwidth]{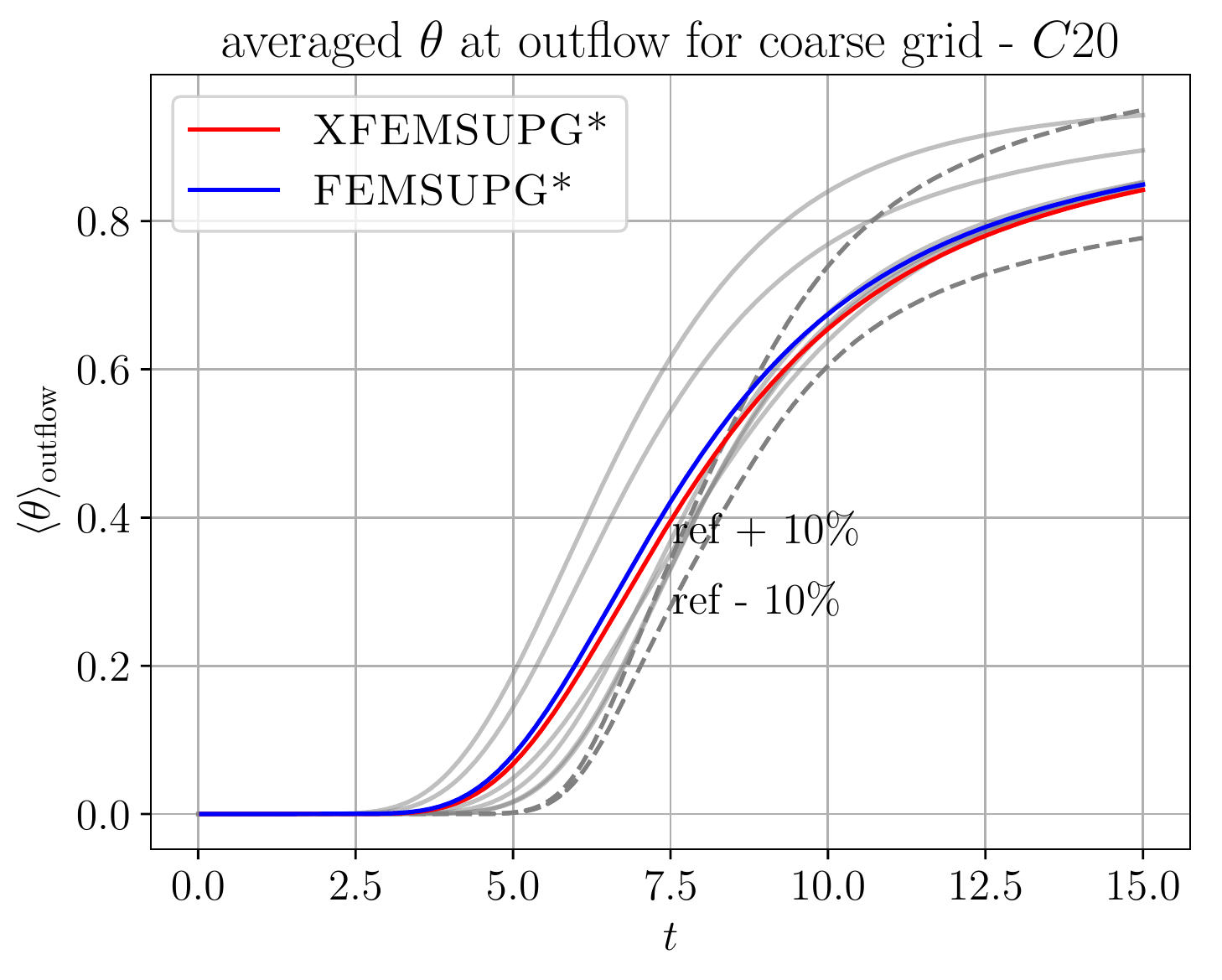}
\includegraphics[width=0.66\textwidth]{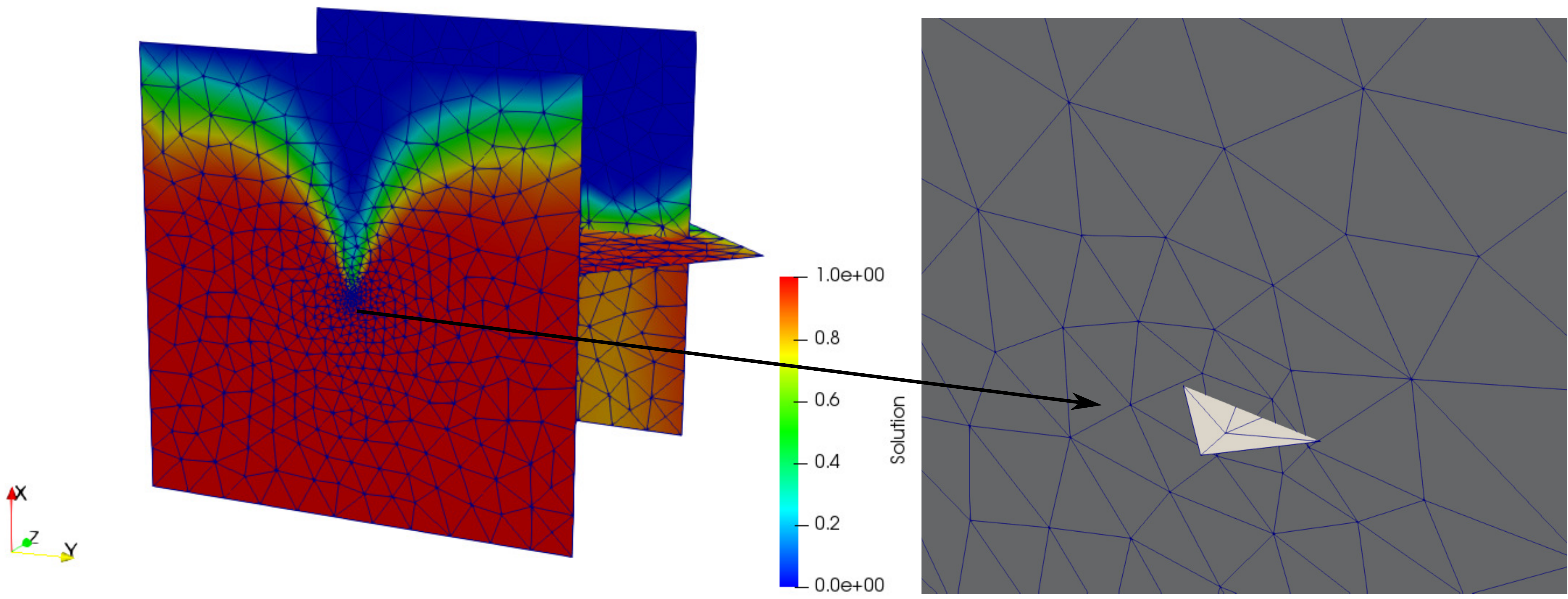}
\caption{Outflow average temperature curve with \OPTXFEM{} and \OPTFEM{} approaches (left), and
solution with \OPTXFEM{} at $t=15$ with a zoom of the mesh around $\Gamma_0$
(right), both on on the \textit{perturbed} mesh of configuration $C20$.}
\label{fig:example1_perturbed_mesh_results}
\end{figure}

In order to quantify what is the effect of mesh adaptation due to conformity
requirements on the quality of the solution with respect to the effect of the
approximations introduced by the non-matching schemes themselves, in
Figure~\ref{fig:example1_perturbed_mesh_results} (left) the outflow average
temperature for
the methods \OPTXFEM{} and \OPTFEM{} is reported against time, for configuration
$C20$ on a \textit{perturbation} of the coarse mesh used for the conforming
methods for configuration $C20$, overlapped to the curves of the other schemes,
in transparency, on the original meshes. We observe that the curves of
methods \OPTFEM{} and \OPTXFEM{} on the perturbed adapted mesh are now much
closer to the curves of other conforming approaches, thus clearly showing the
effect of mesh adaptation on the quality of the solution. The solution obtained
with the \OPTXFEM{} approach on the perturbed mesh at $t=15$ is reported in
Figure~\ref{fig:example1_perturbed_mesh_results}, on the right, along with a
detail of the mesh around $\Gamma_0$, clearly showing the non conformity
of the mesh. Two aspects are to be remarked: first, the non-conforming
approaches are capable of producing reasonable approximations of the solution
also on the coarse mesh, which can be greatly improved refining the mesh and
still using a fraction of the degrees of freedom required by the conforming
approaches, see, e.g., the last two plots in the third column of
Figure~\ref{fig:example1_plots}; second, non-conforming approaches allow to
freely choose the refinement level of the mesh, thus allowing to efficiently use
mesh adaptation strategies, only refining the mesh where required, independently
of the geometrical constraints \cite{BBV2019}.

\begin{figure}
    \centering
    \includegraphics[width=0.33\textwidth]{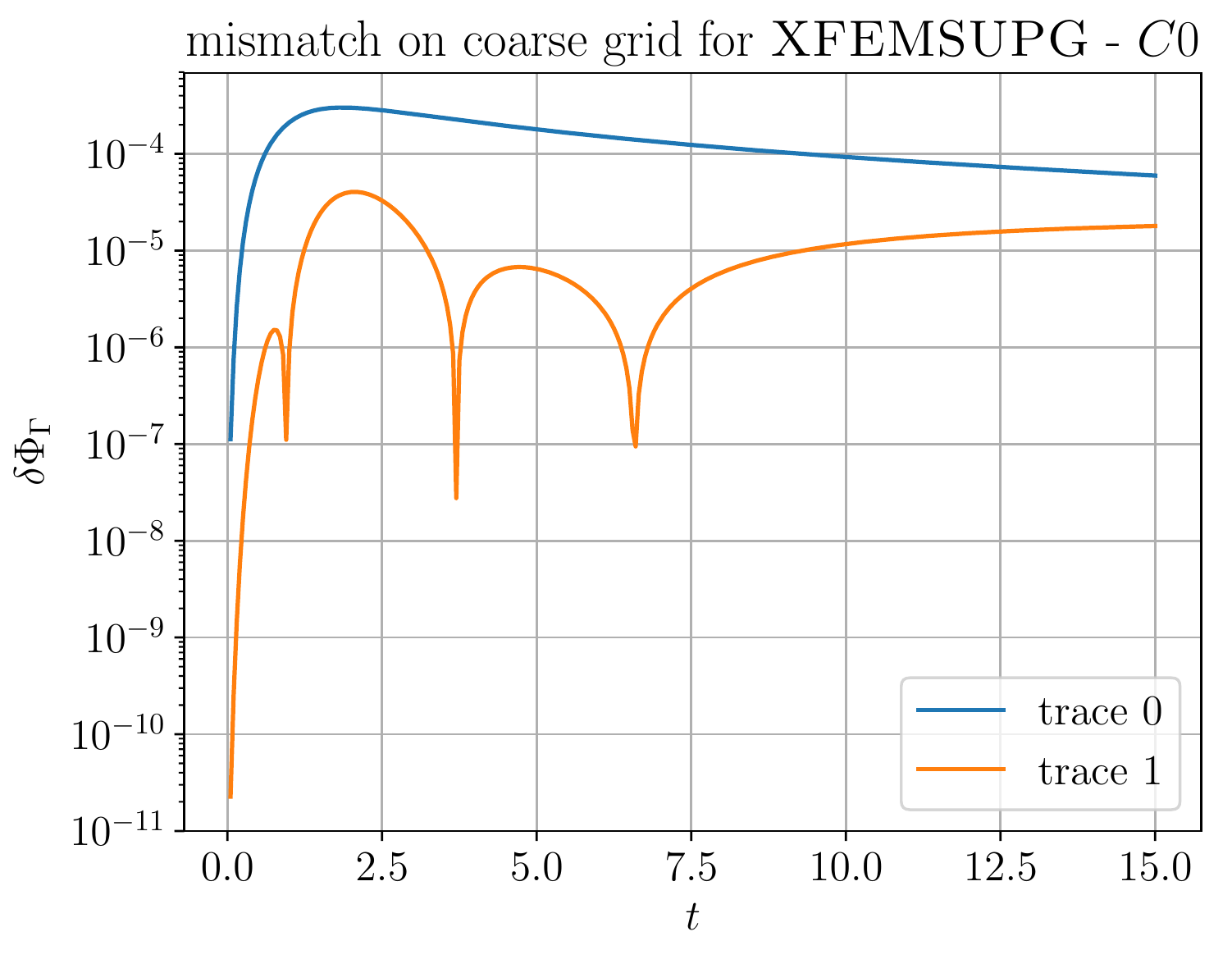}%
    \includegraphics[width=0.33\textwidth]{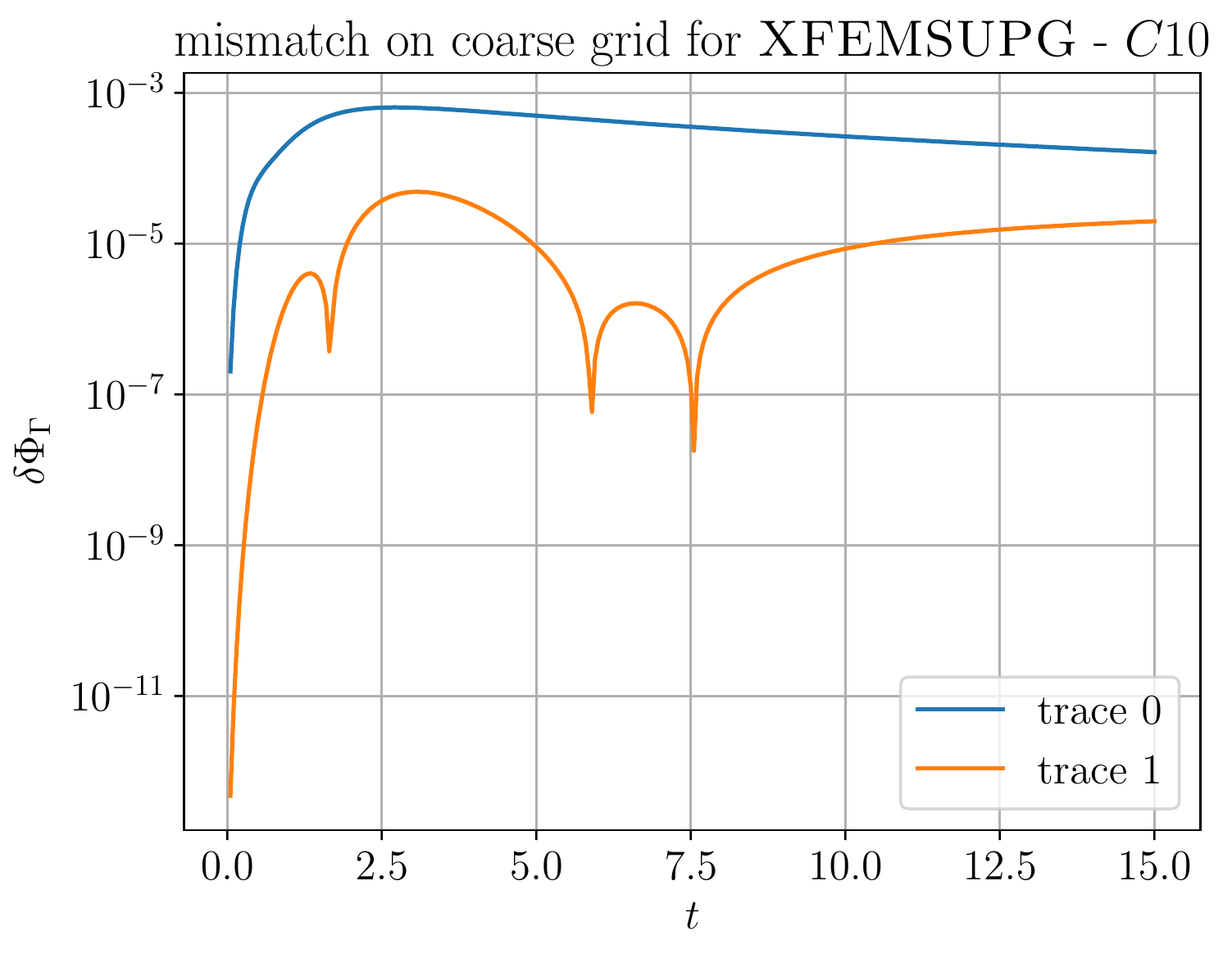}%
    \includegraphics[width=0.33\textwidth]{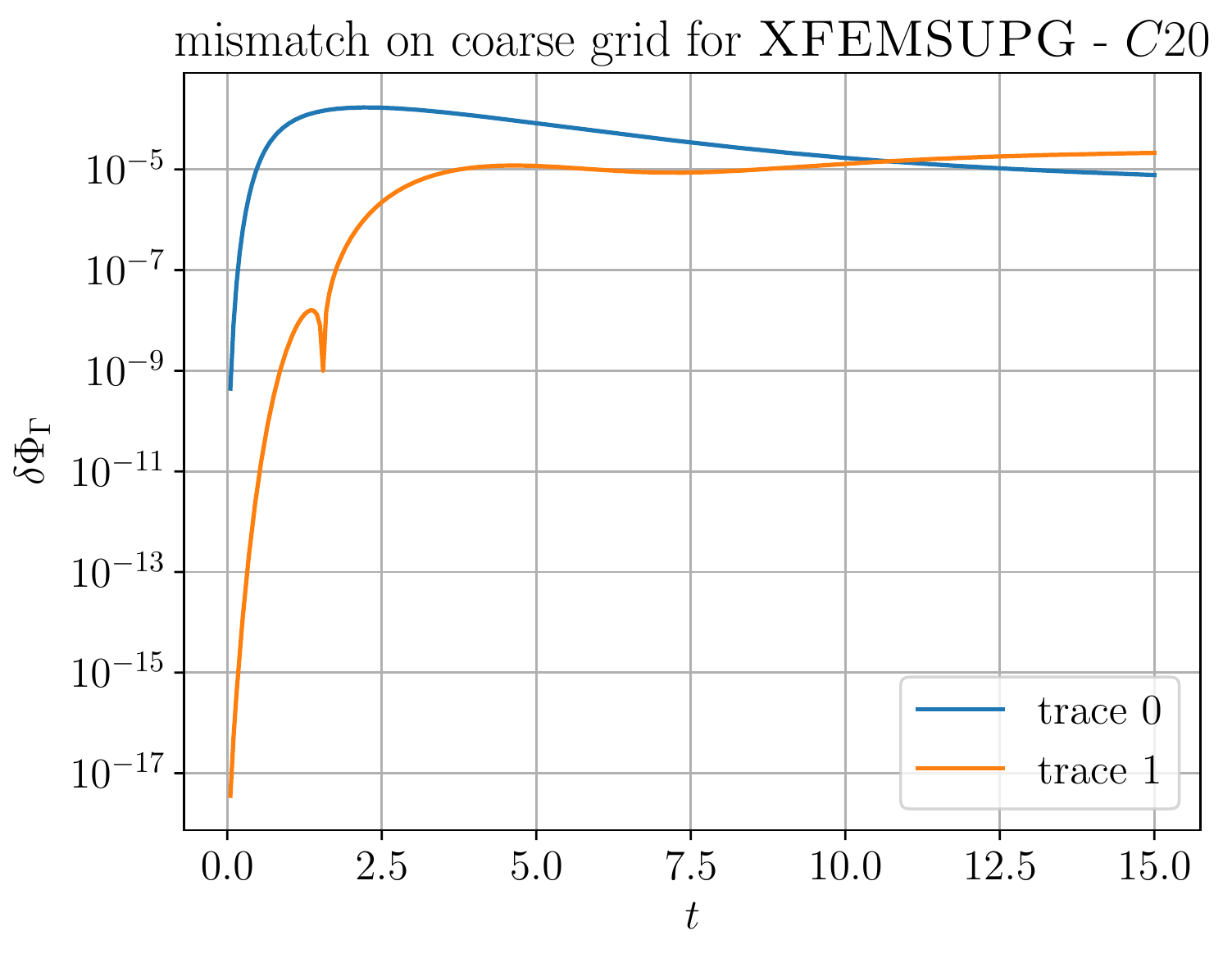}\\
       \includegraphics[width=0.33\textwidth]{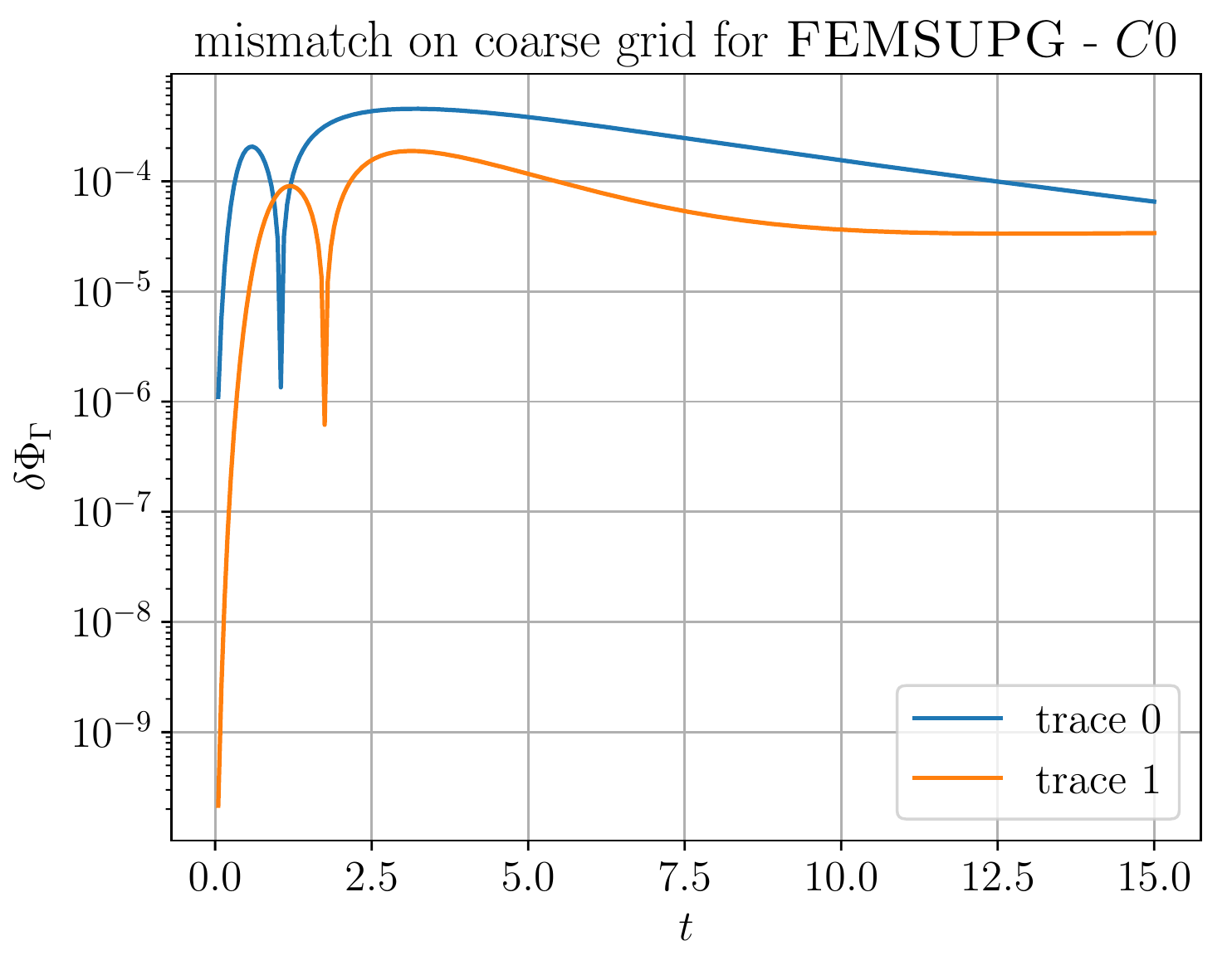}%
    \includegraphics[width=0.33\textwidth]{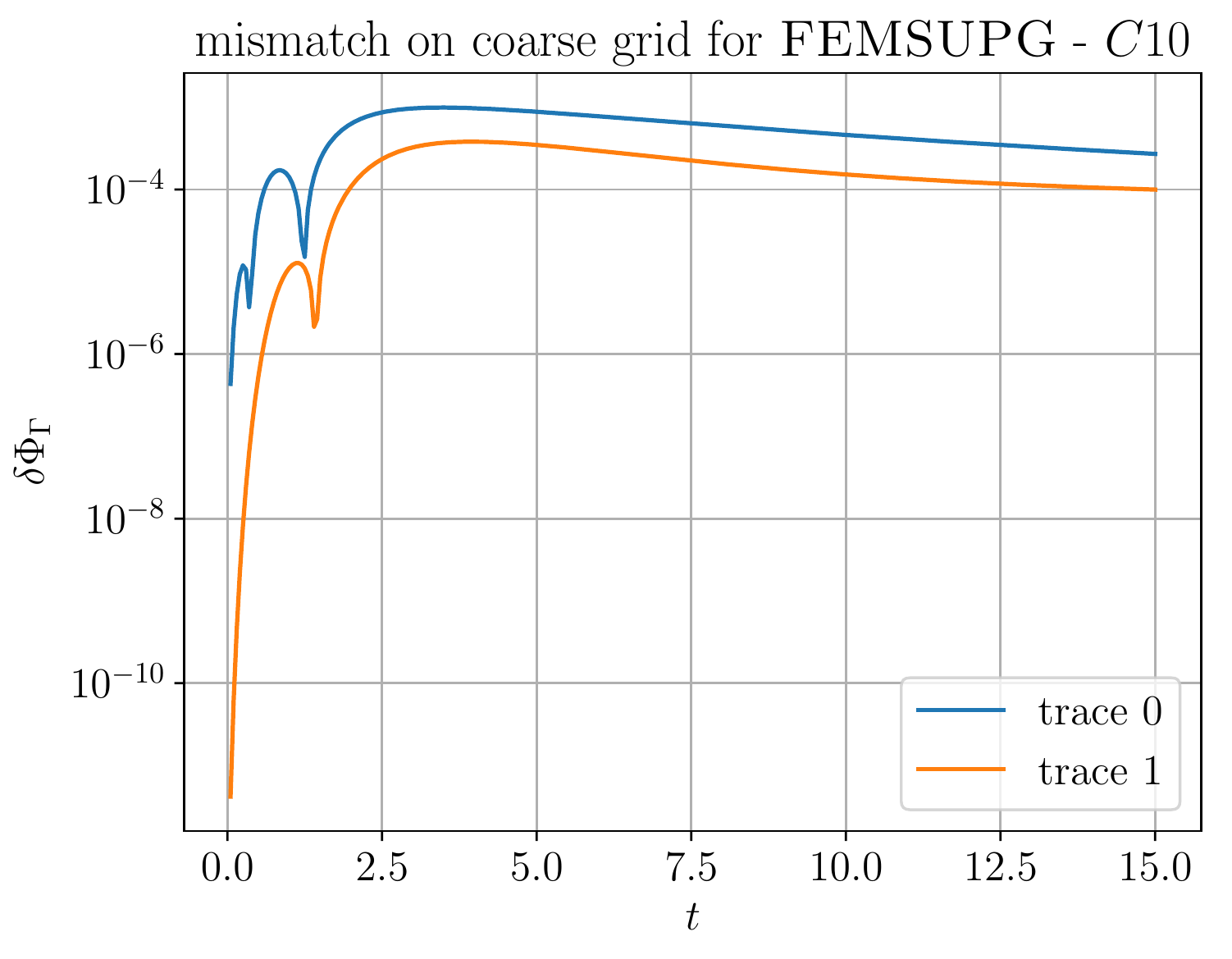}%
    \includegraphics[width=0.33\textwidth]{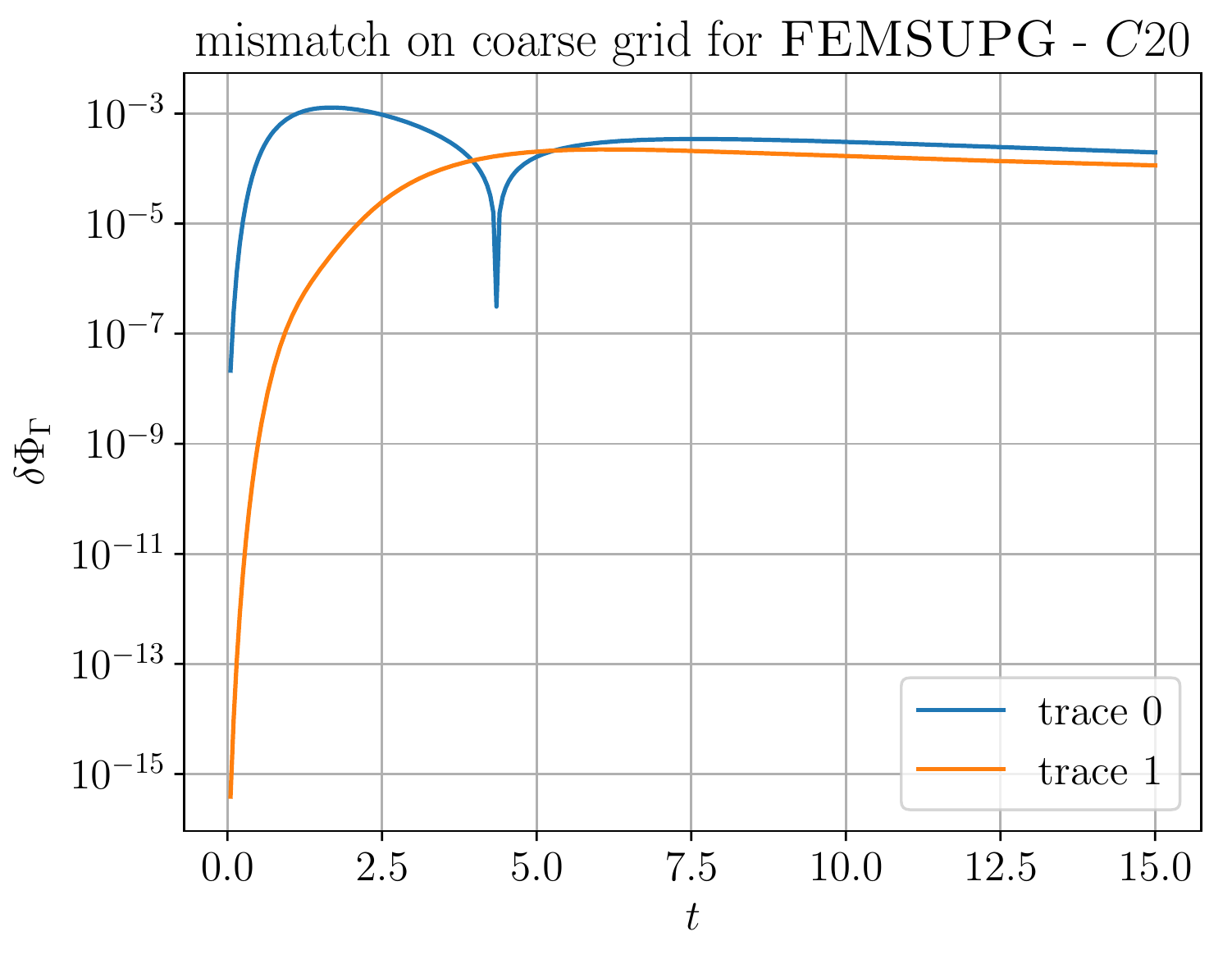}\\
    \caption{Total flux mismatch against time for Test Case 1 on traces
    $\Gamma_0$ (top) and $\Gamma_1$ (bottom) for configuration $C0$ (left),
    $C10$ (middle) and $C20$ (right)}
    \label{fig:example_1_mismatch}
    \end{figure}

As the non-conforming approaches \OPTXFEM{} and \OPTFEM{} are non locally
conservative, Figure~\ref{fig:example_1_mismatch} report the value of $\delta
\Phi$ against time on the coarsest mesh for the two traces of this example at
configurations $C0$, $C10$ and $C20$, from left to right, respectively, results
with the \OPTXFEM{} approach are on the top, results with \OPTFEM{} on the
bottom. The maximum-in-time absolute values of the total flux on traces
$\Gamma_0$ and $\Gamma_1$ are reported in Table~\ref{tab:example1_flux_values},
computed on the finest mesh as
$\Phi_{\Gamma_0}=\max_{t}\frac12\left(\Phi_{\Gamma_0,l}+\Phi_{\Gamma_0,c}\right)$,
and
$\Phi_{\Gamma_1}=\max_t\frac12\left(\Phi_{\Gamma_1,c}+\Phi_{\Gamma_1,r}\right)$.
We can see that relative values of less than $1\%$ are obtained for all times for
both methods, with the higher values corresponding to the configuration $C20$,
as expected. Moreover, for the larger times, mismatch values tend to decrease or
to remain constant at a fixed value. Thus this non-local-conservation has a
negligible impact on the computed solution. Further, mismatch errors can be
reduced by refining the mesh.

The mesh P\'eclet number for this problem ranges between a maximum of about
$6\times 10^2$ to a minimum of about $100$ on the computational meshes for
\OPTXFEM{}, \OPTFEM{} and \MFEM{} methods, thus a
Streamline-Upwind-Petrov-Galerkin (\SUPG{}) stabilization strategy was adopted. As
a consequence small overshoots/undershoots in the solution are observed in
\SUPG{}
stabilized methods, as well known in the literature, whereas the intrinsic
diffusive behavior of methods using upwinding for advection prevents this kind
of phenomena.

\begin{table}
\centering
\caption{Maximum value in time of total fluxes $\Phi_{\Gamma_0}$ and $\Phi_{\Gamma_1}$ on traces $\Gamma_0$ and $\Gamma_1$ of Test Case 1, respectively, computed with the \XFEMSUPG{} method on configurations $C0$, $C10$ and $C20$.}
\label{tab:example1_flux_values}
\begin{tabular}{c|ccc}
& $C0$ & $C10$ & $C20$\\
\hline
$\Phi_{\Gamma_0}$ &0.6072 & 0.4725 & 0.251 \\
$\Phi_{\Gamma_1}$ & 0.876 & 0.912 & 1.300\\
\hline
\end{tabular}
\end{table}


\subsection{Synthetic network}\label{subsec:example2}

\begin{figure}[htbp]
    \centering
    \resizebox{0.31\textwidth}{!}{\fontsize{20pt}{8}\selectfont%
\begingroup%
  \makeatletter%
  \providecommand\color[2][]{%
    \errmessage{(Inkscape) Color is used for the text in Inkscape, but the package 'color.sty' is not loaded}%
    \renewcommand\color[2][]{}%
  }%
  \providecommand\transparent[1]{%
    \errmessage{(Inkscape) Transparency is used (non-zero) for the text in Inkscape, but the package 'transparent.sty' is not loaded}%
    \renewcommand\transparent[1]{}%
  }%
  \providecommand\rotatebox[2]{#2}%
  \newcommand*\fsize{\dimexpr\f@size pt\relax}%
  \newcommand*\lineheight[1]{\fontsize{\fsize}{#1\fsize}\selectfont}%
  \ifx\svgwidth\undefined%
    \setlength{\unitlength}{293.39062788bp}%
    \ifx\svgscale\undefined%
      \relax%
    \else%
      \setlength{\unitlength}{\unitlength * \real{\svgscale}}%
    \fi%
  \else%
    \setlength{\unitlength}{\svgwidth}%
  \fi%
  \global\let\svgwidth\undefined%
  \global\let\svgscale\undefined%
  \makeatother%
  \begin{picture}(1,0.9077297)%
    \lineheight{1}%
    \setlength\tabcolsep{0pt}%
    \put(0,0){\includegraphics[width=\unitlength,page=1]{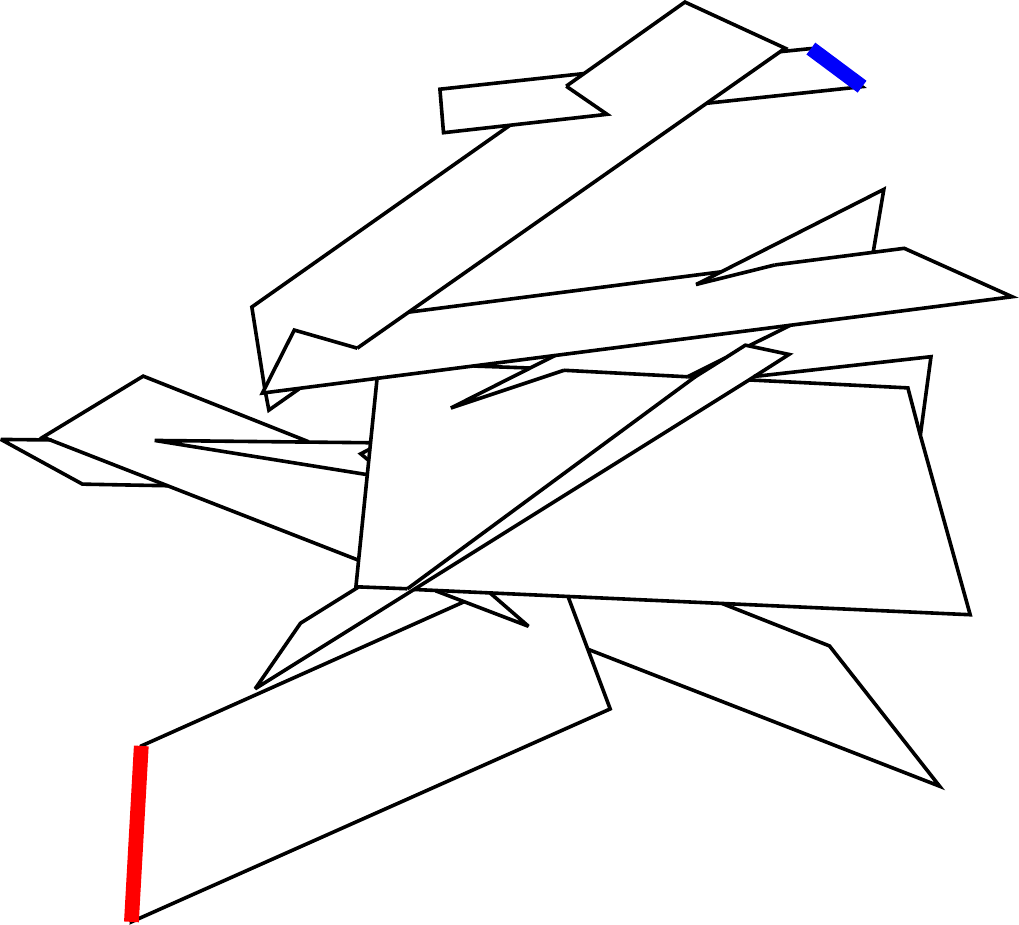}}%
    \put(0.15736293,0.09296195){\color[rgb]{0,0,0}\makebox(0,0)[lt]{\lineheight{1.25}\smash{\begin{tabular}[t]{l}inflow\end{tabular}}}}%
    \put(0.7039781,0.75463738){\color[rgb]{0,0,0}\makebox(0,0)[lt]{\lineheight{1.25}\smash{\begin{tabular}[t]{l}outflow\end{tabular}}}}%
  \end{picture}%
\endgroup%
}
    \includegraphics[width=0.33\textwidth]{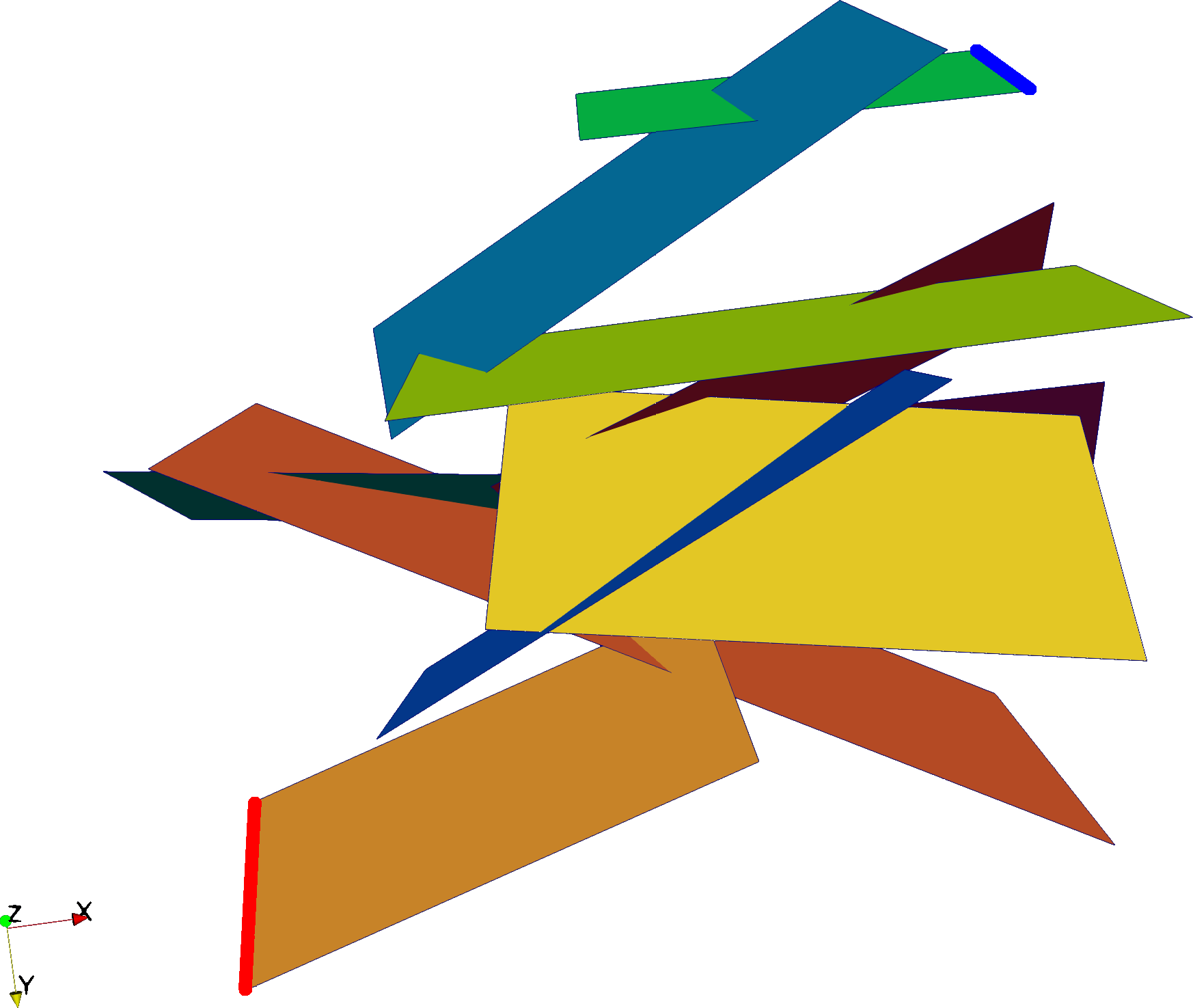}%
    \includegraphics[width=0.33\textwidth]{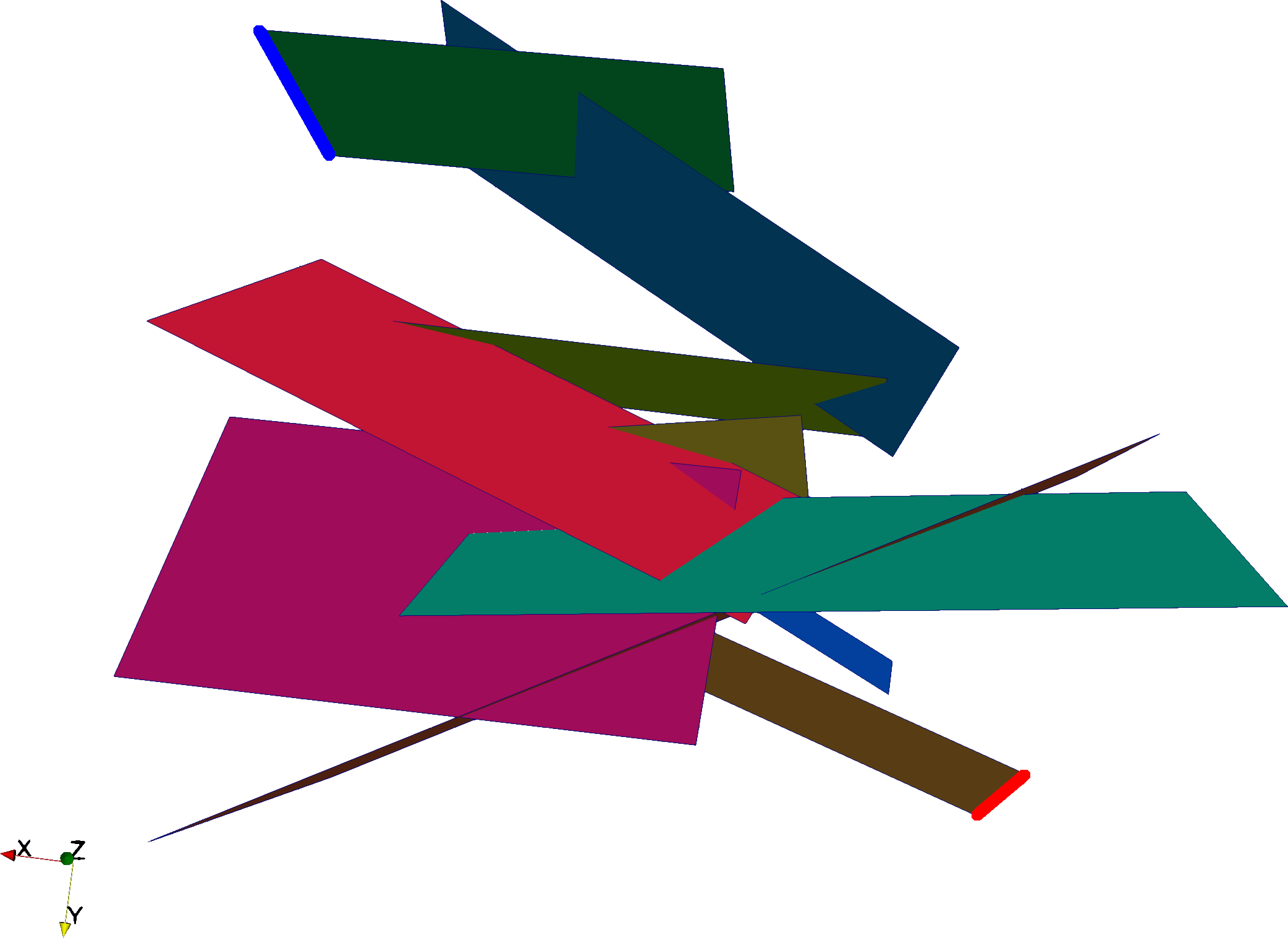}
    \caption{Geometry of Test Case 2 and two views of the network.
    A red line represents the inflow and a blue thick line the outflow part of the
    boundary, respectively.}%
    \label{fig:geometry_example_2}
\end{figure}

In the second test, labeled Test Case 2, we consider a more complex network composed of
10 fractures with 14 traces, thus being more similar to (a portion of) realistic
DFNs, still remaining simple enough to perform analyses on the solutions
obtained with the considered schemes. The network is represented in
Figure~\ref{fig:geometry_example_2}, where the inflow and outflow portions of
the boundary are also marked. Also in this case the Darcy velocity $\bm{u}$ is
first computed solving problem \ref{pb:model_darcy_weak} or
\ref{pb:model_darcy_primal_weak}, depending on the chosen numerical scheme, with
constant hydraulic conductivity equal to one on all fractures and null source term. A
unitary pressure head drop is imposed between the inflow and outflow portions of the
border, all other fracture edges being insulated. The Darcy velocity is then
used as an input for an advection-diffusion problem for a scalar quantity
$\theta$, as formulated in \ref{pb:model_heat_weak}, with null source and
reaction terms. A coefficient $\zeta=1$ is chosen, whereas the diffusion
coefficient is $D=10^{-4}$. A unitary Dirichlet boundary condition is prescribed
on the inflow border and all other edges are insulated.

\begin{figure}[htbp]
    \centering
    \includegraphics[width=0.25\textwidth]{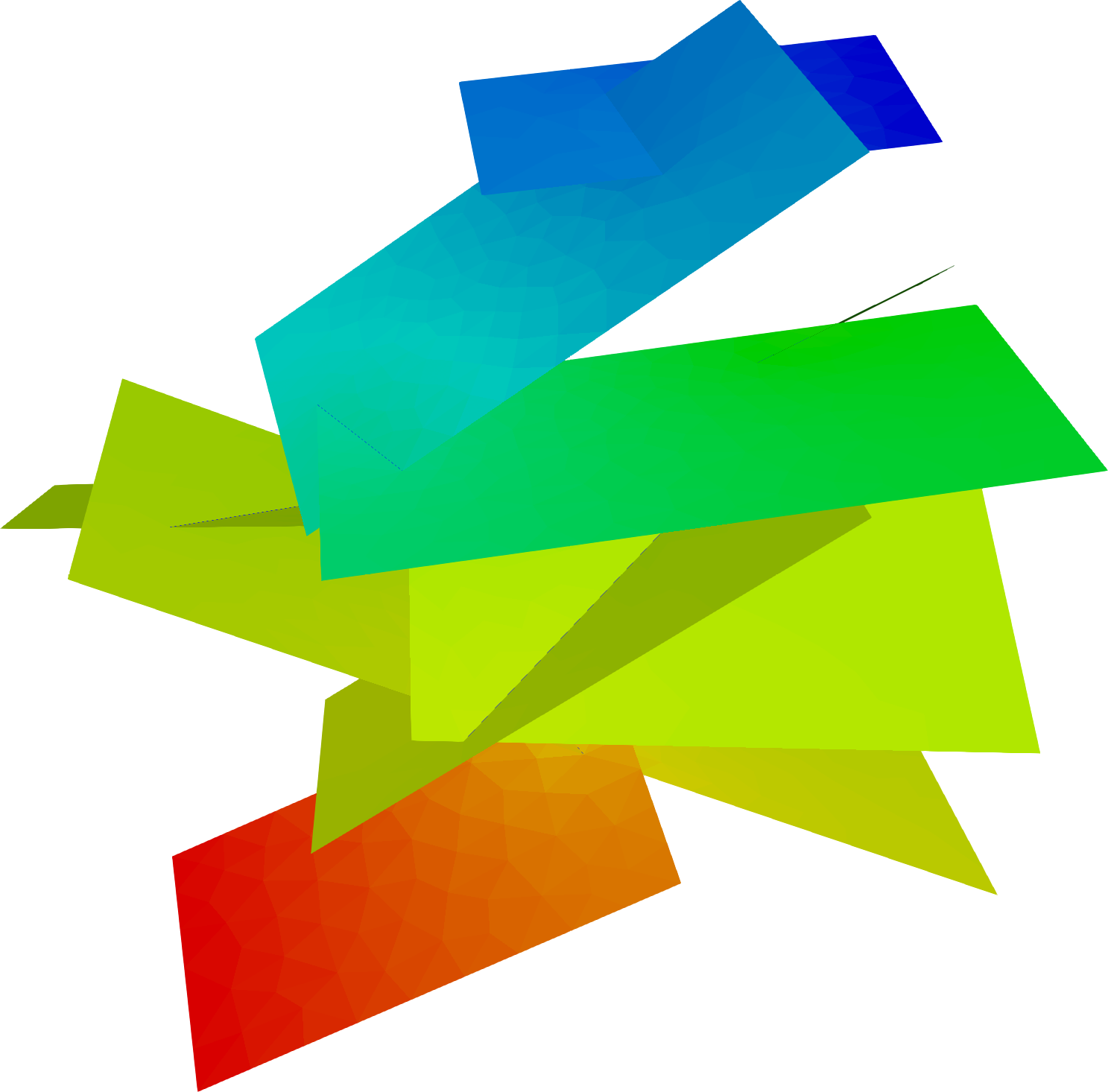}%
    \includegraphics[width=0.25\textwidth]{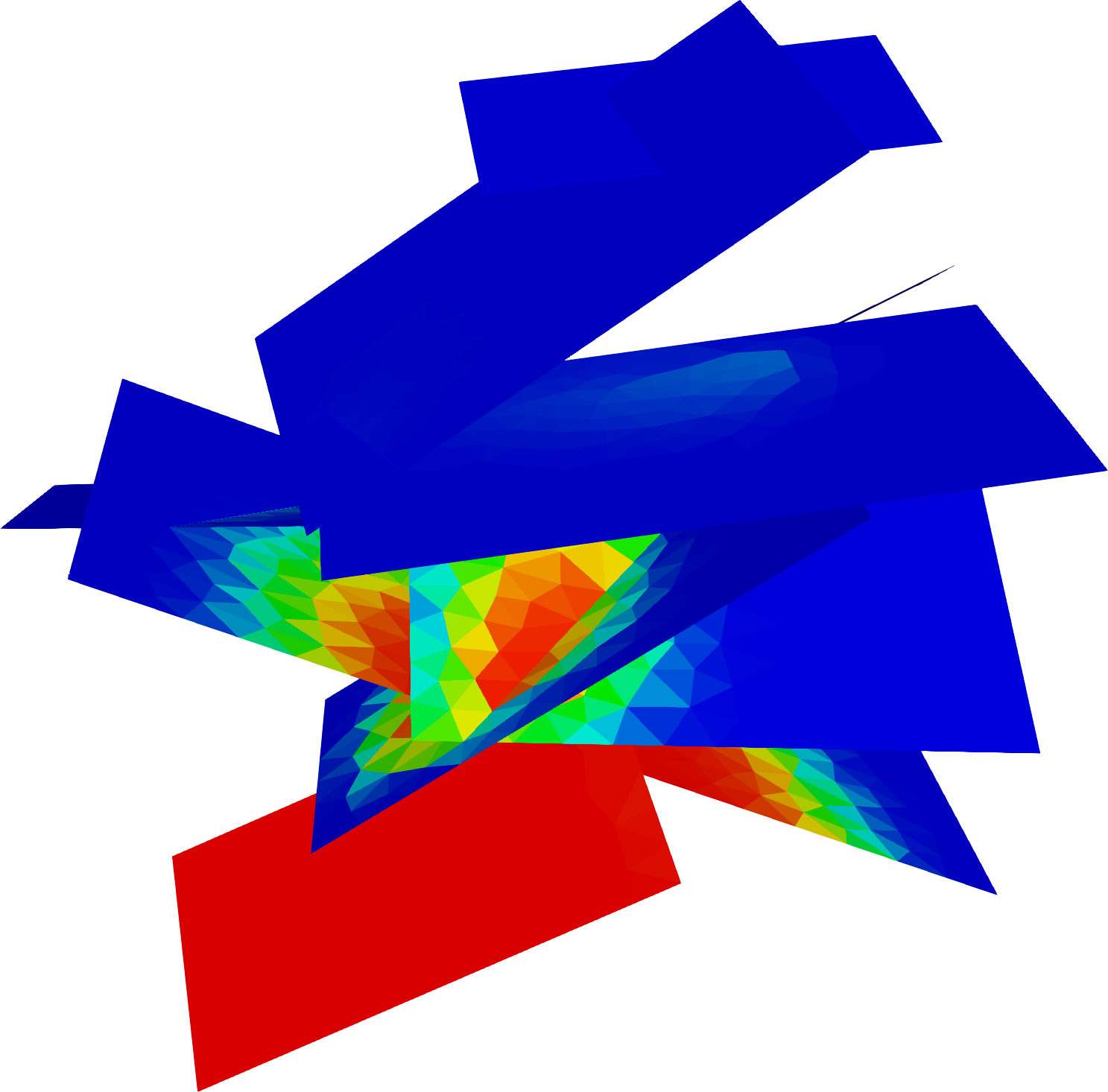}%
    \includegraphics[width=0.25\textwidth]{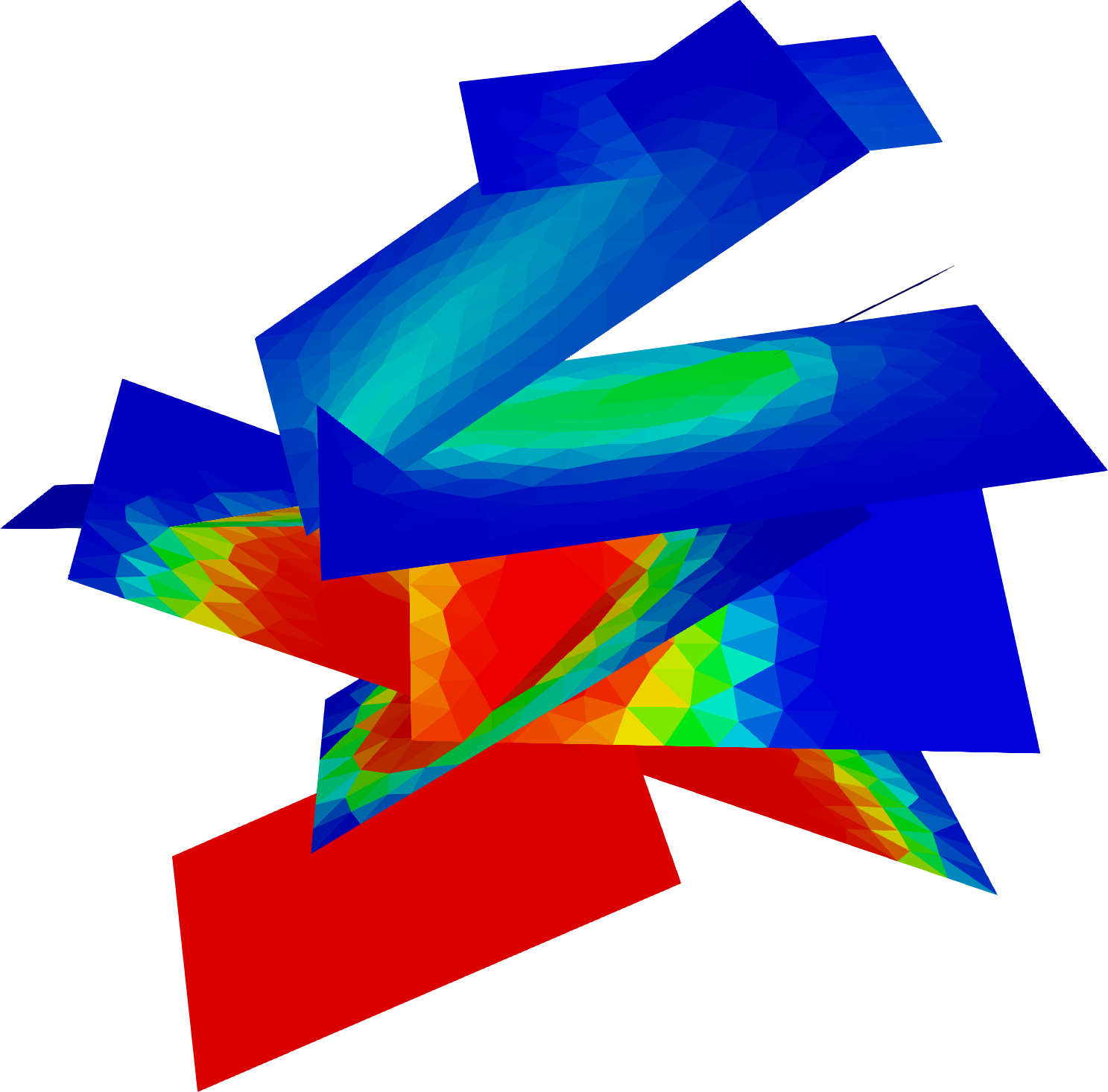}%
    \includegraphics[width=0.25\textwidth]{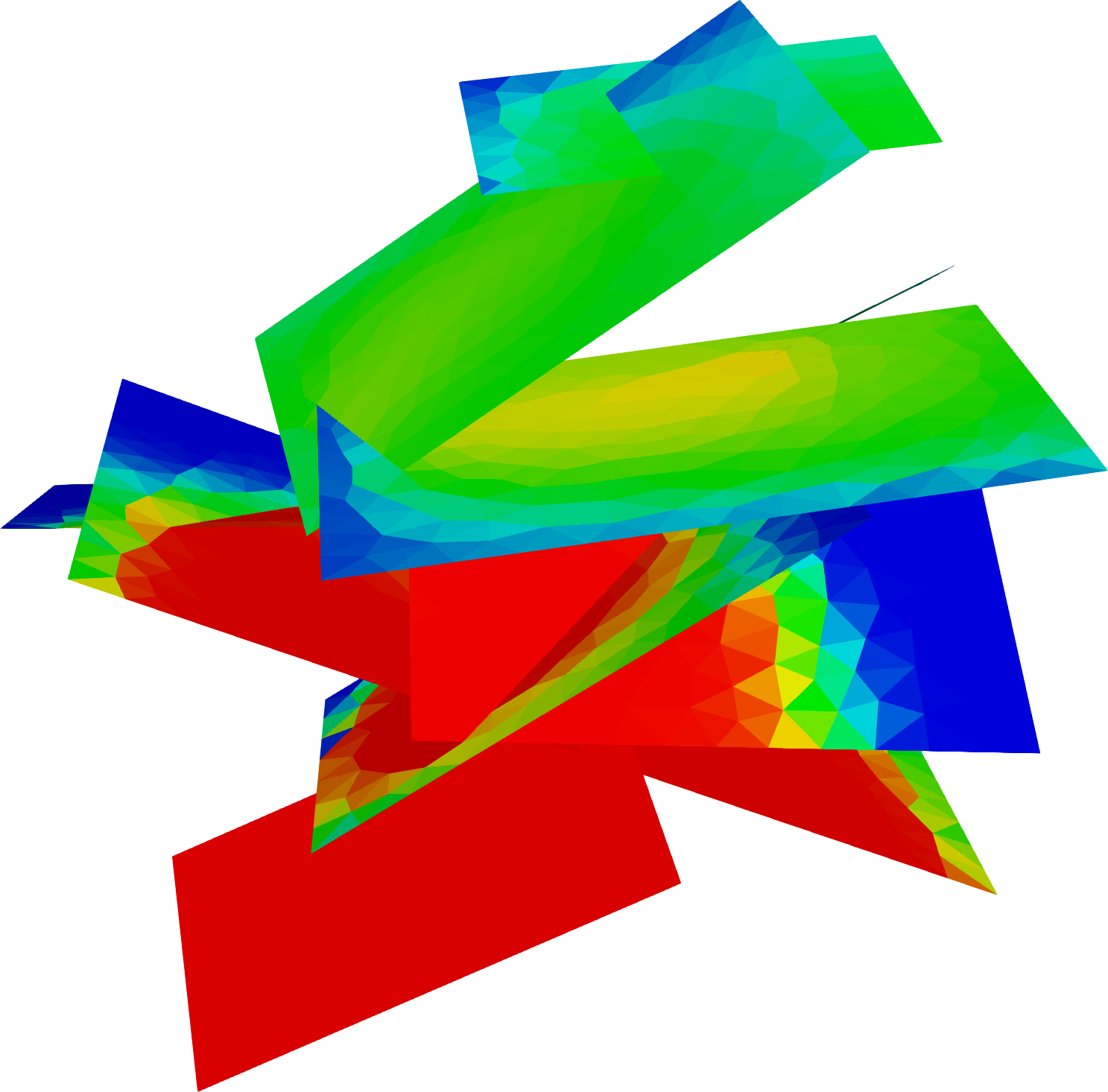}\\
    \includegraphics[width=0.25\textwidth]{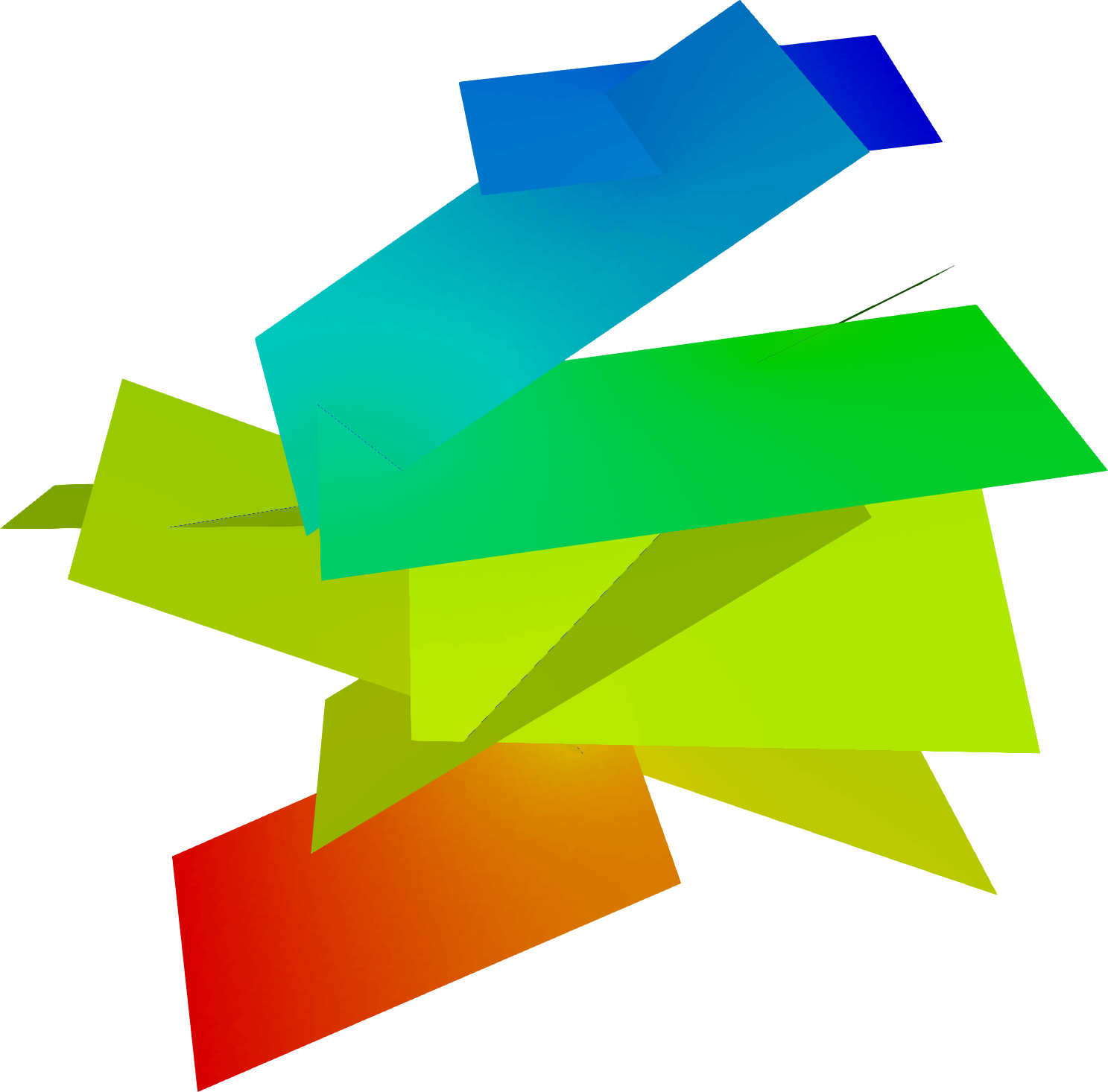}%
    \includegraphics[width=0.25\textwidth]{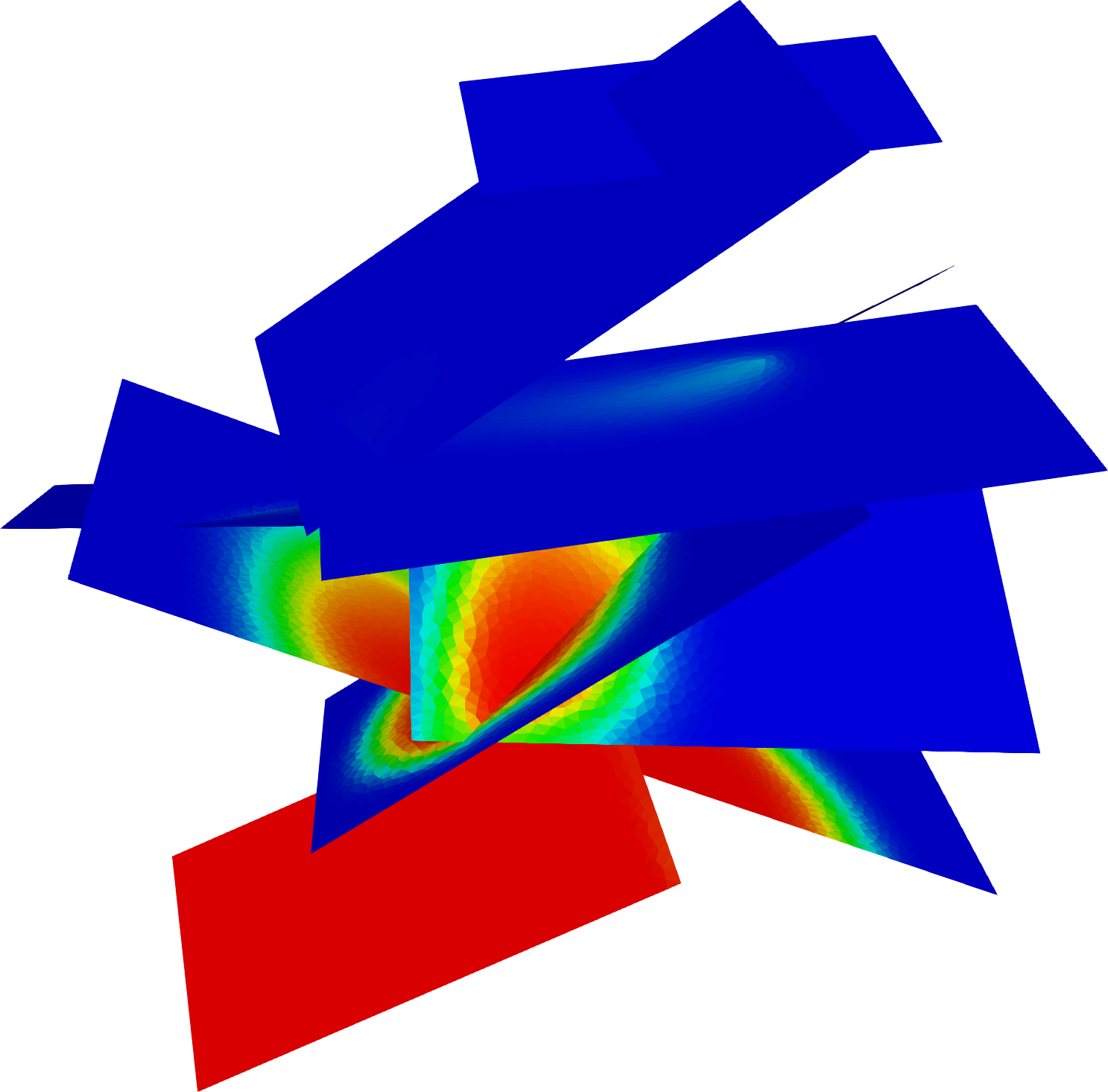}%
    \includegraphics[width=0.25\textwidth]{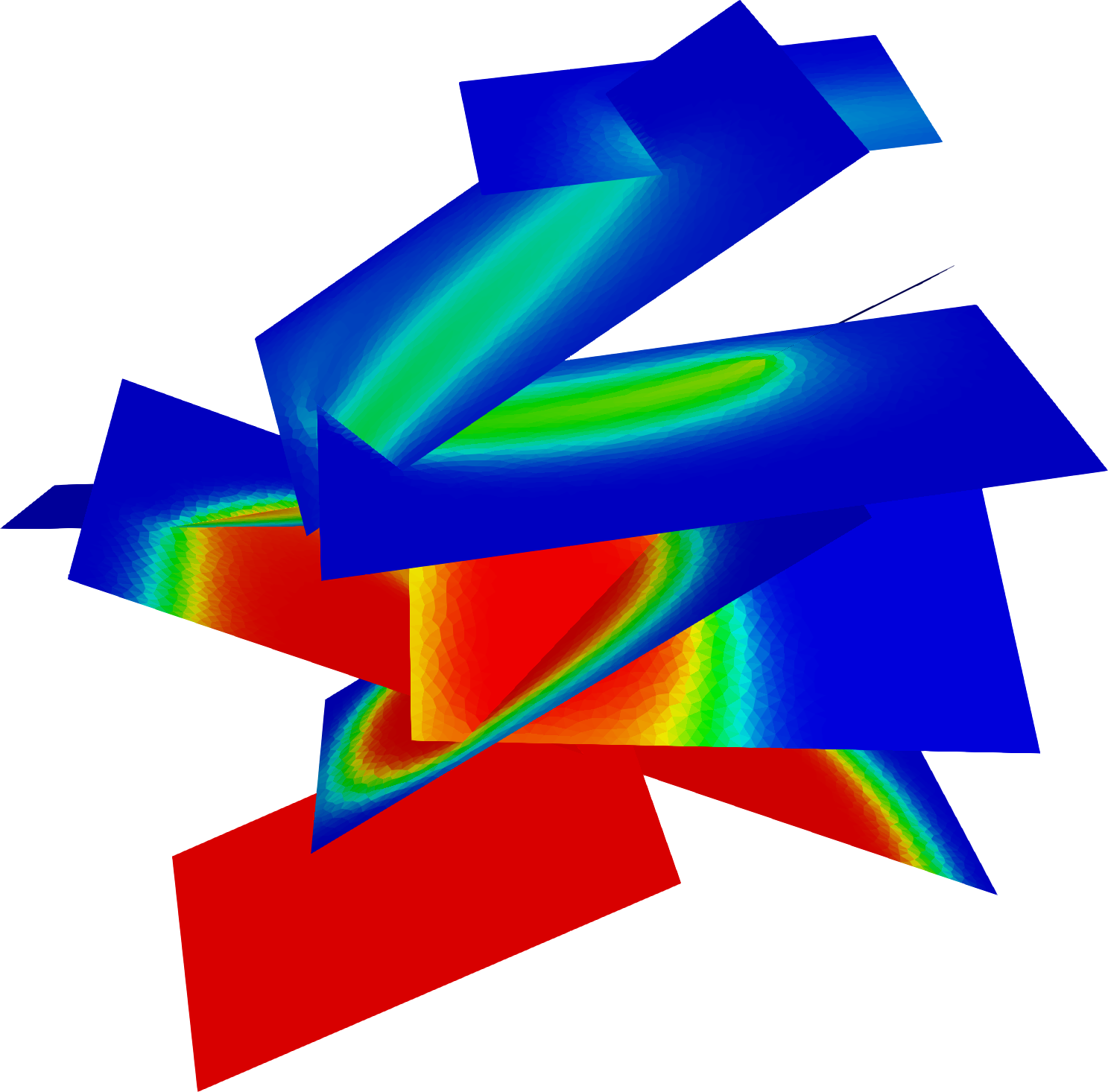}%
    \includegraphics[width=0.25\textwidth]{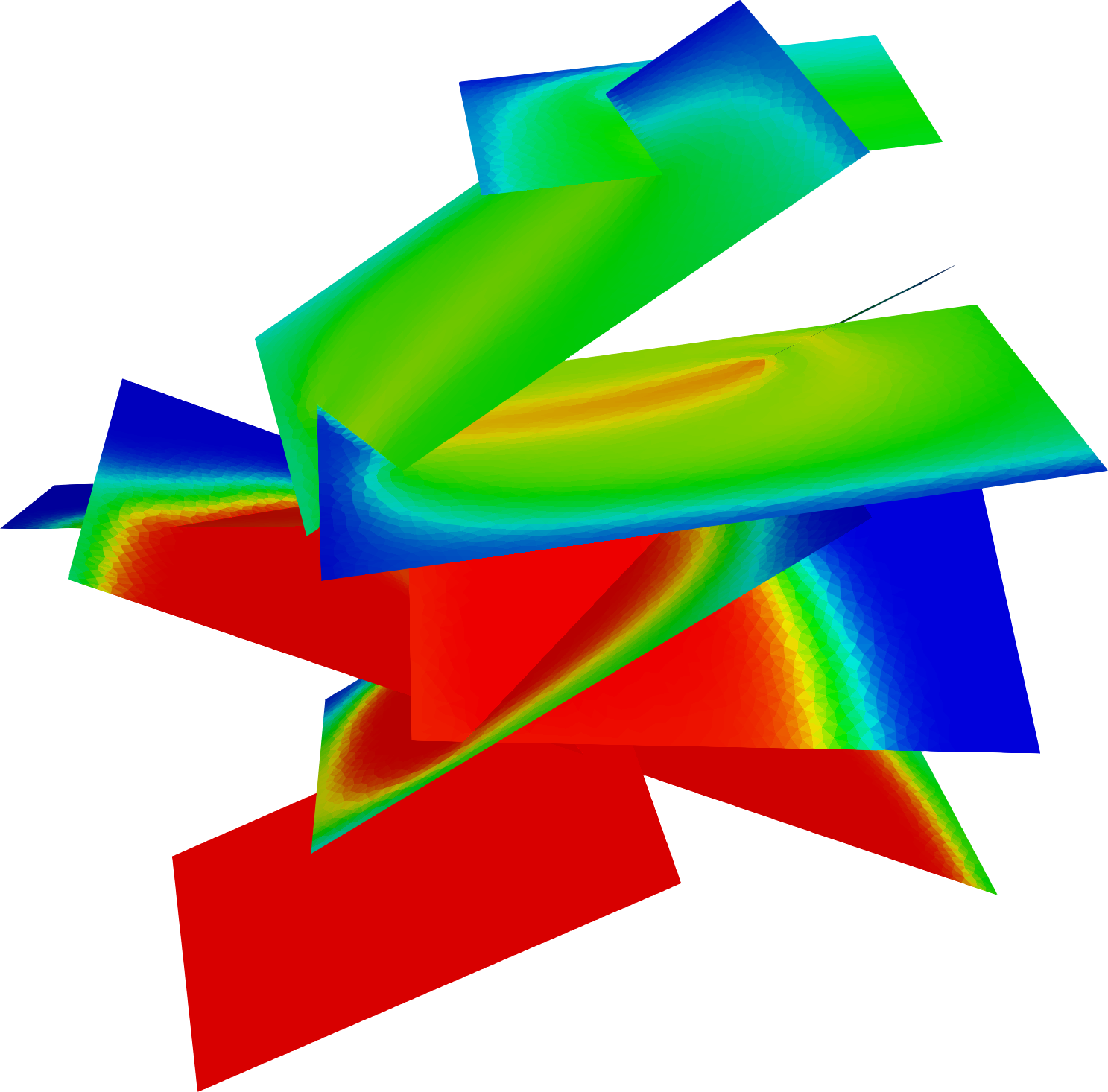}
    \caption{Solution of the test case in Subsection \ref{subsec:example2}. On the first column
    the pressure head solution for the two grid refinements. The other
    columns show $\theta$ at different time: 3.35, 6.25, and 12.5, respectively,
    for both level of refinement. In
    all cases the solution is rescaled in the range $[0, 1]$.}%
    \label{fig:solution_example_2}
\end{figure}

Two meshes are used also for Test Case 2 counting about $10^3$ elements (coarse
mesh) and $4\times 10^4$ elements (fine mesh), respectively. The mesh P\'eclet number,
related to the coarse mesh for \SUPG{} stabilized methods, is of about $100$. An equally-spaced
mesh is then used for time discretization with $500$ steps of length $0.05$. Also in this case, the initial solution is the null function.
An example of the obtained numerical solution with the \TPFA{} method is shown in
Figure~\ref{fig:solution_example_2}: in the first column on the left the
pressure head distribution in the network is represented, solution of the Darcy
problem on the coarse (top) and fine (bottom) meshes; the remaining columns
depict the solution $\theta$ of the dispersion problem at three selected time
frames, corresponding to $t=3.35$, $t=6.25$ and $t=12.5$, again on the coarse
(top) and fine (bottom) mesh. Coherently the heat flows from the inflow to the
outflow by following a tortuous path given by the complex disposition of the
fractures and traces. The solution on the coarse grid has a spreader front than
on the fine grid due to the artificial diffusivity of the scheme.

\begin{figure}
\centering
    \includegraphics[width=0.32\textwidth]{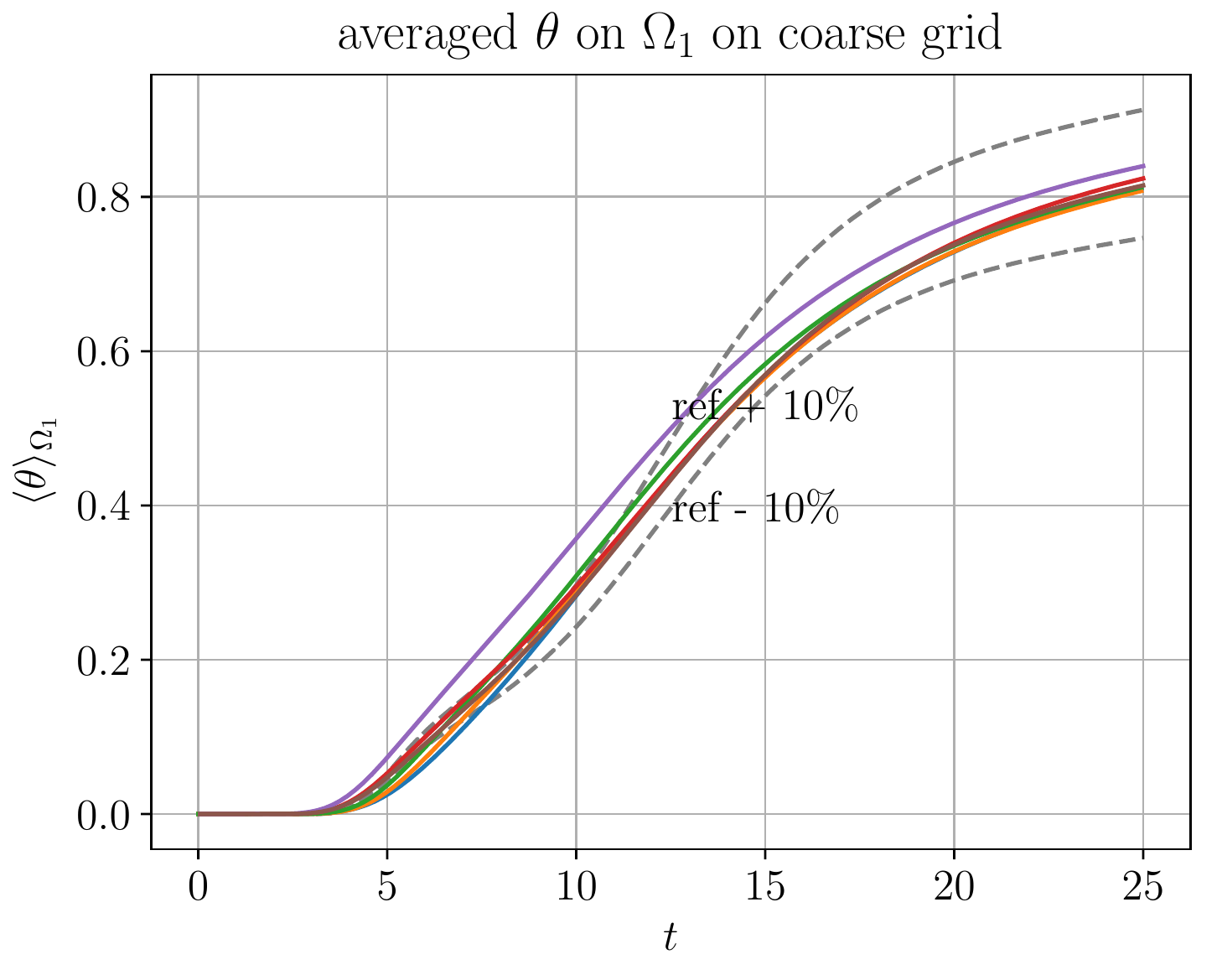}
    \includegraphics[width=0.33\textwidth]{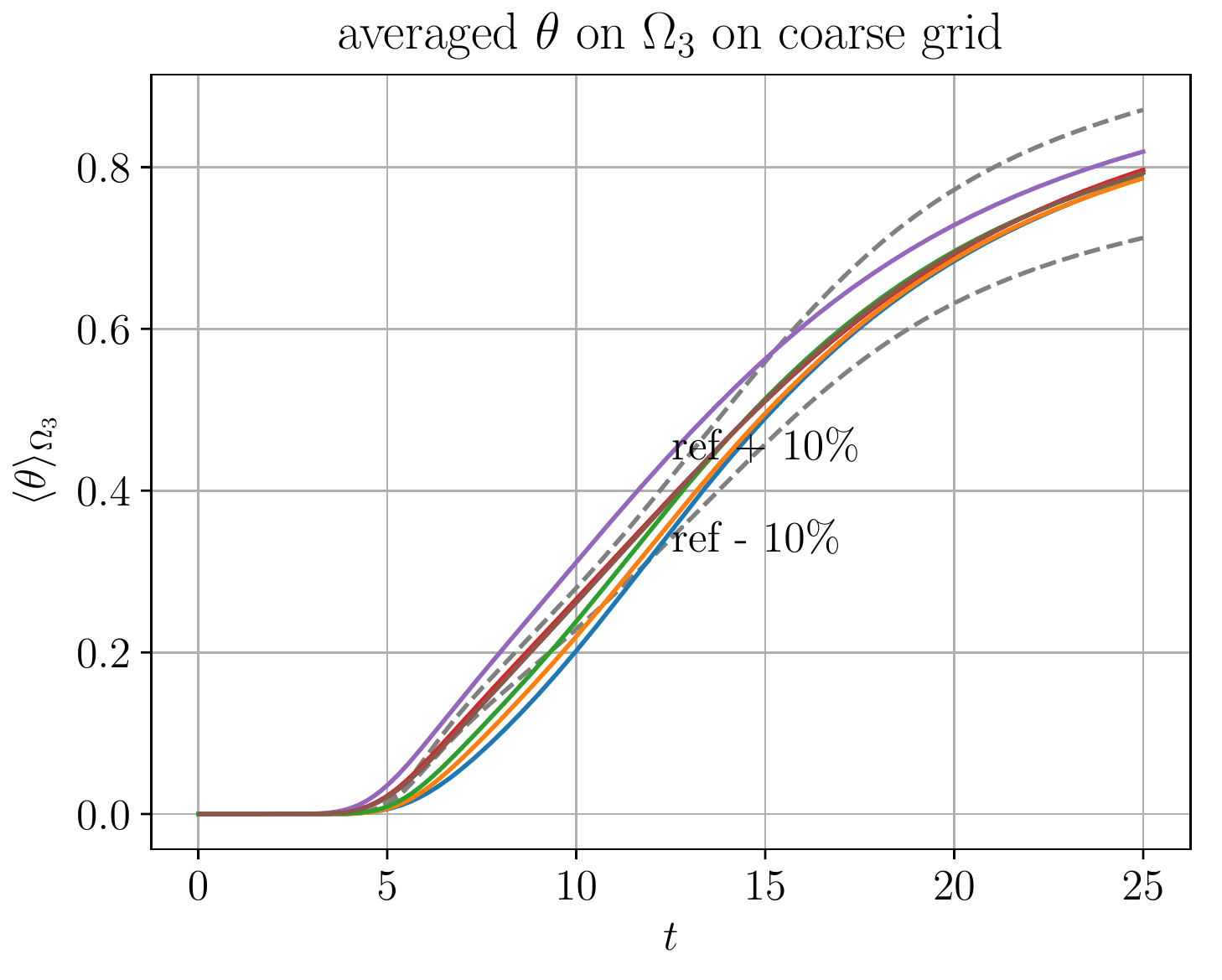}
    \includegraphics[width=0.33\textwidth]{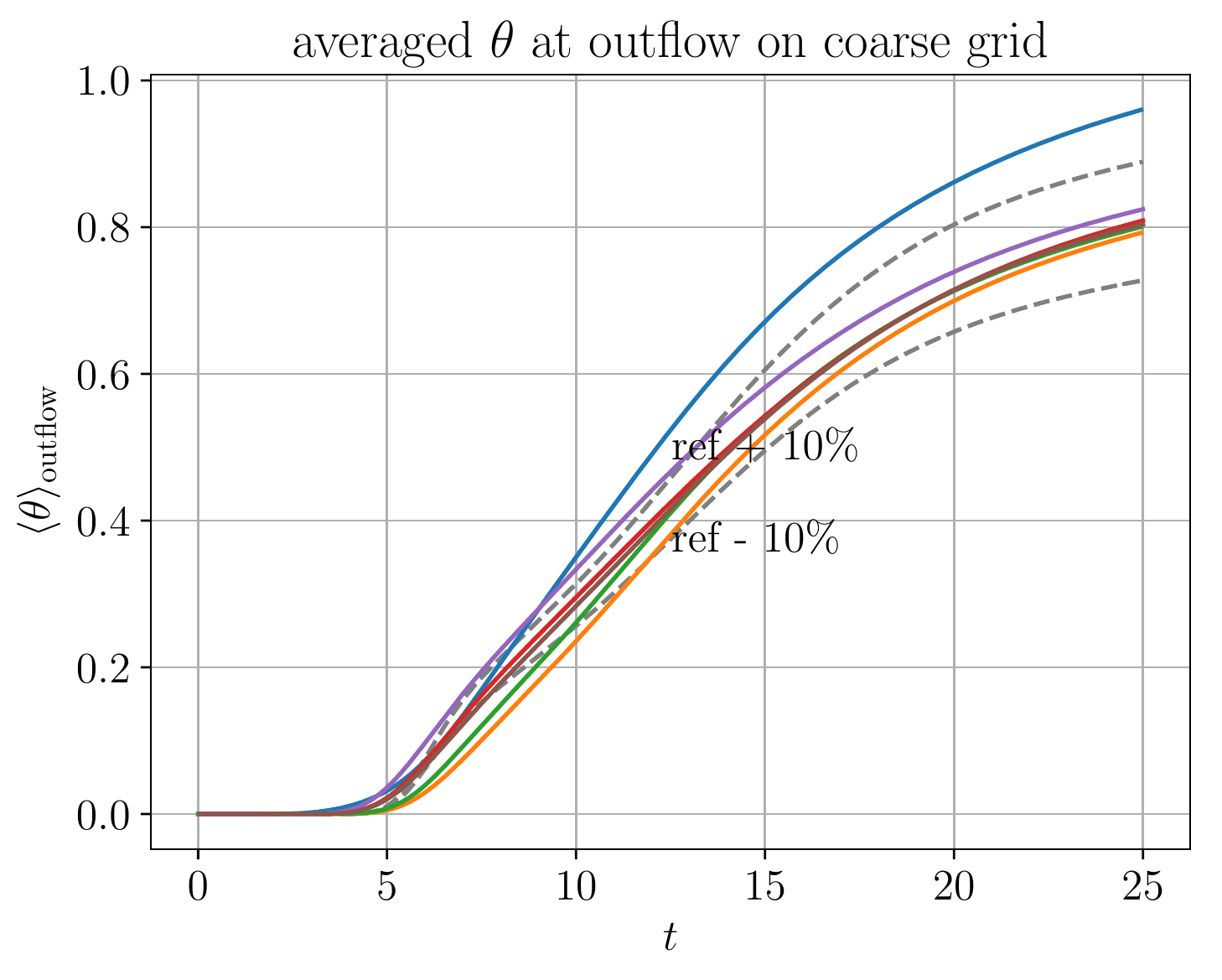}\\
    \includegraphics[width=0.32\textwidth]{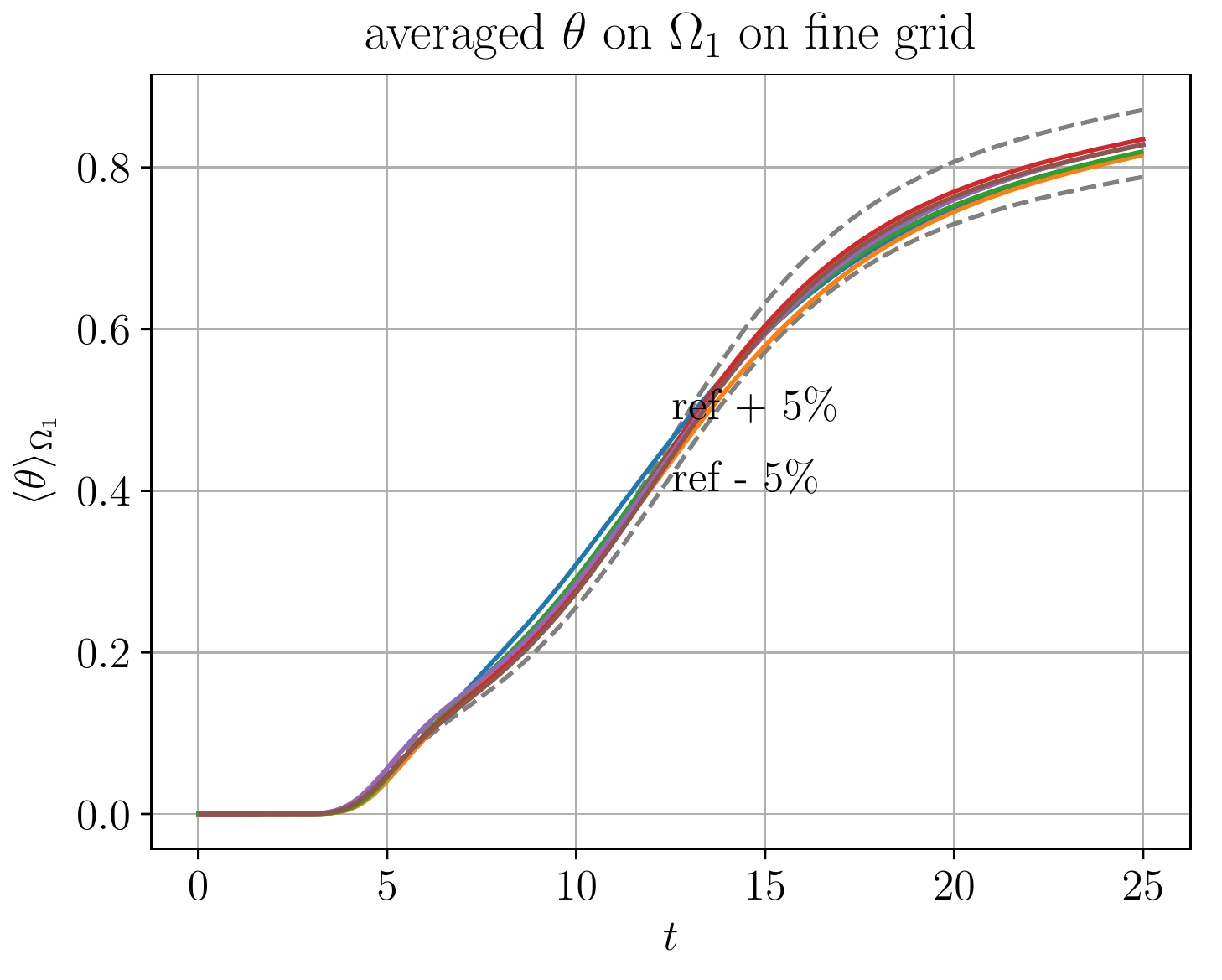}
    \includegraphics[width=0.33\textwidth]{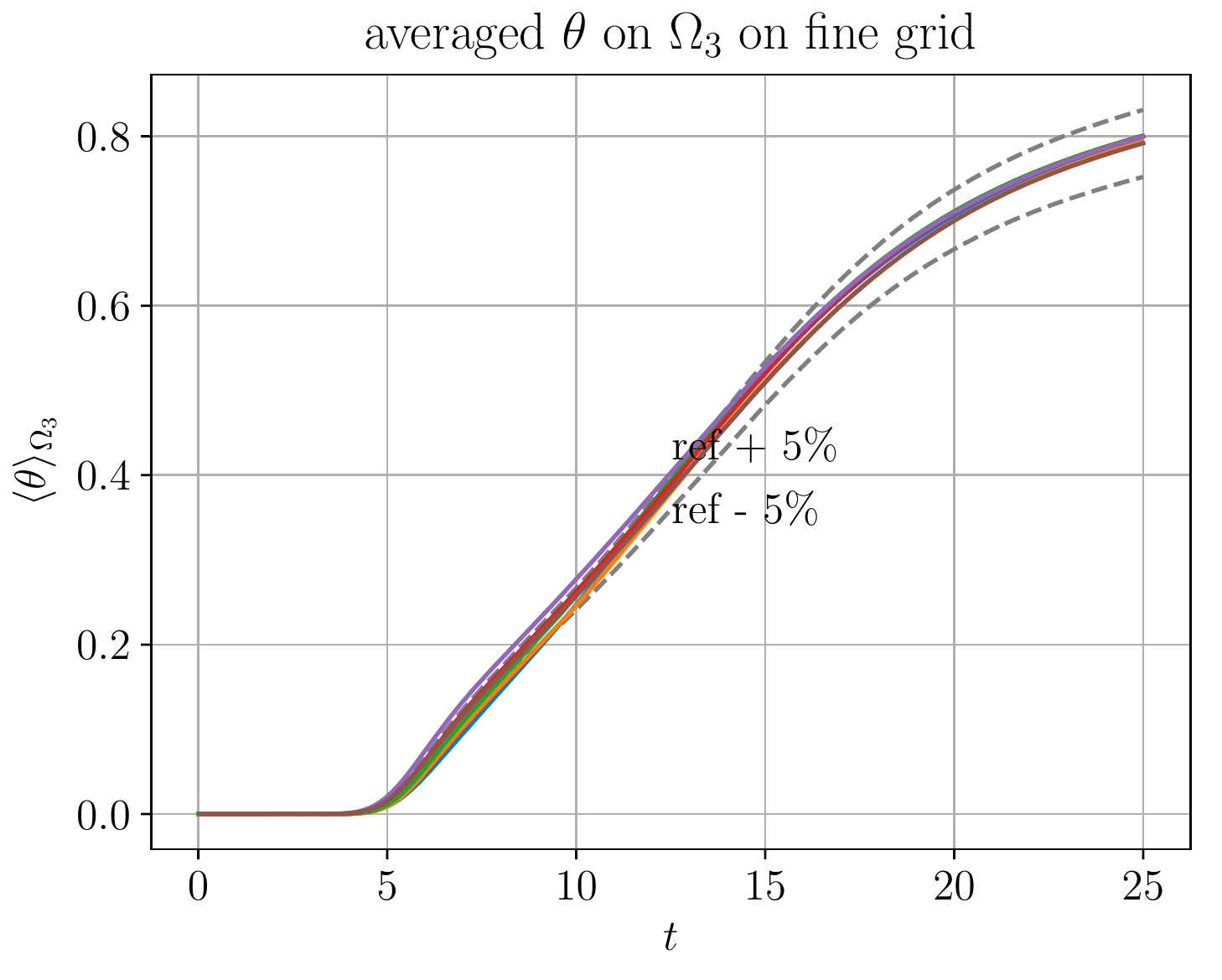}
    \includegraphics[width=0.33\textwidth]{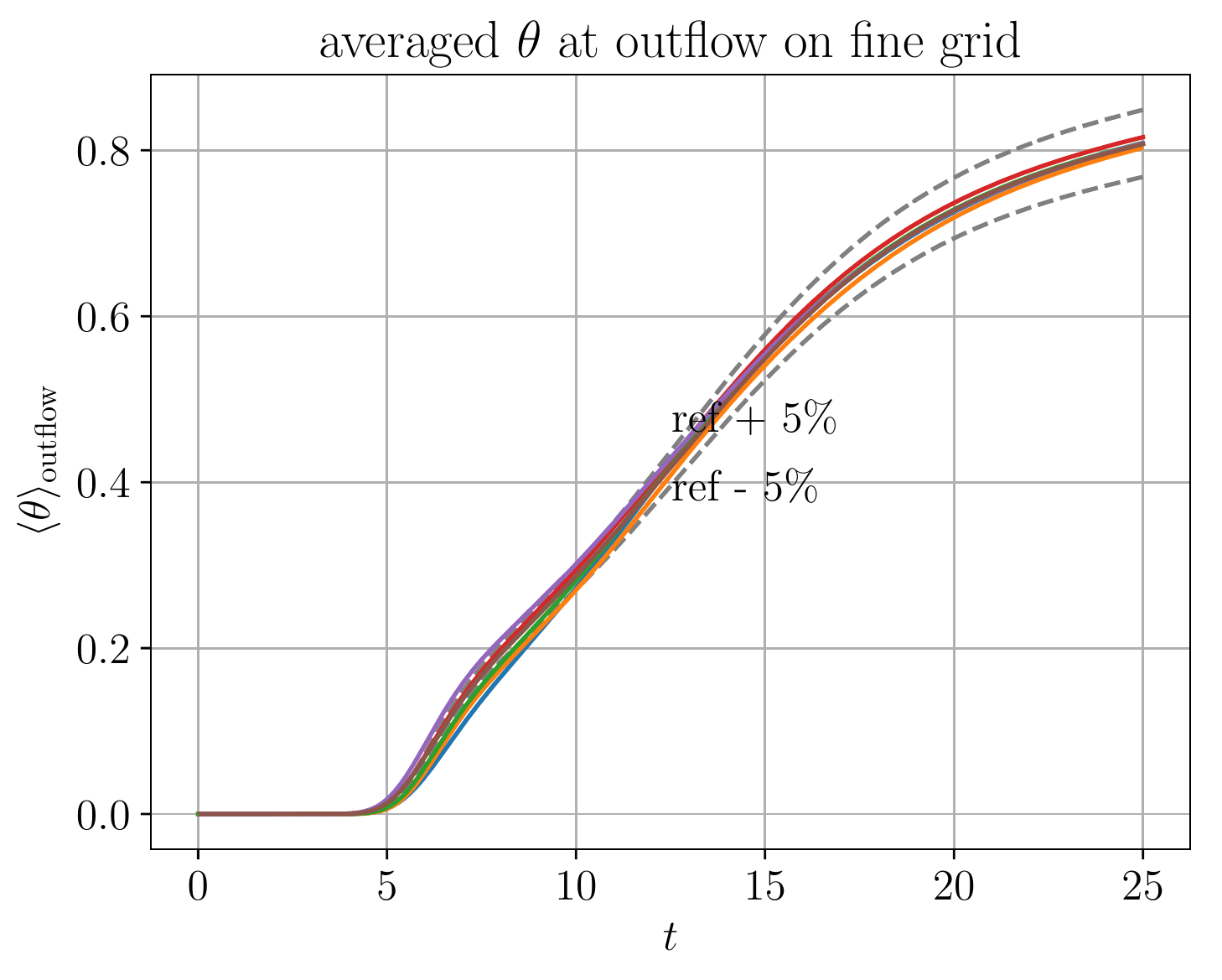}\\
    \includegraphics[width=\textwidth]{label}
\caption{Average temperature on two selected fractures (left and middle columns)
    and at outflow (right column) curves against time for Test Case 2. Coarse mesh on top, fine mesh at the bottom.}
\label{fig:example2_avg_prod}
\end{figure}

As previously, the average temperature $\langle{\theta}\rangle_\Omega$ on
selected fractures, the average outflow temperature $\langle{\theta}\rangle_{\mathrm{outflow}}$
and flux mismatch $\delta \Phi_\Gamma$ at traces are used to compare and assess
the approximation capabilities of the various schemes. The curves of
$\langle{\theta} \rangle_{\Omega_1}$, $\langle{\theta} \rangle_{\Omega_3}$ and
$\langle{\theta} \rangle_{\mathrm{outflow}}$ are reported in
Figure~\ref{fig:example2_avg_prod}. A reference solution is computed with the
\MFEMSUPG{} method on a mesh counting about $2\times 10^4$ cells, and relative
error curves with respect to this solution are shown in dashed lines in the
Figure~\ref{fig:example2_avg_prod}. We can observe that, despite the network has a larger number of
fractures and fracture intersections with respect to Test Case 1, all the
methods, in absence of severe geometrical features have good approximation
properties, that further improve with mesh refinement. The larger differences
are observed for the average outflow temperature curve related to the \MVEMUP{} approach on the
coarse mesh and for the average temperature curves of \FEMSUPG{} again on the
coarse mesh. In both cases differences slightly exceed $10\%$ of the reference,
and are reduced by mesh refinement. Concerning the \MVEMUP{} the difference is
caused by the coarsening process which reduces the number of elements of the
original mesh to about one half. Also, the \TPFA{} method used for advection is
observed to have poor performances on irregular polygonal cells, as the ones
generated by the coarsening method. Concerning the \FEMSUPG{} method, some
discrepancies with respect to the reference are to be expected, as the approach
is designed to be computationally inexpensive, nonetheless it is capable of
providing satisfactory predictions.

\begin{figure}
\centering
\includegraphics[width=0.48\textwidth]{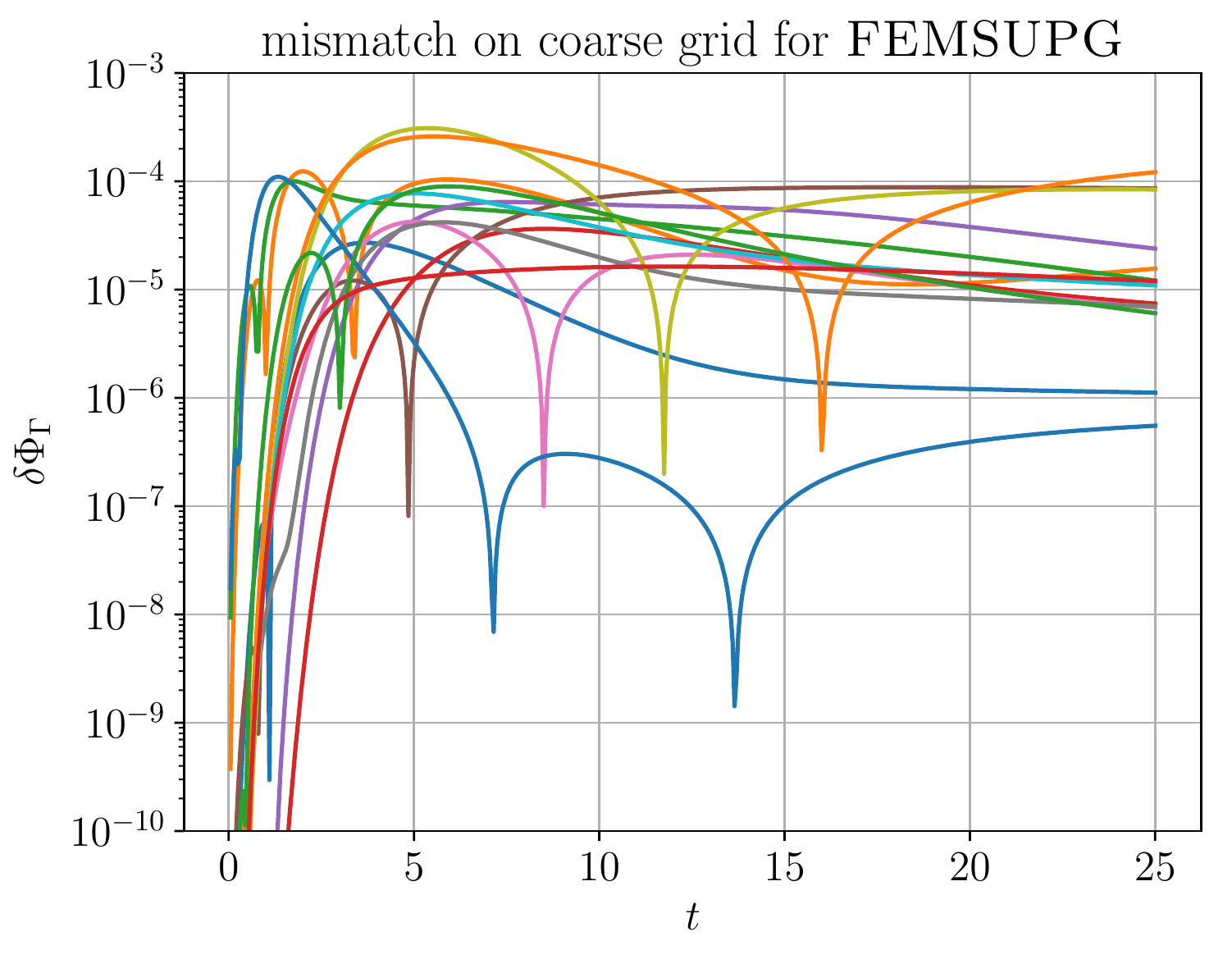}
\includegraphics[width=0.48\textwidth]{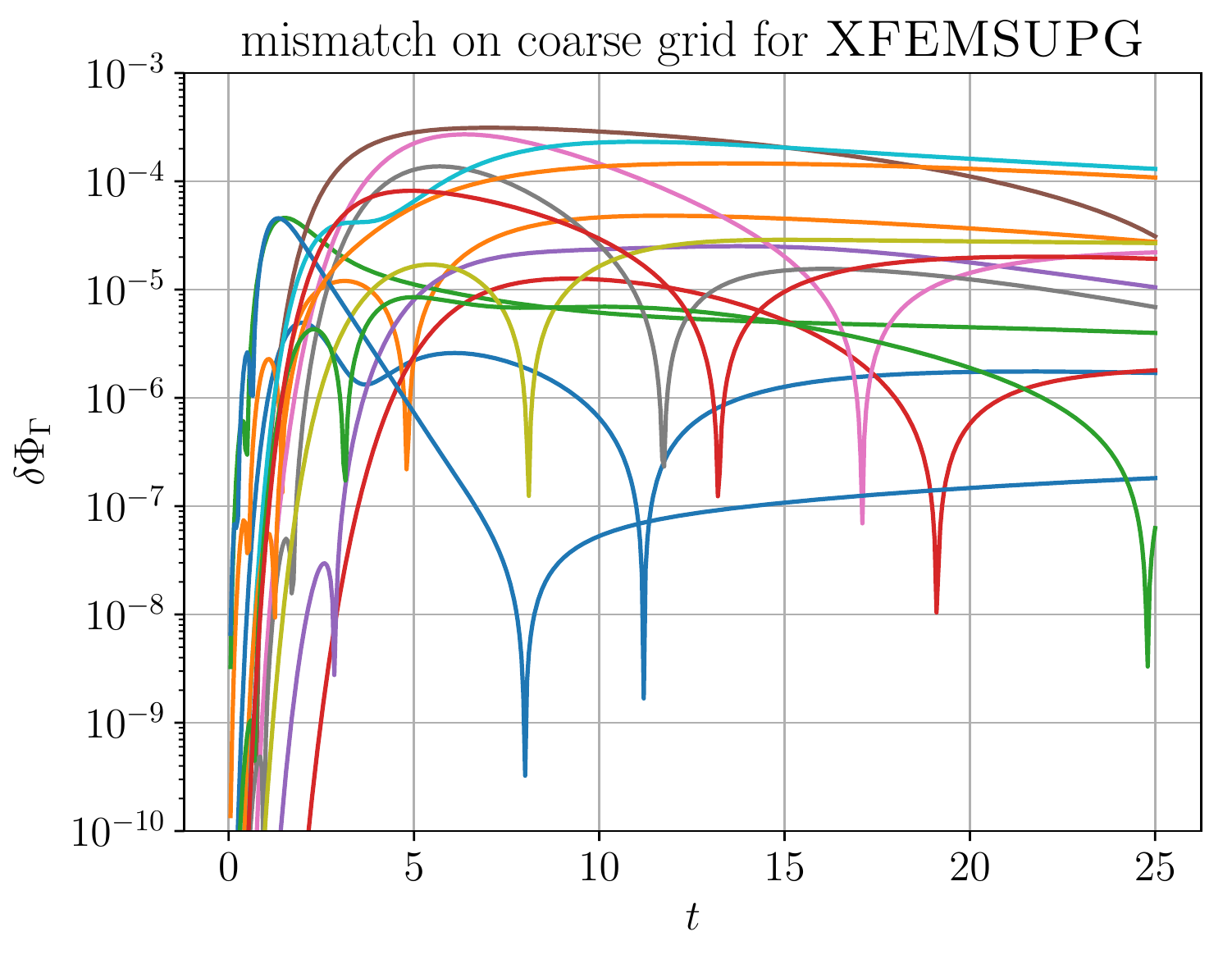}
\caption{Total flux mismatch against time for all the traces of the network of Test Case 2 on the coarse mesh with methods \FEMSUPG{} (left) and  \XFEMSUPG{} (right).}
\label{fig:example2_flux_mismatch}
\end{figure}

Curves of the total flux mismatch at the traces are reported against time in
Figure~\ref{fig:example2_flux_mismatch} for the \FEMSUPG{} and \XFEMSUPG{}
methods. In this picture values of $\delta \Phi_\Gamma$ are shown, without
labels, for all the traces in the network, highlighting that, in all cases, the
errors remain limited in time. Further, maximum-in-time mismatch values are
lower than $1\%$ of the total flux, for all the traces; values of the maximum
flux with respect to time on each trace are reported in
Table~\ref{tab:example2_flux_on_traces}, computed with the \XFEMSUPG{} method on
the finest mesh.

\begin{table}
\centering
\caption{Maximum value in time of total flux across each trace of the DFN for Test Case 2 computed with the \XFEMSUPG{} method on the fine mesh.}
\label{tab:example2_flux_on_traces}
\begin{tabular}{c|c|c|c|c|c|c}
$\Phi_{\Gamma_1}$ & $\Phi_{\Gamma_2}$ & $\Phi_{\Gamma_3}$ & $\Phi_{\Gamma_4}$ & $\Phi_{\Gamma_5}$ & $\Phi_{\Gamma_6}$ & $\Phi_{\Gamma_7}$ \\
\hline
0.040 & 0.079 & 0.150 & 0.260 & 0.284 & 0.047 & 0.074\\
\hline
\hline
$\Phi_{\Gamma_8}$ & $\Phi_{\Gamma_9}$ & $\Phi_{\Gamma_10}$ & $\Phi_{\Gamma_11}$ & $\Phi_{\Gamma_12}$ & $\Phi_{\Gamma_13}$ & $\Phi_{\Gamma_14}$ \\
\hline
0.020 & 0.393 & 0.163 & 0.443 & 0.061 & 0.058 & 0.063

\end{tabular}
\end{table}


\subsection{Extruded real outcrop}\label{subsec:example3}

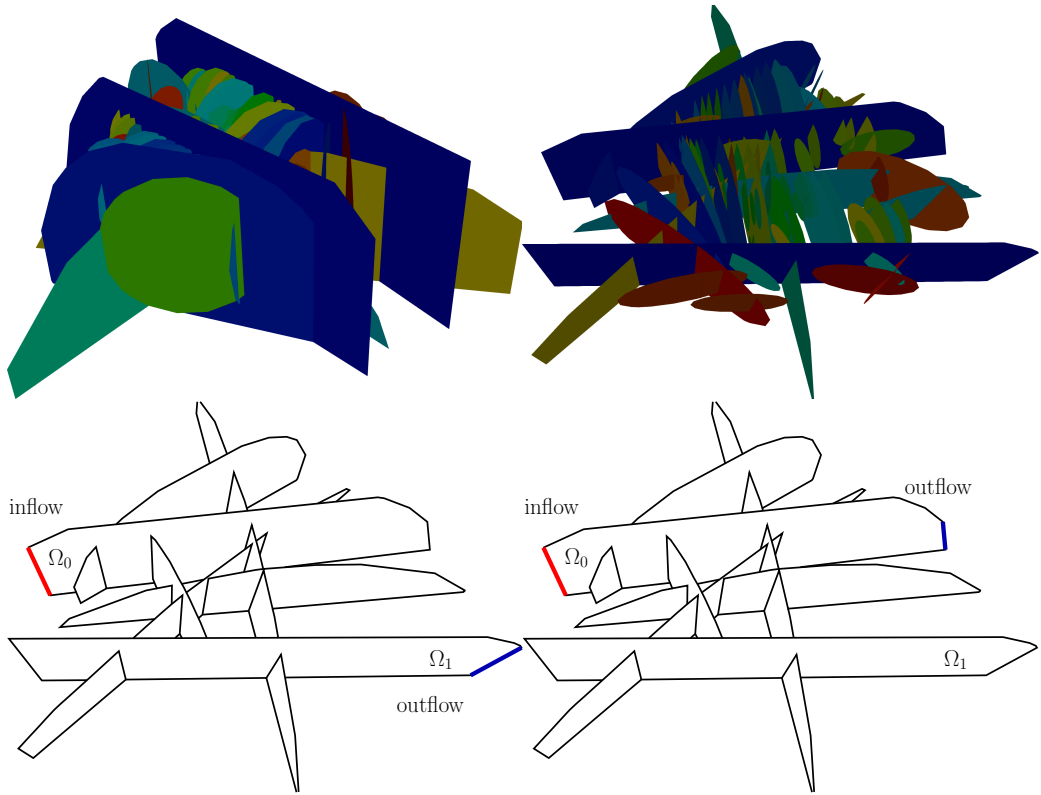
\begin{figure}[htbp]
  \centering
  \includegraphics[width=0.5\textwidth]{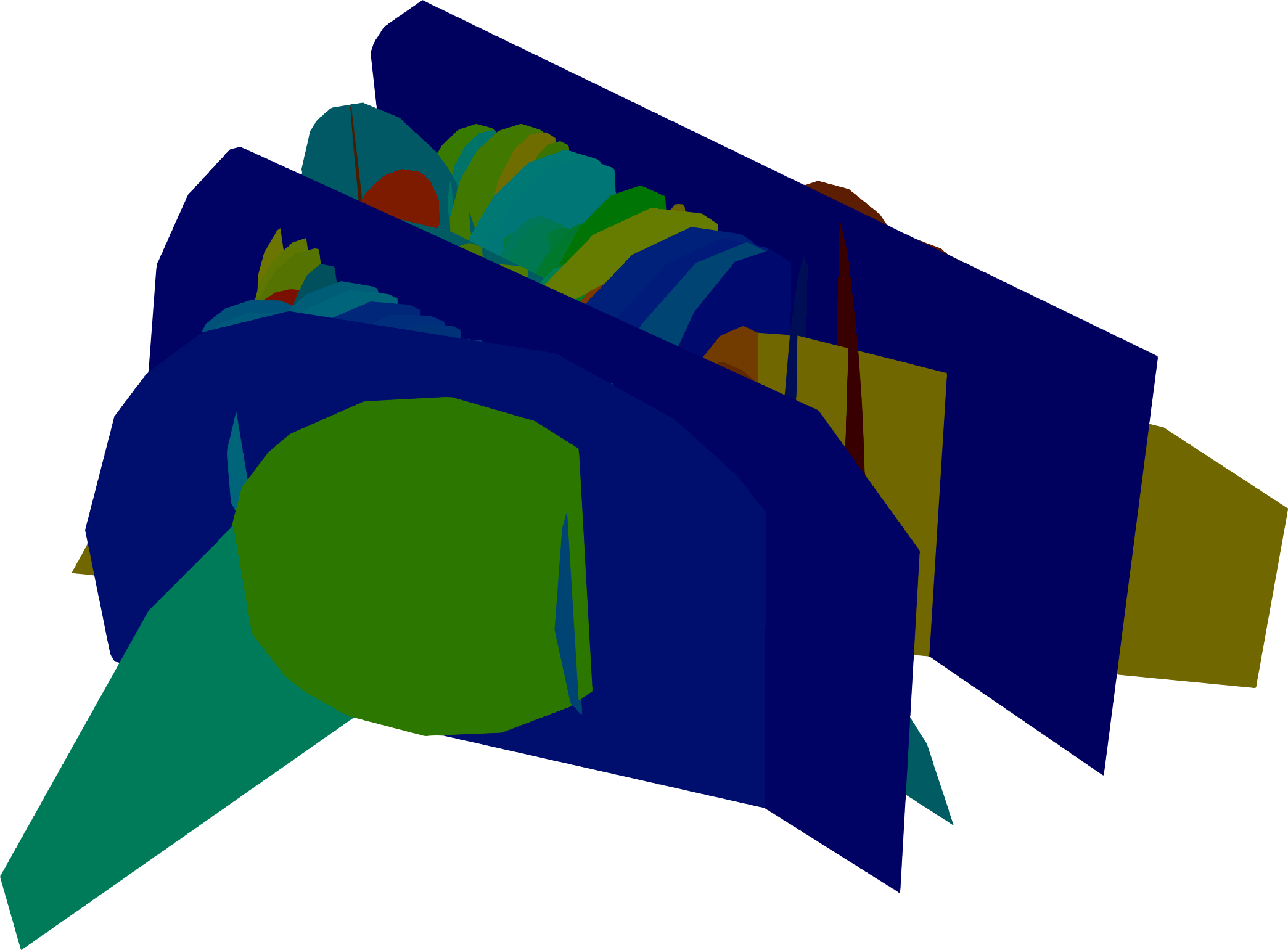}%
  \includegraphics[width=0.5\textwidth]{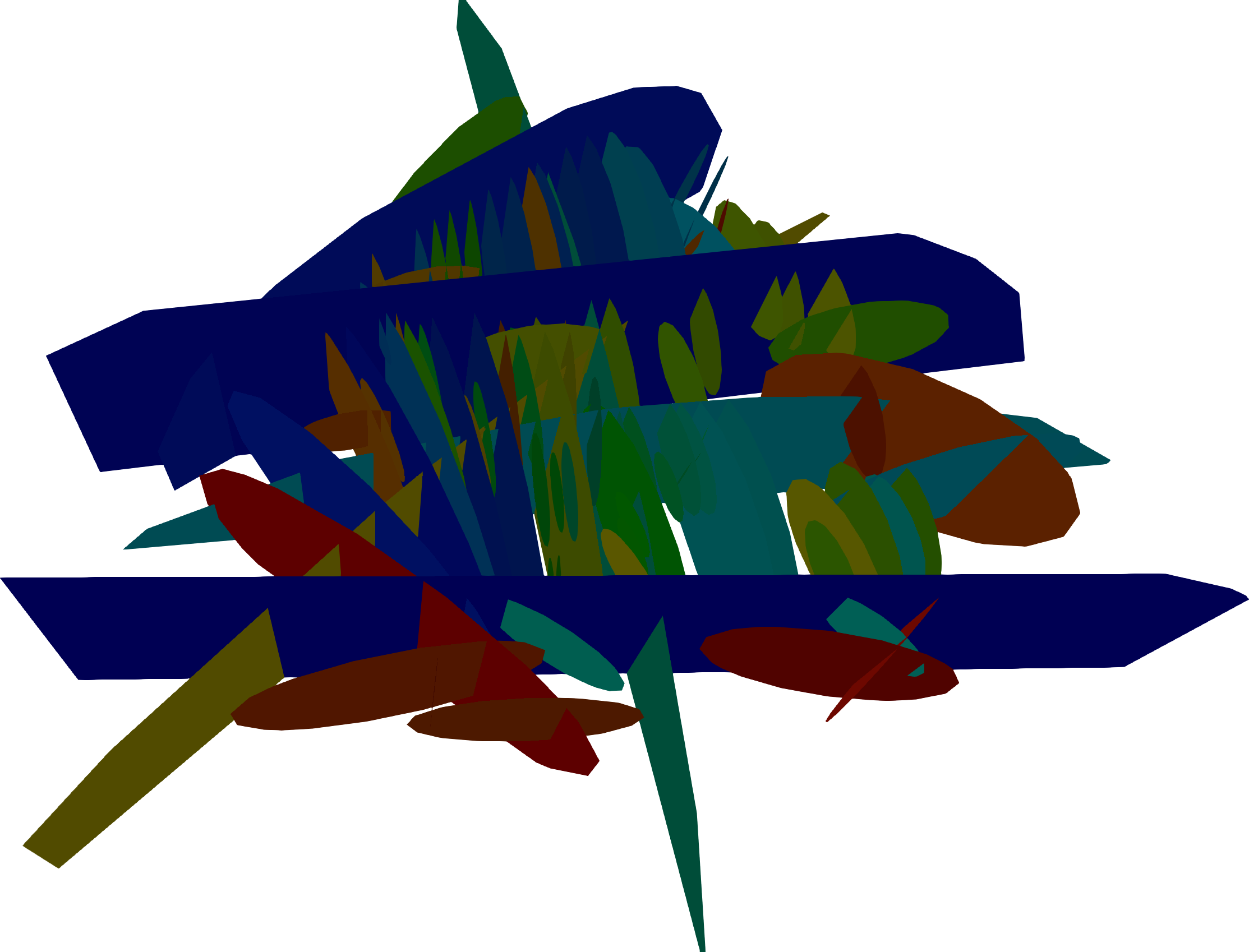}\\
  \resizebox{0.5\textwidth}{!}{\fontsize{20pt}{8}\selectfont%
\begingroup%
  \makeatletter%
  \providecommand\color[2][]{%
    \errmessage{(Inkscape) Color is used for the text in Inkscape, but the package 'color.sty' is not loaded}%
    \renewcommand\color[2][]{}%
  }%
  \providecommand\transparent[1]{%
    \errmessage{(Inkscape) Transparency is used (non-zero) for the text in Inkscape, but the package 'transparent.sty' is not loaded}%
    \renewcommand\transparent[1]{}%
  }%
  \providecommand\rotatebox[2]{#2}%
  \newcommand*\fsize{\dimexpr\f@size pt\relax}%
  \newcommand*\lineheight[1]{\fontsize{\fsize}{#1\fsize}\selectfont}%
  \ifx\svgwidth\undefined%
    \setlength{\unitlength}{463.20314422bp}%
    \ifx\svgscale\undefined%
      \relax%
    \else%
      \setlength{\unitlength}{\unitlength * \real{\svgscale}}%
    \fi%
  \else%
    \setlength{\unitlength}{\svgwidth}%
  \fi%
  \global\let\svgwidth\undefined%
  \global\let\svgscale\undefined%
  \makeatother%
  \begin{picture}(1,0.75933571)%
    \lineheight{1}%
    \setlength\tabcolsep{0pt}%
    \put(0,0){\includegraphics[width=\unitlength,page=1]{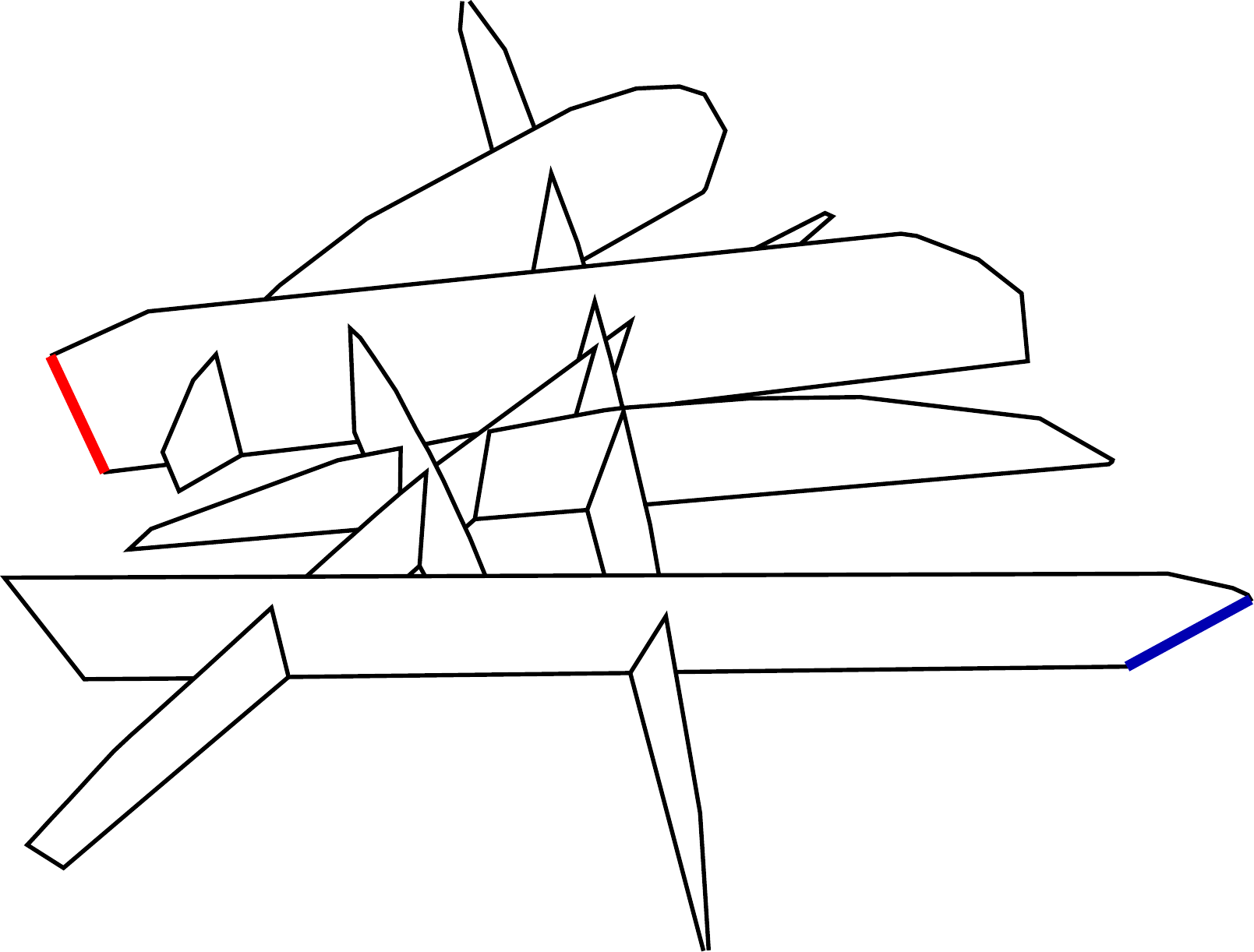}}%
    \put(0.00357184,0.54096145){\color[rgb]{0,0,0}\makebox(0,0)[lt]{\lineheight{1.25}\smash{\begin{tabular}[t]{l}inflow\end{tabular}}}}%
    \put(0.75581443,0.15595061){\color[rgb]{0,0,0}\makebox(0,0)[lt]{\lineheight{1.25}\smash{\begin{tabular}[t]{l}outflow\end{tabular}}}}%
    \put(0.08034002,0.4420005){\color[rgb]{0,0,0}\makebox(0,0)[lt]{\lineheight{1.25}\smash{\begin{tabular}[t]{l}$\Omega_0$\end{tabular}}}}%
    \put(0.81915647,0.24636488){\color[rgb]{0,0,0}\makebox(0,0)[lt]{\lineheight{1.25}\smash{\begin{tabular}[t]{l}$\Omega_1$\end{tabular}}}}%
  \end{picture}%
\endgroup%
}%
  \resizebox{0.5\textwidth}{!}{\fontsize{20pt}{8}\selectfont%
\begingroup%
  \makeatletter%
  \providecommand\color[2][]{%
    \errmessage{(Inkscape) Color is used for the text in Inkscape, but the package 'color.sty' is not loaded}%
    \renewcommand\color[2][]{}%
  }%
  \providecommand\transparent[1]{%
    \errmessage{(Inkscape) Transparency is used (non-zero) for the text in Inkscape, but the package 'transparent.sty' is not loaded}%
    \renewcommand\transparent[1]{}%
  }%
  \providecommand\rotatebox[2]{#2}%
  \newcommand*\fsize{\dimexpr\f@size pt\relax}%
  \newcommand*\lineheight[1]{\fontsize{\fsize}{#1\fsize}\selectfont}%
  \ifx\svgwidth\undefined%
    \setlength{\unitlength}{463.20314422bp}%
    \ifx\svgscale\undefined%
      \relax%
    \else%
      \setlength{\unitlength}{\unitlength * \real{\svgscale}}%
    \fi%
  \else%
    \setlength{\unitlength}{\svgwidth}%
  \fi%
  \global\let\svgwidth\undefined%
  \global\let\svgscale\undefined%
  \makeatother%
  \begin{picture}(1,0.75933571)%
    \lineheight{1}%
    \setlength\tabcolsep{0pt}%
    \put(0,0){\includegraphics[width=\unitlength,page=1]{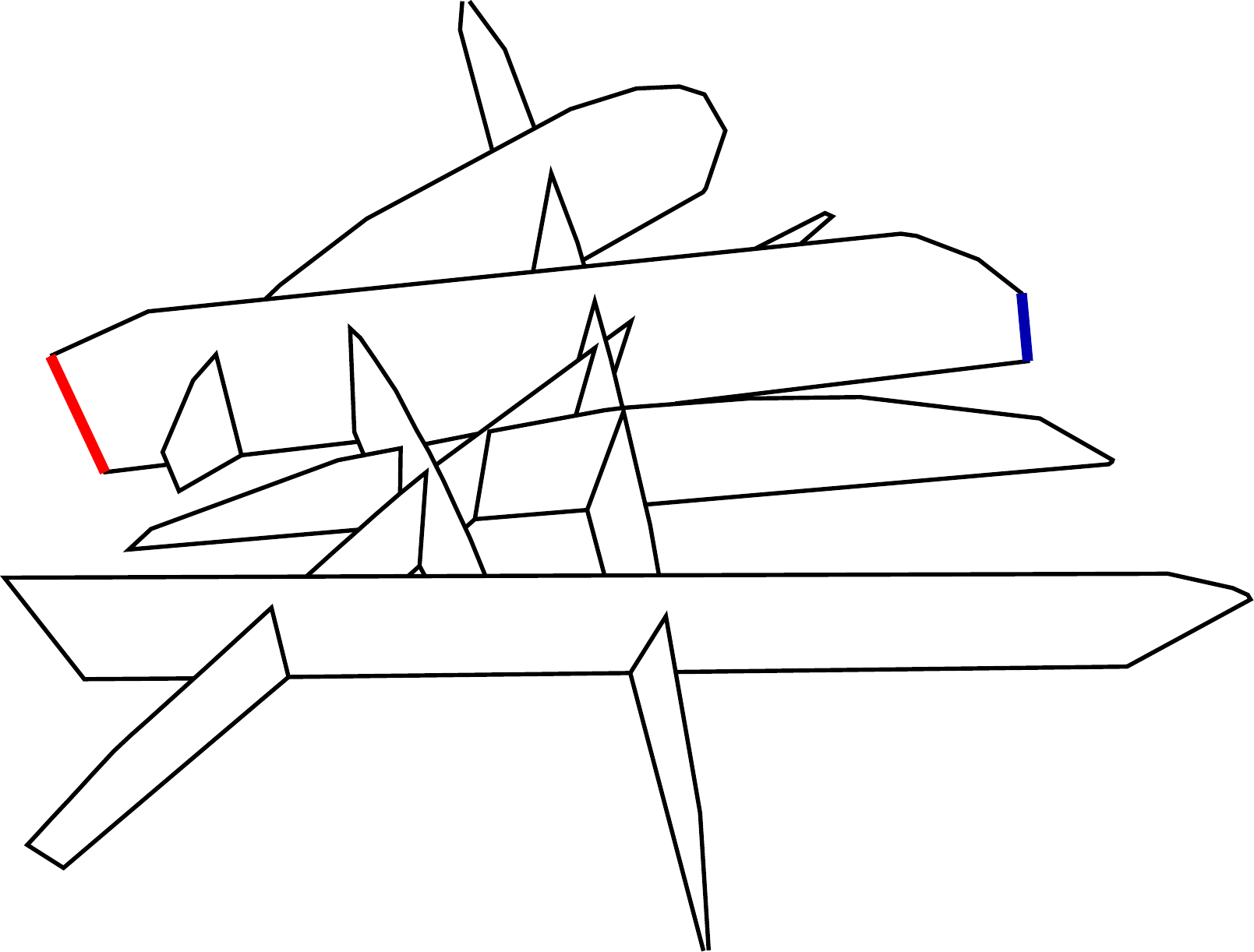}}%
    \put(0.00357184,0.54096145){\color[rgb]{0,0,0}\makebox(0,0)[lt]{\lineheight{1.25}\smash{\begin{tabular}[t]{l}inflow\end{tabular}}}}%
    \put(0.74109406,0.57957067){\color[rgb]{0,0,0}\makebox(0,0)[lt]{\lineheight{1.25}\smash{\begin{tabular}[t]{l}outflow\end{tabular}}}}%
    \put(0.08150151,0.4442725){\color[rgb]{0,0,0}\makebox(0,0)[lt]{\lineheight{1.25}\smash{\begin{tabular}[t]{l}$\Omega_0$\end{tabular}}}}%
    \put(0.81752085,0.24636497){\color[rgb]{0,0,0}\makebox(0,0)[lt]{\lineheight{1.25}\smash{\begin{tabular}[t]{l}$\Omega_1$\end{tabular}}}}%
  \end{picture}%
\endgroup%
}
  \caption{Geometry of the test case in Subsection
    \ref{subsec:example3}. On the top two screenshot of the geometry
    of the network. On the bottom a sketch of the network with only
    the largest fractures to present the boundary conditions: on the
    left in the case of inflow and outflow imposed on two different
    fractures, on the right on the same fracture. The inflow is
    represented in red and the outflow in blue.}%
  \label{fig:geometry_example_3}
\end{figure}

In this last test case we consider a fracture network generated from
an extruded outcrop, located in Western Norway. The test case is
inspired by Section 4.4 of \cite{Scialo2017}. The network is composed
by 89 intersecting fractures resulting in 166 traces. There are 7
non-connected fractures and with no flow boundary conditions, which will
not contribute to the solution. The geometry is depicted in Figure
\ref{fig:geometry_example_3}.  The aim of this test case is to
validate the proposed numerical schemes in presence of realistic
physical parameters of a real fracture network. However, we assume
that all the fractures share the same values of hydraulic conductivity
and heat diffusion coefficient.

We consider two distinct problems, that differ from each other from
the position of the inflow and outflow boundaries. In the first case
(denoted as \textit{case 1}) the inflow and outflow are imposed on two
different fractures, while in the second case (\textit{case 2}) they
belong to the same fracture. See Figure \ref{fig:geometry_example_3}
where a sketch of the network is shown with the position of the inflow
and outflow boundaries. On the other portions of the boundaries a no
flow boundary condition is given. In both cases, we require a
simulation grid with roughly 70k elements.

Fractures are immersed in granite and we assume that at the beginning of the simulation the water contained
in the fractures is at \SI{353.15}{\kelvin} (\SI{80}{\celsius}). The relations to
compute the physical parameters for the simulations are the ones
presented in Subsection \ref{subsec:heat}. We assume
$\epsilon_i = \SI{2}{\milli\metre}$ $\forall i$, and $\phi_i = 0.95$
$\forall i$. The water and rock physical parameters are
reported in Table \ref{tab:values}.
\begin{table}
  \centering
  \begin{tabular}{l|cc}
    \hline
    & Water & Rock
    \\
    \hline
    Dynamic viscosity & $\mu = \SI{3.55}{\pascal\second}$ &  --
    \\
    Thermal conductivity & $\lambda_w =
                           \SI{0.667}{\watt\per\metre\per\kelvin}$
            & $\lambda_m = \SI{3.07}{\watt\per\metre\per\kelvin}$
    \\
    Density & $\rho_w = \SI{1000}{\kilo\gram\per\cubic\metre}$
            & $\rho_m = \SI{2700}{\kilo\gram\per\cubic\metre}$
    \\
    Specific heat capacity & $c_w = \SI{4099}{\joule\per\kilo\gram\per\kelvin}$
            & $c_m = \SI{790}{\joule\per\kilo\gram\per\kelvin}$
    \\
    Heat transfer coefficient & \multicolumn{2}{c}{$\gamma =
                                \SI{1.25e-3}{\watt\per\square\meter\per\kelvin}$}
    \\
    \hline
  \end{tabular}
  \caption{List of rock and water coefficients for the example in
    Subsection \ref{subsec:example3}.}%
  \label{tab:values}
\end{table}
From these data we obtain, $\forall i = 1,\ldots, N_\Omega$,
\begin{align*}
  K_i &\approx \SI{1.84e-6}{\square\metre\per\second} \,,
  &  c_{e,i} &\approx \SI{4000700}{\joule\per\cubic\metre\per\kelvin} \,,
  \\
  \lambda_{e,i} &\approx \SI{0.72}{\watt\per\metre\per\kelvin} \,,
  & \zeta_i &\approx \SI{1.95e-3}{\metre} \,,
  \\
  D_i &\approx \SI{0.35e-9}{\cubic\metre\per\second} \,,
  & \iota_i &\approx \SI{3.05e-10}{\metre\per\second} \,.
\end{align*}

Regarding boundary conditions, for the Darcy problem we impose a
pressure head equal to $\SI{2.5}{\kilo\metre}$ at the inflow boundary
and $\SI{0}{\metre}$ at the outflow boundary, while for the heat
problem we impose \SI{303.15}{\kelvin} (\SI{30}{\celsius}) at the
inflow and zero diffusive flux at the outflow. The simulation time is
a year (\SI{3.154e+7}{\second}), divided in 200 time steps.

\begin{figure}
\includegraphics[width=0.99\textwidth]{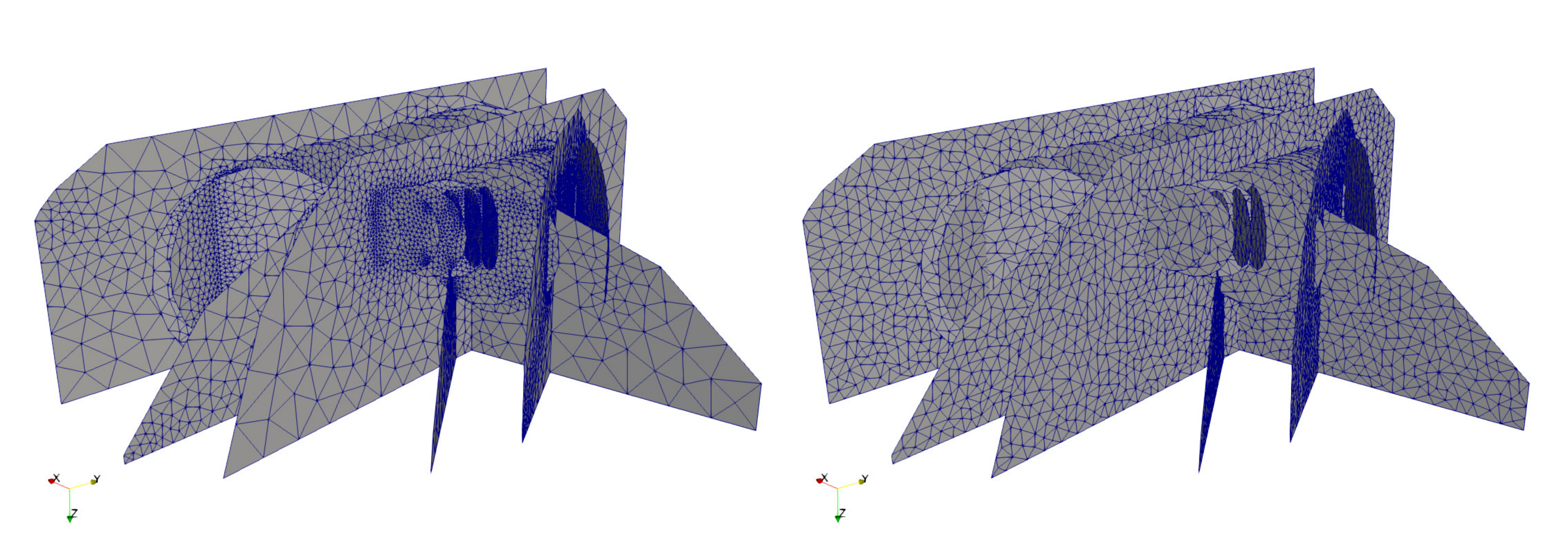}
\caption{Conforming (left) and non-matching (right) mesh for Test Case 3}
\label{fig:example3_meshes}
\end{figure}

The conforming computational mesh counts about $7\times 10^4$ elements, while
the non-matching computational mesh has $2\times 10^4$ cells and are shown in
Figure~\ref{fig:example3_meshes}. It is possible to notice, how, in order to
meet the conformity requirement, mesh elements of the conforming mesh are
concentrated near the traces, whereas, the non-matching mesh has all elements of
equal size evenly distributed in the network. The non matching mesh is
characterized by a mesh P\'eclet number of about $3\times 10^4$ for \textit{case
1} and $6.6\times10^4$ for \textit{case 2}, reached on fracture $\Omega_0$, in
both cases. The non-adapted mesh and the high mesh P\'eclet number make this
example extremely complex for methods built on non-conforming meshes and relying
on stabilization for advection dominated flow regimes.

\begin{figure}[htbp]
  \centering
  \includegraphics[width=0.48\textwidth]{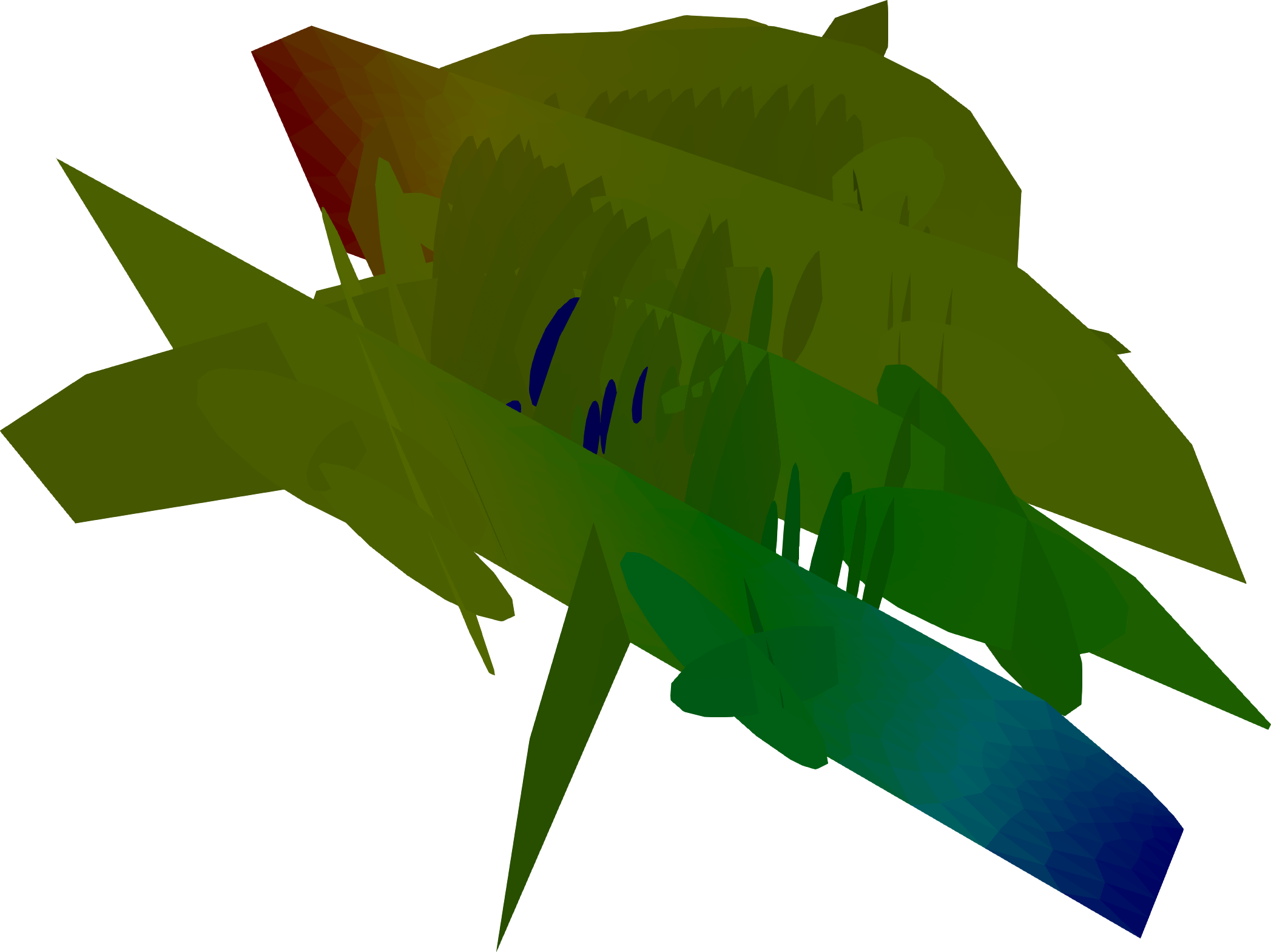}%
  \includegraphics[width=0.48\textwidth]{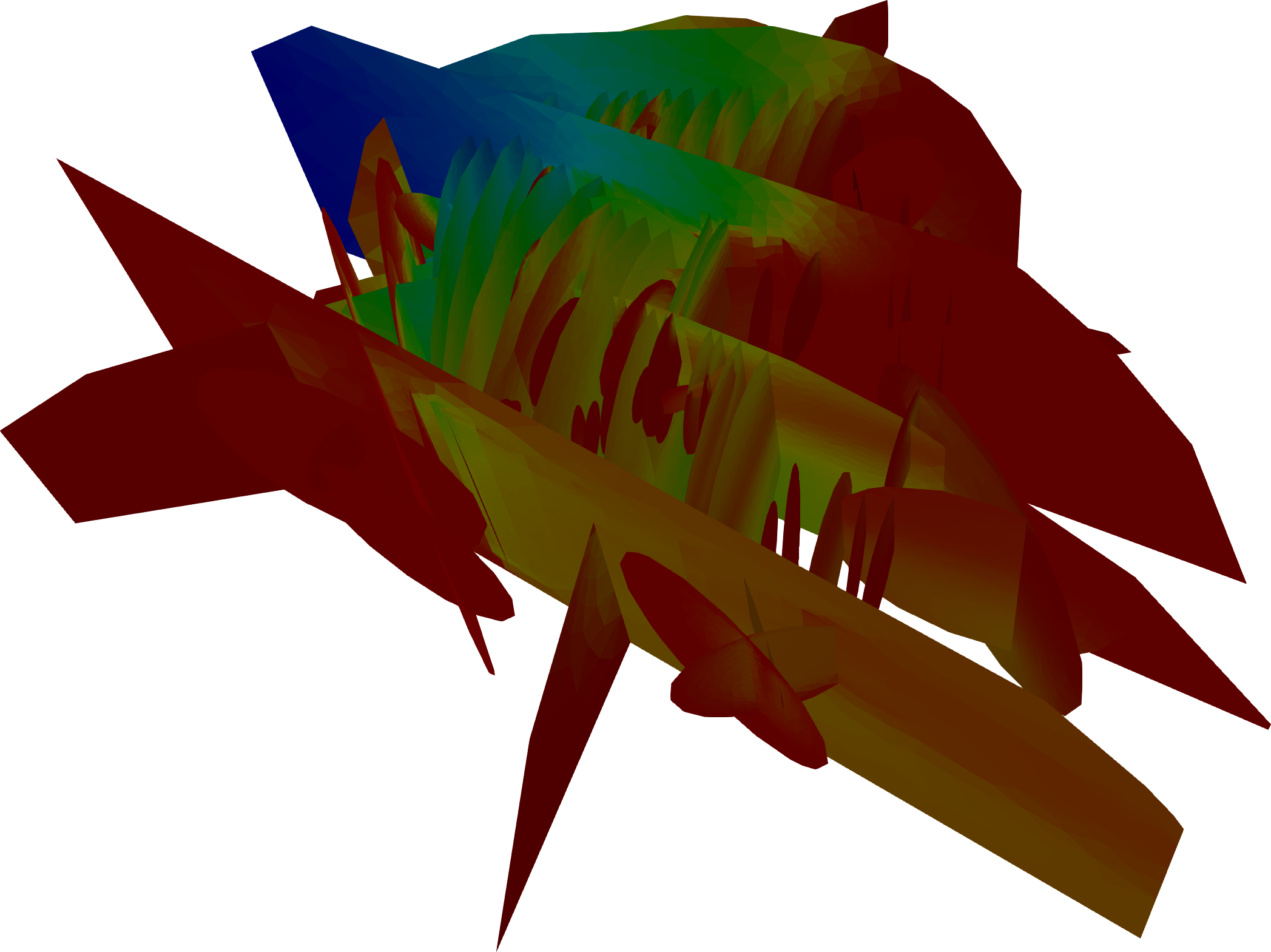}\\
  \includegraphics[width=0.48\textwidth]{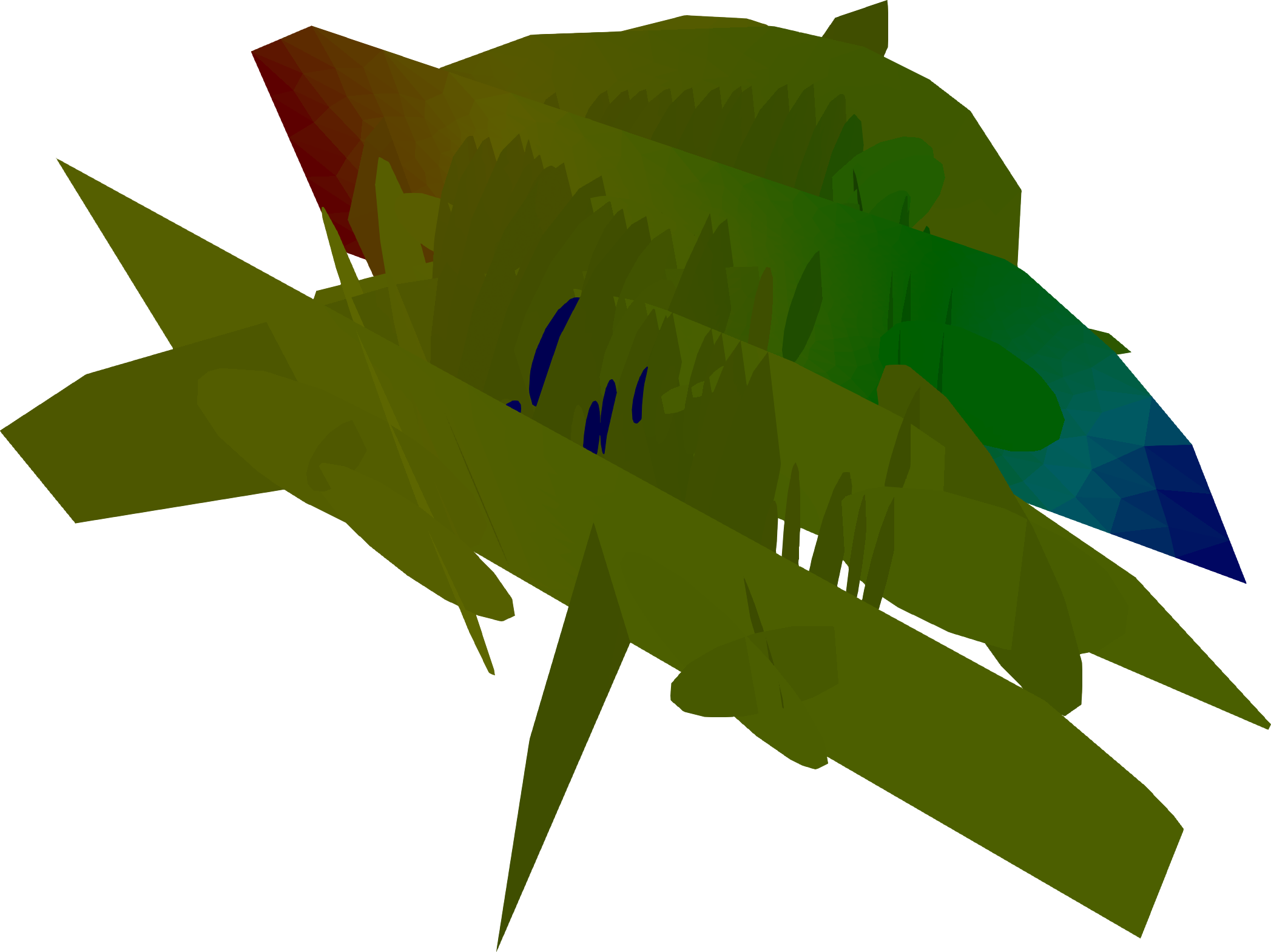}
  \includegraphics[width=0.48\textwidth]{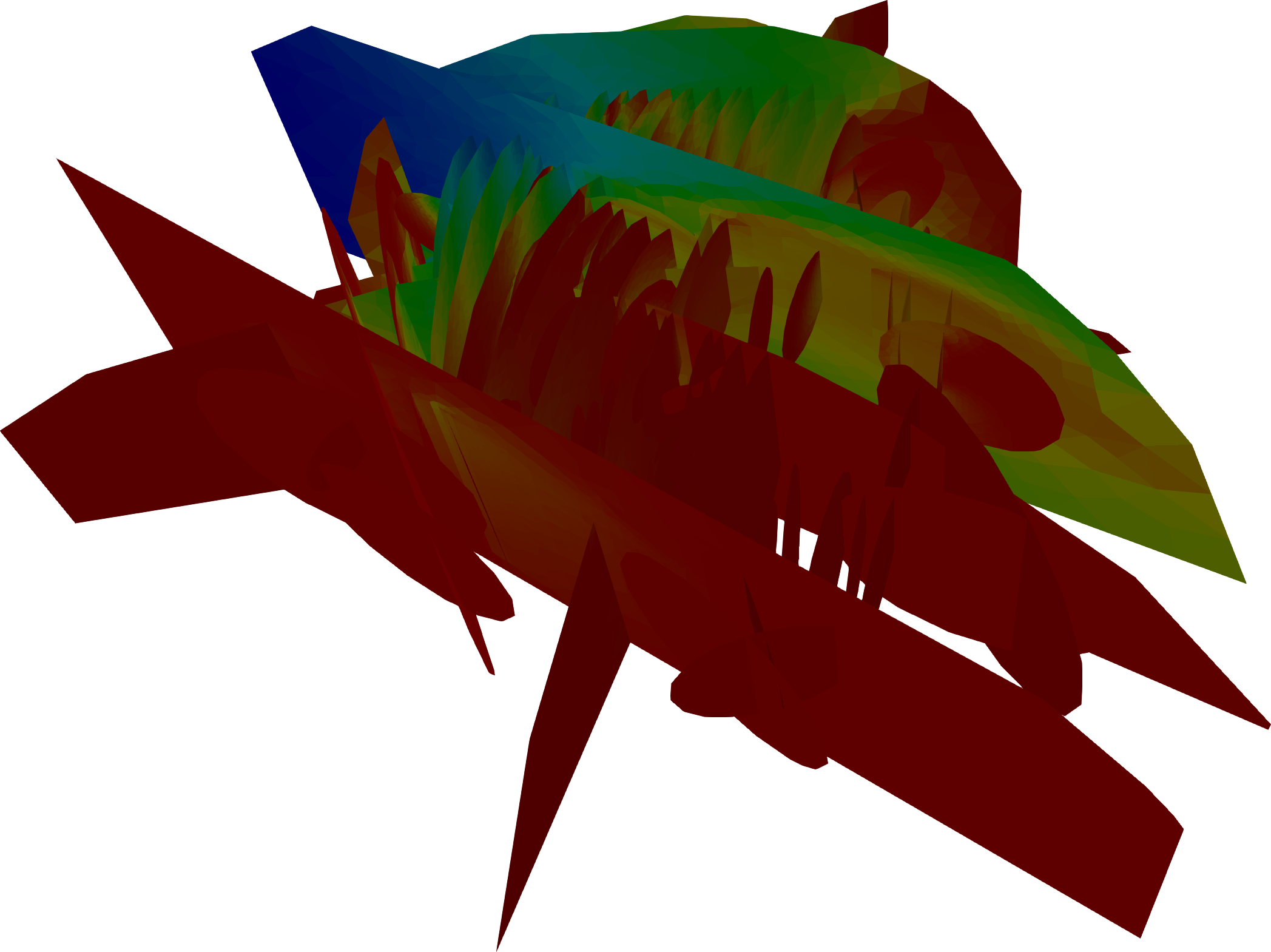}\\
  \caption{Solution with the \MFEMUP{} scheme for the Test Case 3: \textit{case
  1} on the top; \textit{case 2} at the bottom. The first column shows the
  pressure head solution, the second column the temperature distribution at the end of time-evolution.}%
  \label{fig:solution_example_3}
\end{figure}

The computed solution with the \MFEMUP{} scheme is reported in Figure~\ref{fig:solution_example_3} where we display the solution of the Darcy problem on the left and the solution at the end of time evolution on the right, for both the setting of \textit{case 1}, on top and \textit{case 2}, at the bottom.

\begin{figure}
\includegraphics[width=0.32\textwidth]{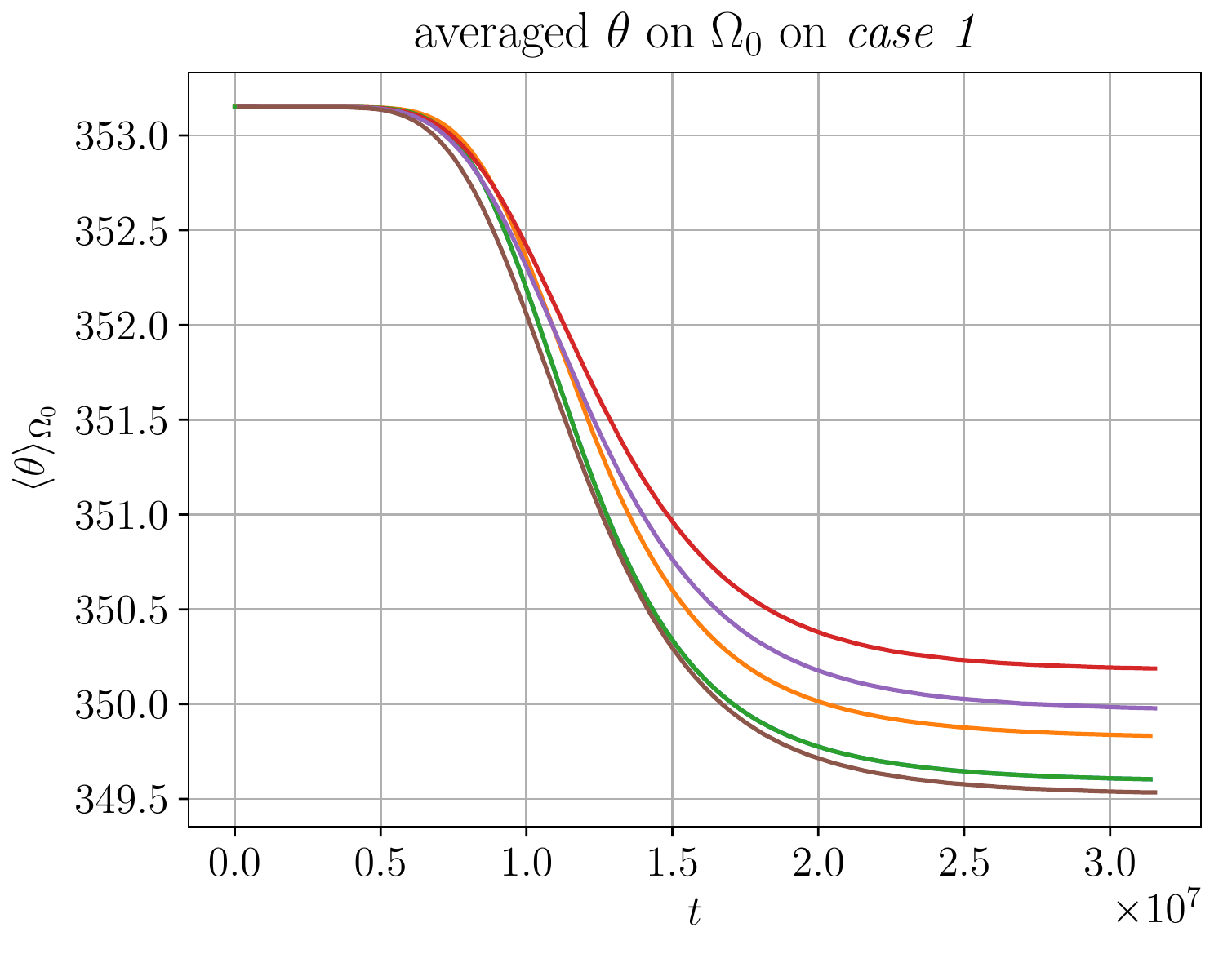}
\includegraphics[width=0.32\textwidth]{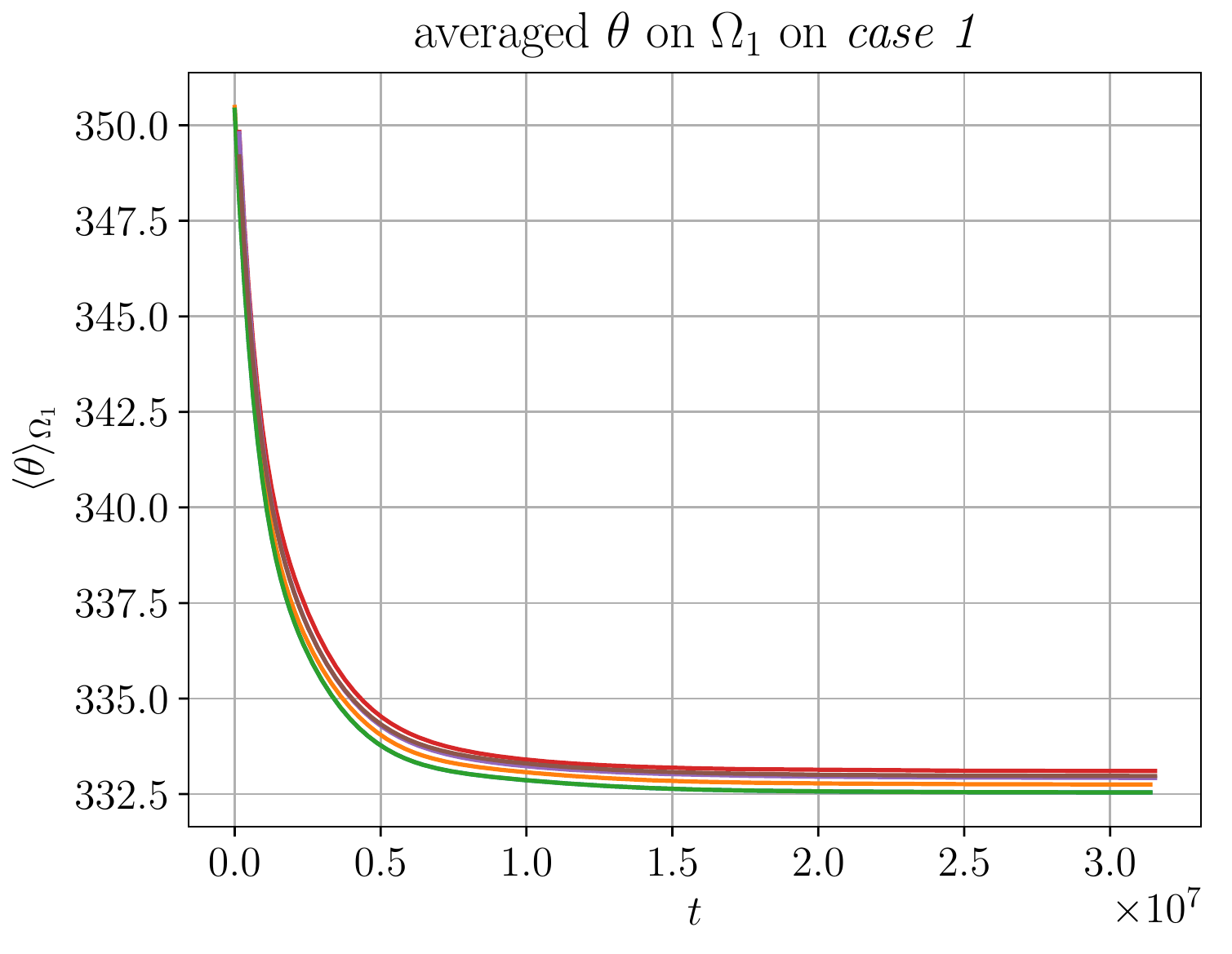}
\includegraphics[width=0.32\textwidth]{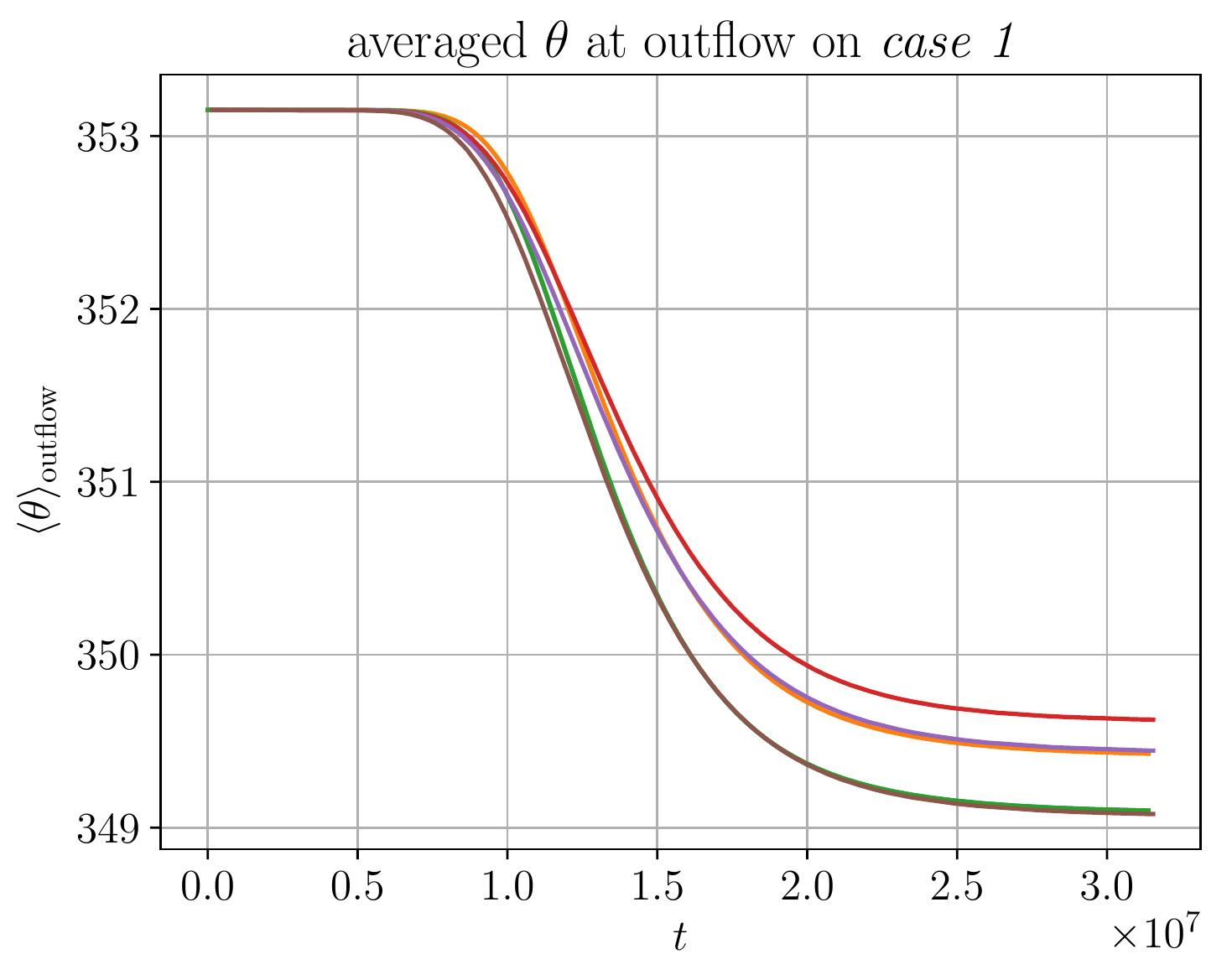}\\
\includegraphics[width=0.32\textwidth]{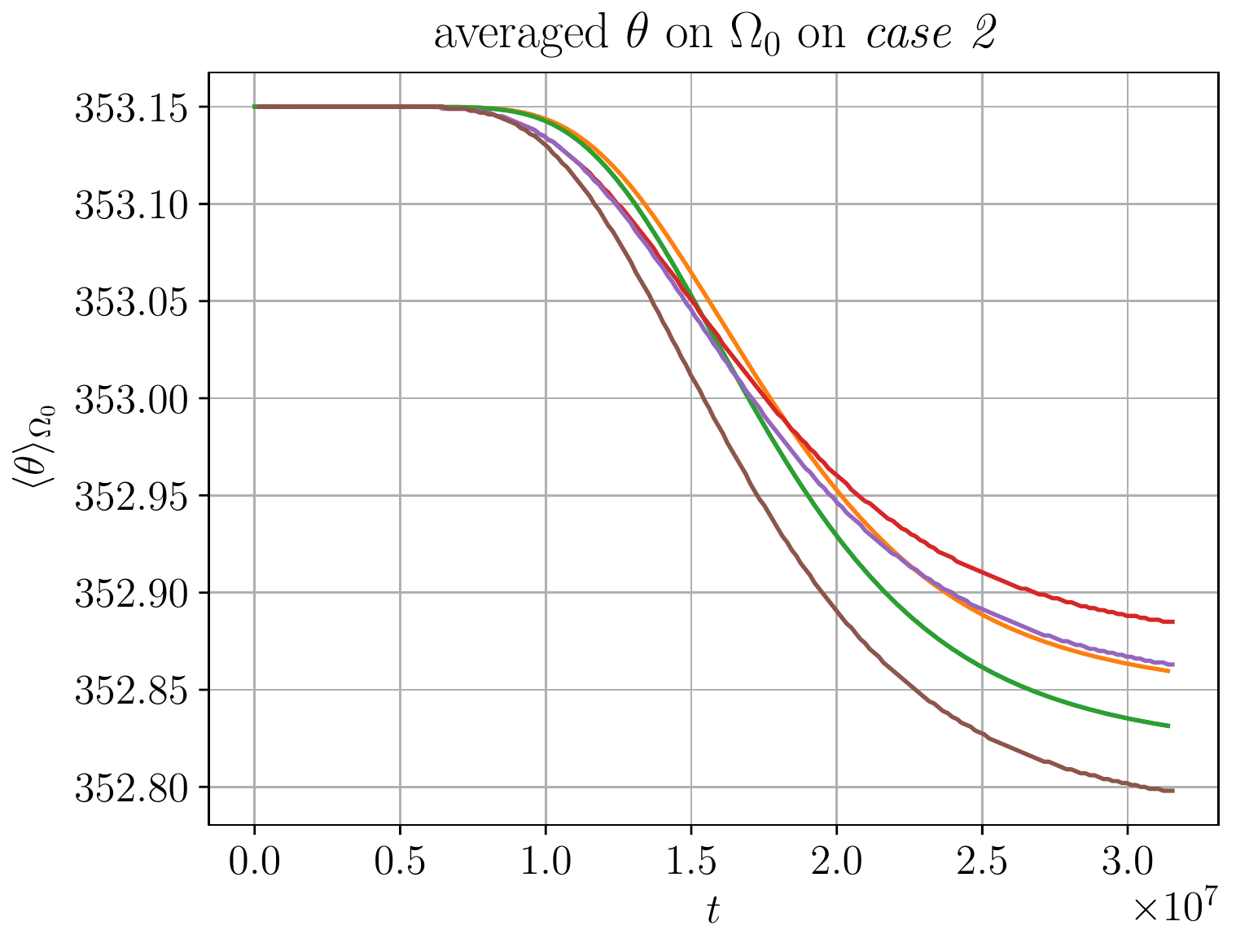}
\includegraphics[width=0.32\textwidth]{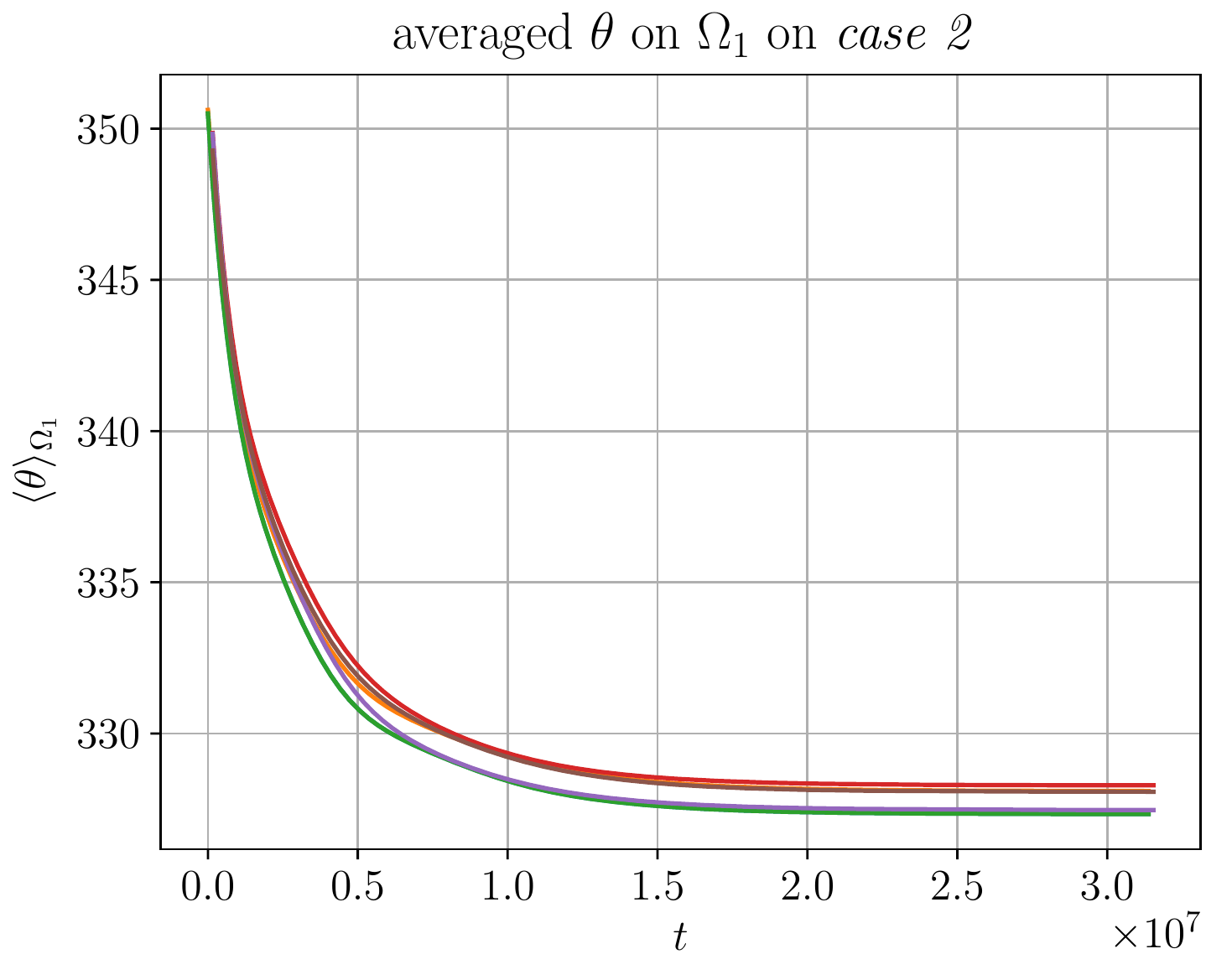}
\includegraphics[width=0.32\textwidth]{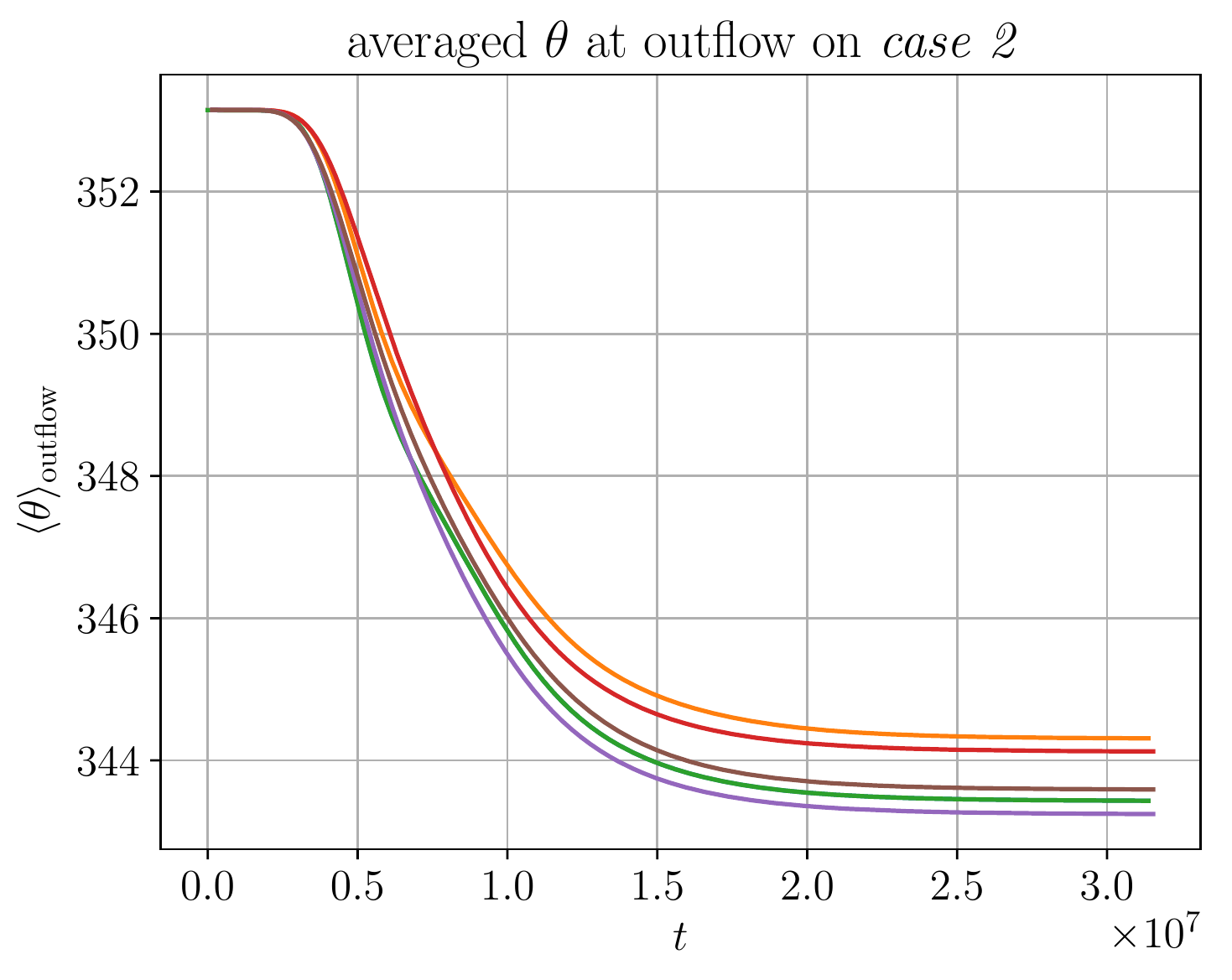}\\
\includegraphics[width=\textwidth]{label}
\caption{Curves of average theta on fractures $\Omega_0$ (left) and $\Omega_1$
and at outflow (right) against time for Test Case 3: \textit{case 1} top, \textit{case 2} bottom.}
\label{fig:example3_avg_prod}
\end{figure}

Figure~\ref{fig:example3_avg_prod} shows the curves against time of the average
temperature on the inflow ($\Omega_0$) and outflow ($\Omega_1$) fractures, along
with the average temperature on the outflow boundary. Despite the complexity of
the geometry and of the model, curves appear in good agreement. For this last
example no coarsening was used for the \MVEMUP{} method, as the poor
performances of the \TPFA{} method on polygonal cells, already observed in
Test Case 2, have a strong impact on the quality of the solution in this more
complex case. The curves related to the \MVEMUP{} approach on triangular meshes
are almost perfectly overlapped to the curves of the \MFEMUP{} method. Curves of
\XFEMSUPG{} and \FEMSUPG{} are in good agreement with those of the other
methods.

\section{Conclusions}\label{sec:conclusion}

In this work, we presented a detailed comparative study of several solution strategies for single-phase flow
and transport in discrete fracture networks. The proposed numerical schemes are challenged with networks of increasing geometrical complexity and with unsteady advection-reaction-diffusion problems. The characteristics of the various approaches are compared in terms of flexibility in handling geometrical complexity, local and global conservativity and stability to high P\'eclet numbers.

Methods based on matching grids at the traces trade simplicity in imposing coupling conditions with the lack of control on the number of mesh elements, which is actually constrained by the conformity requirement. On the other hand, non-matching and polygonal based approaches demand ad-hoc discretization strategies but allow full flexibility in meshing. Non-locally conservative schemes might be susceptible to loss in mass in particular at the intersection of fractures. The proposed examples showed that this quantity is small and, for many practical problems it is possible to conclude that it might be acceptable compared to the usual uncertainty in the model parameters. Finally, we have compared schemes that are naturally stable with respect to high P\'eclet number, with schemes that require a stabilization term to avoid spurious oscillations. Also in this case the obtained solutions are coherent with respect to each other. We can conclude that numerical schemes based on polygonal or non-matching meshes give a good balance in terms of computational cost and accuracy with less limitations compared to the matching cases on more standard grids.



\section*{Acknowledgments}

This research has been partially supported by the MIUR project
``Dipartimenti di Eccellenza 2018-2022'' (CUP E11G18000350001), PRIN project "Virtual Element Methods: Analysis and Applications" (201744KLJL\_004) and by
INdAM - GNCS.  The second author acknowledges financial support for the
ANIGMA project from the Research Council of Norway (project
no. 244129/E20) through the ENERGIX program and from the ``Visiting
Professor Project'' of Politecnico di Torino.  The first and
third author acknowledge financial support from Politecnico di Torino
through project ``Starting grant RTD''.  The authors warmly thank
Luisa F. Zuluaga for constructing and providing the real fracture
network for the example in Subsection \ref{subsec:example3}.  Finally,
the authors warmly thank the PorePy development team for many fruitful
discussions related to the development of this work.









\section*{References}

\providecommand{\noopsort}[1]{}\providecommand{\noopsort}[1]{}

\end{document}